\documentclass[twoside,11pt]{article}
\usepackage{amsmath, mathtools}
\usepackage{smile}
\usepackage[linesnumbered,ruled,vlined]{algorithm2e}
\usepackage{tikz}
\usepackage{bbm}
\usetikzlibrary{shapes}
\usetikzlibrary{plotmarks}
\usepackage{xcolor}
\usepackage{capt-of}
\usepackage{fullpage}
\usepackage[colorlinks=true,
            linkcolor=blue,
            urlcolor=blue,
            citecolor=blue]{hyperref}

\usepackage[colorinlistoftodos]{todonotes}

\SetKwInput{KwInput}{Input}               
\SetKwInput{KwOutput}{Output}            
\SetKwRepeat{KwRepeat}{repeat}{until}

\theoremstyle{definition}
\newtheorem*{definition*}{Definition}

\newtheorem*{remark*}{Remark}

\begin{document}
\begin{center}

{\LARGE Minimaxity and Efficiency in Exponential Family Regression: From Star-Shaped to Convex Constraints}
\vspace*{1.5em}

{
    \renewcommand{\thefootnote}{\fnsymbol{footnote}}
    \large
    
    Yikun Li\textsuperscript{*},\hspace*{1em} Guanghong Yi\textsuperscript{*},\hspace*{1em} Matey Neykov
}
\vspace*{1.5em}

\begingroup
  \renewcommand{\thefootnote}{\fnsymbol{footnote}}
  \footnotetext[1]{These authors contributed equally to this work and are listed in alphabetical order.}
\endgroup

Department of Statistics and Data Science, Northwestern University
\vspace*{1em}

\texttt{\{yikunli2028, guanghongyi2025\}@u.northwestern.edu ~  mneykov@northwestern.edu}

\end{center}

\begin{abstract}

This paper establishes the minimax estimation rate for nonparametric exponential family regression under star-shaped constraints. We consider a parameter space $K$ that is a star-shaped subset of the hypercube $[-M, M]^n$ for a known constant $M > 0$. We operate under the assumption that the underlying exponential family is nonsingular with a twice continuously differentiable log-partition function. Our main result demonstrates that the minimax rate of the $\ell _{2}$ error of such estimation problem is $\epsilon^{*2} \wedge \operatorname{diam}(K)^2$ up to constants exclusively depending on $M$. Here, the critical radius $\epsilon^*$ is defined as
\begin{equation*}
    \epsilon^* = \sup \{\epsilon \left\lvert\right. \epsilon^2 \kappa(M) \le \log N^{\text{loc}}(\epsilon,c)\},
\end{equation*}
where $N^{\text{loc}}(\epsilon,c)$ denotes the local metric entropy of $K$, and $\kappa(M) > 0, c>0$ are constants depending only on $M$. Such minimax rate is established by a match between an information-theoretic lower bound and an upper bound implied by a theoretical algorithm.

Furthermore, we investigate the computational aspects of this estimation problem. Under mildly stronger assumptions on the constraint set $K$, we propose a computationally efficient, polynomial-time algorithm. We prove that the resulting estimator achieves the minimax optimal rate up to poly-logarithmic factors in the dimension $n$ and the geometric parameters of $K$.

Finally, to illustrate the efficacy of our framework, we derive the minimax optimal rates for some concrete examples.
\end{abstract}

\tableofcontents

\newpage

\section{Introduction}\label{section: introduction}

\subsection{Problem Formulation}\label{subsection: problem formulation}

We consider nonparametric exponential family regression under star-shaped constraints. Assume that we observe independent random variables $Y_1, \ldots, Y_n$, where each $Y_i$ follows a distribution with a density function of the form
\begin{equation}\label{eq: exponential family}
    f(y_i;\theta _{i}) = h(y_i) \exp\left\{\theta _{i} T(y_i) - A(\theta _{i})\right\}, i=1,2,\dots,n.
\end{equation}
In \eqref{eq: exponential family}, $\theta =(\theta _{1},\dots,\theta _{n})^{\top } $ is the natural parameter and $T(Y _{i}), i=1,\dots,n$ are the sufficient statistics. Let $\theta_0 \in \mathbb{R}^n$ denote the unknown true natural parameter vector and $\mathcal{D}:=\operatorname{domain}(A)=\left\{\theta \left\lvert\right. \theta \in \mathbb{R}, \displaystyle\int _{\mathbb{R}}^{}h(x)\exp\left\{\theta T(y)\right\}dy<\infty \right\}$ is the domain of the log-partition function $A(\cdot )$.

We note that the representation in \eqref{eq: exponential family} directly implies a fundamental relationship between $\mathbb{E}[T(Y)]$ and $\operatorname{Var}(T(Y))$. For the more general exponential dispersion model (EDM) of the form
\begin{equation*}
    f(y_i; \theta_i, \phi_i) = h(y_i, \phi_i) \exp\left\{ \frac{T(y_i)\theta_i - b(\theta_i)}{a(\phi_i)} \right\}, i=1,2,\dots,n,
\end{equation*}
we assume the dispersion parameters $\phi_i$ are known. Consequently, by mapping $T(y_i) \mapsto T(y_i)/a(\phi_i)$ throughout the paper, we can absorb the dispersion factor and maintain the canonical exponential family structure without loss of generality.

We impose some regularity conditions on the model as follows. We focus on canonical and nonsingular exponential families, meaning that $\mathcal{D}$ is an open set, and $\text{Var}_{\theta _{i}}(T(Y_i)) > 0$ for all $i \in [n]$ and any $\theta _{i}\in \mathcal{D}$. We restrict the true natural parameter $\theta_0$ to a star-shaped set $K$ (see Definition \ref{def_star_shaped_set}). Such $K$ is assumed to be contained within a bounded hyperrectangle $[M _{\text{l}}, M _{\text{u}}]^n$ with known and fixed constants $M _{\text{l}},M _{u}$. For simplicity of the notations and analysis, we assume that $M _{\text{u}}=-M _{\text{l}}=M>0$ with known $M$, and the compact set $[-M,M]$ is contained within $\mathcal{D}$. By the assumption that $\mathcal{D}$ is a open set, there exists a $\delta >0$ such that $[-(1+\delta )M,(1+\delta )M]\subset \mathcal{D}$. We assume $\delta =1$ without loss of generality. One can easily adapt the arguments and proofs thereafter to the case when there is no $M>0$ satisfying $[-M,M]\subset \mathcal{D}$ (for example, when $\mathcal{D}=\mathbb{R}^{+}$). Finally, $A(\cdot)$ is assumed to be twice continuously differentiable, a standard property satisfied by most common exponential family distributions.

Our primary objective is to characterize the minimax rate of $\ell _{2}$ estimation error for $\theta_0$, exact up to constants depending only on $M$. As a central concept in statistical decision theory, the minimax framework captures the best possible worst-case performance of an estimator over a given parameter space. Formally, the minimax risk over the constraint set $K$ is defined as
\begin{equation}\label{eq: minimax risk of the estimation}
    \mathcal{R}_{n}(f,K):=\inf _{\hat{\theta}} \sup_{\theta _{0}\in K} \mathbb{E}_{Y _{1},\dots,Y _{n}}\left\lVert \hat{\theta }(Y)-\theta _{0}\right\rVert _{2}^{2},
\end{equation}
where $\left\lVert \cdot \right\rVert _{2}^{}$ denotes the Euclidean norm, $Y=\left(Y _{1},Y _{2},\dots,Y _{n}\right)^{\top }$, and the infimum is taken over all possible estimators, i.e., all measurable functions $\hat{\theta }$ of the observed data $Y$.

To motivate our study, consider a model with a monotonic shape constraint. Suppose we observe independent Bernoulli random variables $Y_i \sim \operatorname{Bernoulli}(p_i)$ for $i \in [n]$, where the probabilities are monotonically increasing, i.e., $p_1 \le p_2 \le \ldots \le p_n$. We are primarily interested in estimating the natural parameters $\theta_i = \log \left(\frac{p_i}{1-p_i}\right)$. We assume these parameters are bounded such that $\max_{i \in [n]} |\theta_i| \le M$ for a known positive constant $M$. This condition is equivalent to restricting $\alpha \le p_1 \le \ldots \le p_n \le 1 - \alpha$ for some fixed $\alpha > 0$. Similar shape-constrained frameworks have been utilized in previous literature, such as \cite{isotonic_general_dimensions} and \cite{prasadan2024}, which investigated convex, box, and isotonic constraints within the Gaussian sequence model setting. 

The Bernoulli model described above is a specific instance of the broader problem addressed in this paper. To tackle this general setting, we introduce a novel tree-based algorithm that achieves the minimax optimal rate under the regularity assumptions described above. Our approach builds upon the framework in \cite{neykov2023minimax}, which iteratively selects a vector minimizing the $\ell_2$ norm within designed packing sets. In contrast, our algorithm introduces a new selection criterion driven by the data likelihood, coupled with a crucial pruning step inspired by \cite{10.1214/25-AOS2576}. We prove that our algorithmic upper bound matches the information-theoretic lower bound, yielding a minimax rate of $\epsilon^{*2} \wedge \text{diam}(K) ^2$ up to constants depending only on $M$. Here, $\epsilon ^{*}$ is the critical radius (see Definition \ref{def_critical_radius}) and $\kappa(M)$ is a constant determined by $M$ (see Table \ref{table: constants}). We will revisit the monotonic Bernoulli example in Section \ref{section: example} to explicitly instantiate and derive these optimal rates.

\subsection{Related Literature}\label{subsection: literature review}

The minimax rate is an important metric for evaluating the intrinsic difficulty of a statistical estimation problem. Readers may consult books such as \cite{Emery2007,minimaxbook} for general results on minimax theory. Classical minimax theory has been extensively developed for nonparametric regression under various assumptions. The classical theory of minimax rates for nonparametric regression can be traced back at least to \cite{stone1980optimal}, who studied the estimation of a smooth regression function and its derivatives at a fixed point and showed that, when the regression function is \(p\)-smooth on \(\mathbb{R}^d\) and the target is a derivative of order \(m\), the optimal rate of convergence is \(n^{-(p-m)/(2p+d)}\). Building on this classical rate result, \cite{fan1993local} studied local linear regression smoothers for estimating \(m(x_0)=\mathbb{E}(Y\mid X=x_0)\) at a fixed point \(x_0\). The paper considered the minimax risk \(R(n,\mathcal{C}_2)=\inf_{\widehat{m}}\sup_{f\in\mathcal{C}_2}\mathbb{E}_f[(\widehat{m}(x_0)-m_f(x_0))^2]\), where \(n\) is the sample size and \(\mathcal{C}_2\) is a class of data-generating distributions whose regression functions satisfy a bounded second-order smoothness condition. The local linear regression smoother attains the optimal rate \(n^{-4/5}\). Moreover, with an appropriate choice of kernel and bandwidth, it is nearly optimal in terms of the asymptotic minimax constant and asymptotically optimal within a broad class of linear smoothers. In high-dimensional settings, \cite{raskutti2009lower} studied sparse additive nonparametric regression models of the form \(f^*(X_1,\ldots,X_p)=\sum_{j\in S}h_j^*(X_j)\), where \(S\subset\{1,\ldots,p\}\) is an unknown active subset with \(|S|=s\), and each component \(h_j^*\) belongs to a univariate function class \(\mathcal H\). They established a minimax lower bound under squared \(L_2(\mathbb P)\) loss of order \(\max\left\{\frac{s\log(p/s)}{n},\;s\epsilon_n^2(\mathcal H)\right\}\), where \(n\) is the sample size, \(p\) is the total number of predictors, \(\mathbb P\) is the covariate distribution, and \(s\epsilon_n^2(\mathcal H)\) reflects the difficulty of estimating \(s\) univariate component functions. This result was further developed by \cite{JMLR:v13:raskutti12a}, who proposed a computationally tractable convex-regularized estimator over reproducing kernel Hilbert spaces and established matching upper bounds. Consequently, the minimax-optimal squared \(L_2(\mathbb P)\) rate is \(\frac{s\log(p/s)}{n}+s\epsilon_n^2(\mathcal H)\) up to constant factors. For an \(m\)-th order Sobolev class, where \(m\) is the smoothness order, this becomes \(\frac{s\log(p/s)}{n}+s n^{-2m/(2m+1)}\). \cite{yang2015minimax} further studied minimax \(L_2\) risks for high-dimensional nonparametric regression, where the number of predictors \(p\) may grow with the sample size \(n\). They considered sparse regression functions that depend only on \(d\) important predictors among the \(p\) available predictors, with \(d=\mathcal{O}(\log n)\) and \(\log p=o(n)\), as well as sparse additive models whose additive components may involve different subsets of predictors and different smoothness levels. They showed that suitable Bayesian Gaussian process procedures can adaptively  attain the optimal minimax rates up to logarithmic factors. Related minimax results have also been obtained for spline estimators \cite{10.1214/aos/1176349650,10.1214/aos/1176349651}, isotonic regression \cite{chatterjee2015risk,10.1214/18-AOS1792}, and other shape-constrained regression problems. \cite{Cai_2012} provided a comprehensive review of minimax theory and adaptive inference in nonparametric function estimation, pointwise estimation, and nonparametric confidence intervals.

However, minimax results for nonparametric regression in the exponential family setting, especially those concerning the estimation of natural parameters, appear to be much less developed. The exponential family is fundamental to statistical theory and practice due to its wide applicability and well-defined theoretical framework. For classical results on the exponential family, interested readers may consult books such as \cite{mccullagh1989generalized,Casella2002,Bickel2015}. Nevertheless, there are relatively fewer works that discuss minimax rates for parameter estimation in exponential families. Among the most relevant are \cite{cai2010nonparametric} and \cite{brown2010nonparametric}. \cite{cai2010nonparametric} studied nonparametric regression within natural exponential family models and derived minimax rates for estimating the natural parameter function. \cite{brown2010nonparametric} further developed this direction by introducing a mean-matching variance stabilizing transformation, which allows techniques from Gaussian nonparametric regression to be applied to a broader class of natural exponential family models. These works are closely related to ours in that they also study minimax estimation in exponential family models. However, they mainly focus on nonparametric regression problems where the natural parameter is modeled as a smooth function. In contrast, our work studies the estimation of a finite-dimensional natural parameter vector \(\theta \in K \subset \mathbb{R}^n\) under a general star-shaped constraint. The difficulty of our problem is characterized by the local geometry of the constraint set \(K\), through local packing or local entropy quantities, rather than by classical smoothness assumptions on an underlying regression function.

Geometric quantities such as metric entropy and moduli of continuity have played an important role in the characterization of minimax rates \cite{donoho1991geometrizing, yang1999information}, and a recent line of work has shown that local geometric quantities can sharply characterize minimax rates. \cite{neykov2023minimax} derived the minimax rate of the Gaussian sequence model under bounded convex constraints purely in terms of the local geometry of the constraint set \(K\). Specifically, the minimax risk under squared \(\ell_2\) loss is, up to constant factors, $ \epsilon^{*2} \wedge \operatorname{diam}(K)^2, \epsilon^*  = \sup\left\{ \epsilon: \frac{\epsilon^2}{\sigma^2} \leq \log M^{\mathrm{loc}}(\epsilon,c) \right\}$, where \(\log M^{\mathrm{loc}}(\epsilon,c)\) denotes the local entropy (see Definition \ref{def_local_entropy}). More recently, \cite{prasadan2024characterizing} characterized the minimax rate of nonparametric regression over bounded star-shaped function classes. They showed that the minimax risk under squared population \(L_2\) loss is, up to constant factors,
$ \epsilon^{*2} \wedge \operatorname{diam}(\mathcal F)^2, $
where \(\epsilon^*\) is determined by the local metric entropy of the function class. Their work demonstrates that the local geometry of a star-shaped constraint can determine the intrinsic difficulty of nonparametric regression.

Our work is closely related in spirit to these local geometry based minimax results, but differs in both the statistical model and the estimation target. Instead of estimating a regression function in a Gaussian or sub-Gaussian regression model, we study the estimation of the natural parameter vector in an exponential family model. The main contribution of this paper is to extend the local-geometry characterization of minimax rates to a broad class of sufficiently smooth nonsingular exponential families under star-shaped constraints.

\subsection{Our Contributions}\label{subsection: contributions}

Although some special cases of our problem have been studied extensively (see the \hyperref[subsection: literature review]{related literatures} above), the optimal rate of estimation for the natural parameter $\theta$ in general exponential families remains an open problem. In this paper, we bridge this gap by making the following primary contributions:

\begin{itemize}
    \item We establish the information-theoretic lower bound for the minimax risk \eqref{eq: minimax risk of the estimation} of estimating the natural parameter $\theta_{0}$ under any bounded, star-shaped constraint set $K$ via the characterization of its localized geometry. (Section \ref{section: theoretical results})
    \item Inspired by \cite{10.1214/25-AOS2576}, we propose a novel tree-based algorithm incorporating iterative pruning layers. We prove that the sequence of estimators constructed from this procedure achieves a minimax risk that tightly matches the theoretical lower bound, provided the number of iterations is appropriately controlled. (Section \ref{section: theoretical results})
    \item We address the inherent computational bottlenecks by developing a computationally efficient variant of our approach. Under slightly stronger structural conditions on the constraint set $K$, this algorithm runs in polynomial time with respect to the dimensional parameter $n$ and the geometric parameters of $K$. Furthermore, we demonstrate that this efficient estimator is near-optimal, achieving a minimax risk that matches the theoretical lower bound up to a poly-logarithmic factor. (Section \ref{section: efficient estimator})
\end{itemize}

\subsection{Organization}\label{subsection: organization}

The rest of the paper is organized as follows. Section \ref{section: theoretical results} presents our main theoretical results. Specifically, in Section \ref{subsection: lower bound}, we introduce bounds on the Kullback-Leibler (KL) divergence alongside Fano's inequality, which serve as the foundation for establishing the information-theoretic minimax lower bound. In Section \ref{subsection: upper bound}, we propose the theoretical procedure, Algorithm \ref{algorithm: theoretical algorithm}, and derive the corresponding minimax upper bound based on its statistical performance. Section \ref{SecMinimax} synthesizes these results, demonstrating that the upper and lower bounds match up to a constant depending only on $M$, thereby characterizing the exact minimax rate. To address the computational intractability of Algorithm \ref{algorithm: theoretical algorithm}, Section \ref{section: efficient estimator} introduces a computationally efficient alternative. We prove that the estimator $\hat{\theta} _{\text{p}}$ outputted by Algorithm \ref{algorithm: main algorithm} achieves near-optimal statistical guarantees relative to the established minimax rate. In Section \ref{section: example}, we instantiate our framework with some concrete examples, namely, simple one-dimensional estimation, constrained Bernoulli regression, and constrained Logistic regression. We comprehensively analyze these models and derive their corresponding minimax rates. Finally, Section \ref{section: discussion} concludes the paper with a brief discussion and some possible directions for future works. For clarity of exposition, the proofs of certain lemmas, theorems, and corollaries are deferred to Appendix \ref{section: deferred proofs}.

\subsection{Notations and Definitions}\label{subsection: notations and definitions}

In this section, we introduce some commonly used notation. We use $n$ to denote the full dimension of the model and hence the dimension of the natural parameter $\theta $. In addition, $\theta _{0}$ represents the true parameter. For any $a, b \in \mathbb{R}$, we write $a \vee b := \max\limits \left\{a,b\right\}$ and $a \wedge b := \min\limits \left\{a,b\right\}$. For any integer $m \in \mathbb{N}$, we use the shorthand $[m] = \left\{1,2,\dots,m\right\}$. Throughout the paper, $\log$ denotes the natural logarithm, and $\left\lVert \cdot \right\rVert_{2}^{}$ denotes the Euclidean norm. To simplify the notations, we use $\operatorname{KL}(\theta _{i}\left\lVert\right. \theta _{i}^{\prime })$ to denote the Kullback-Leibler (KL) divergence between two probability laws from the studied exponential family with natural parameter $\theta _{i}$ and $\theta _{i}^{\prime }$, and we let $\operatorname{KL}(\theta \left\lVert\right. \theta ^{\prime }):=\sum\limits _{i=1}^{n}\operatorname{KL}(\theta _{i}\left\lVert\right. \theta _{i}^{\prime })$ for $\theta =(\theta _{1},\dots,\theta _{n})^{\top } ,\theta ^{\prime }=(\theta _{1}^{\prime },\dots,\theta _{n}^{\prime })^{\top } $. We use $B_{2}(\theta, r)$ to represent a closed Euclidean ball centered at $\theta$ with radius $r>0$. Let $d=\text{diam}(K) $ represent the diameter of the constraint $K$, which is trivially defined as $\text{diam}(K) :=\sup\limits _{\theta ,\theta ^{\prime }\in K}\left\lVert \theta -\theta ^{\prime }\right\rVert_{2}^{}$. For two sequences $a_n$ and $b_n$, we write $a_n \lesssim b_n$ (and equivalently, $b_n \gtrsim a_n$) to denote inequalities up to a universal constant factor. We write $a_n \asymp b_n$ if both $a_n \lesssim b_n$ and $a_n \gtrsim b_n$ hold. Similarly, $a _{n}\lesssim _{M}b _{n}$ (and also for $a _{n}\gtrsim _{M}b _{n}, a _{n}\asymp _{M}b _{n}$) means the rates between $a _{n}$ and $b _{n}$ except for a constant possibly depending only on $M$. The absolute constants may vary from line to line without explicit mention.

It is also important to note that thought some crucial quantities in the analysis and expressions hereafter are denoted to be dependent only on $M$, by our argument in the first paragraph of Section \ref{subsection: problem formulation}, they are the results of simplification, and rigorously speaking, should be understood as depending on $M _{l},M _{u},A(\cdot )$. As we provide detailed derivations for these quantities, it is not hard for the reader to derive them under the general conditions (i.e., $K \subset [M _{\text{l}},M _{\text{u}}]^{n}\subset \mathcal{D}$ with $\mathcal{D}$ being an open set).

Finally, we provide several key definitions that are extensively used in the subsequent analysis.

\begin{definition}[$\delta$-covering set and covering number of $K$]\label{def_covering}
A $\delta$-covering set of a set $K$ with respect to the norm $\left\lVert \cdot \right\rVert_{}^{}$ is a set $\left\{\theta _{1},\theta _{2},\dots,\theta _{N}\right\}\subset K$ such that for every $\theta \in K$, there exists some $i \in [N]$ such that $\left\lVert \theta -\theta _{i}\right\rVert_{}^{}\le \delta $. The minimal possible $N$ is defined as the covering number of $K$, denoted as $N(\left\lVert \cdot \right\rVert_{}^{},\delta ,K)$. In this work, we exclusively study the Euclidean norm and the covering number is denoted as $N(\delta ,K)$ for convenience.
\end{definition}

\begin{definition}[$\delta$-packing set and packing number of $K$]\label{def_packing}
A $\delta$-packing set of a set $K$ with respect to the norm $\left\lVert \cdot \right\rVert_{}^{}$ is a set $\left\{\theta _{1},\theta _{2},\dots,\theta _{M}\right\}\subset K$ such that $\left\lVert \theta _{i}-\theta _{j}\right\rVert_{}^{}$ for all distinct $i, j \in [M]$. The maximal possible $M$ is defined as the packing number of $K$, denoted as $M(\left\lVert \cdot \right\rVert_{}^{},\delta ,K)$. Similar to the covering number, we use $M(\delta ,K)$ for convenience.
\end{definition}

\begin{definition}[Local Entropy]\label{def_local_entropy}
    Given a set $K$ and a universal constant $c > 0$, the local entropy of $K$ is defined as
    \begin{equation}\label{eq: definition of local entropy}
        N^{\text{loc}}(\epsilon, c) := \sup_{\theta \in K} N(\epsilon/c, B_{2}(\theta, \epsilon) \cap K).
    \end{equation}
    Throughout the paper, the letter $c$ is reserved to denote the constant introduced in this definition. Furthermore, we assume $c$ to be sufficiently large.
\end{definition}
\begin{remark*}[]
    We note that for star-shaped set, it is not hard to show that $N ^{\text{loc}}(\epsilon ,c)$ is a non-increasing function in $\epsilon $ and non-decreasing in $c$. Interested readers can refer \cite{10.1214/25-AOS2576} for the proof. We can analogously define $M ^{\text{loc}}(\epsilon ,c)$. By the well-known properties of the maximum packing and minimal covering, we know that 
    \begin{equation*}
        N ^{\text{loc}}(\epsilon ,c)\le M ^{\text{loc}}(\epsilon ,c)\le N ^{\text{loc}}(\epsilon ,2c).
    \end{equation*}
\end{remark*}

As a result of the local entropy, a fundamental quantity in this work, the critical radius $\epsilon ^{*}$, is defined below.
\begin{definition}[critical radius]\label{def_critical_radius}
    Given the constraint $K$ and some constants $c,\kappa (M)$ only depending on $M$, the critical radius is defined as
    \begin{equation}\label{eq: def_epsilon_star}
        \epsilon ^{*}:= \sup \{\epsilon \left\lvert\right. \epsilon^2 \kappa(M) \le \log N^{\text{loc}}(\epsilon,c)\}.
    \end{equation}
    The choice of $c$ and $\kappa (M)$ will be specified later.
\end{definition}

The following well-known definition of the star-shaped set is a generalization of the convex set.

\begin{definition}[Star-shaped set]\label{def_star_shaped_set}
A set $K \subset \mathbb{R}^{n}$ is called star-shaped with respect to a point $x _{0}\in K$ if for each $x \in K$, the line segment connecting $x _{0}$ and $x$ is entirely contained within $K$:
\begin{equation*}
    (1-t) x_0 + tx \in K, \forall x \in K, \forall t \in [0, 1].
\end{equation*}
We refer to any such point $x _{0}$ as a star center of $K$.
\end{definition}

\begin{remark*}[]
    We note that the class of star-shaped constraints is highly expressive, encompassing a broad range of commonly studied structures. For example, this framework naturally captures both standard convex constraints and non-convex sets, such as $s$-sparse constraints.
\end{remark*}

To conclude this section, we summarized the constants used in this paper in the following table.

\begin{table}[htbp]\label{table: constants}
\centering
\caption{Information of the Constants Appeared in the Paper}
\begin{tabular}{|c|c|l|}
\hline
\textbf{Constants} & \textbf{Description or Explicit Expression} & \textbf{Original Definition} \\
\hline
$N ^{\text{loc}}(\epsilon ,c)$ & $\sup_{\theta \in K} N(\epsilon/c, B_{2}(\theta, \epsilon) \cap K)$ & Definition \ref{def_local_entropy}\\ 
\hline
$C(M)$     & $\max\limits _{-M \le \theta \le M}A ^{\prime \prime }(\theta )$ & Lemma \ref{lemma: bounds on KL divergence} \\
\hline
$c(M)$       & $\min\limits _{-M \le \theta \le M}A ^{\prime \prime }(\theta )$ & Lemma \ref{lemma: bounds on KL divergence} \\
\hline
$C'(M)$  & $\max\limits _{-2M \le \theta \le 2M}A ^{\prime \prime }(\theta )$ & Appendix \ref{subsubsection: proof of sub-exponential of T} \\
\hline
$C$            & A large constant required to be greater than $2\sqrt{\frac{C(M)}{c(M)}+8}$ & Lemma \ref{lemma: high probability for close distance between theta prime and theta prime prime} \\
\hline
$c$          & A large constant required to be greater than $2(C+1)$ & Definition \ref{def_local_entropy} \\
\hline
$\kappa(M)$        & $\frac{c(M)^2 C^2}{32 C'(M)}\wedge  \frac{|C(M) - c(M)(C-1)^2|}{4}$ & Appendix \ref{subsubsection: proof of high probability for close distance between theta prime and theta prime prime} \\
\hline
$\epsilon ^{*}$ & $\sup \{\epsilon \left\lvert\right. \epsilon^2 \kappa(M) \le \log N^{\text{loc}}(\epsilon,c)\}$ & Definition \ref{def_critical_radius}\\ 
\hline
$\zeta$     & $\frac{1}{8C ^{\prime }(M)c(M)}\wedge \frac{1}{4}$ & Remark \ref{remark: comparison between kappa and C(M)} \\
\hline
$\tau _{K}$ & The ``small'' constant of the set $K$ & Theorem \ref{theorem: main theorem of the cited paper}, Lemma \ref{lemma: kolmogorov width control}\\
\hline
$T _{2}(K)$ & The type-$2$ constant of the set $K$ & Definition \ref{def_type_2_condition}\\ 
\hline
$r$ & The known constant such that $B _{2}(0,r) \subset K$ & Definition \ref{def_well_balanced_set}\\
\hline
$R$ & The known constant such that $K \subset B _{2}(0,R)$ & Definition \ref{def_well_balanced_set}\\
\hline
\end{tabular}
\end{table}

\section{Theoretical Results}\label{section: theoretical results}

In this section, our goal is to derive matching lower and upper bounds on the $\ell _{2}$ estimation error, thereby establishing the minimax rate (up to universal constants) for the class of star-shaped constrained, sufficiently smooth, and nonsingular canonical exponential family distributions defined in the \hyperref[subsection: problem formulation]{problem formulation}. We begin by establishing the fundamental lower bound via a precise control of the KL divergence and Fano's inequality.

\subsection{Lower Bound}\label{subsection: lower bound}

This subsection is devoted to establishing the minimax lower bound. Our analysis begins with Lemma \ref{lemma: bounds on KL divergence}, which provides a tight control on the Kullback-Leibler (KL) divergence. Specifically, this lemma quantifies the statistical discrepancy between any two parameters $\theta$ and $\theta ^{\prime }$ directly in terms of their squared Euclidean distance.

\begin{lemma}\label{lemma: bounds on KL divergence}
  For any two parameters $\theta, \theta' \in [-M, M]^n$ associated with a nonsingular exponential family possessing a twice continuously differentiable log-partition function, the KL divergence satisfies
  \begin{equation*}
    c(M)\left\lVert \theta ^{\prime }-\theta \right\rVert_{2}^{2}\le \operatorname{KL}(\theta \left\lVert\right. \theta ^{\prime })\le C(M)\left\lVert \theta ^{\prime }-\theta \right\rVert_{2}^{2},
  \end{equation*}
  where $c(M)$ and $C(M)$ are defined as 
  \begin{equation*}
    \begin{aligned}
        & C(M):=\max\limits _{-M \le \theta \le M}A ^{\prime \prime }(\theta ),\\ 
        & c(M):=\max\limits _{-M \le \theta \le M}A ^{\prime \prime }(\theta ).
    \end{aligned}
  \end{equation*}
\end{lemma}

In Lemma \ref{lemma: bounds on KL divergence}, we require the exponential family distribution to be non-degenerate, meaning $\text{Var}_{\theta_i}(T(Y_i)) > 0$ (or equivalently, $A''(\theta_i) > 0$) for all $\theta_i \in [-M, M]$. This condition, combined with the compactness of the parameter space $[-M, M]$, is necessary to ensure that the sufficient statistic $T(Y_i)$ is not almost surely constant and has bounded variance. Consequently, by the continuity of $A''(\cdot)$, the bounding constants $c(M)$ and $C(M)$ are respectively strictly positive and finite (i.e., $0 < c(M) \le A''(\theta_i) \le C(M) < \infty$). We defer the proof to Appendix \ref{subsubsection: proof of bounds on KL divergence}.

To establish the minimax lower bound, we first reduce the estimation task to a multiple hypothesis testing problem over a finite, separated parameter set. Next, we introduce Fano's inequality in Theorem \ref{theorem: fano inequality}. This fundamental information-theoretic result relates the probability of making an incorrect hypothesis selection to the mutual information between the observed data $Y$ and the true parameter index $J$. Crucially, by coupling Fano's inequality with the Kullback-Leibler divergence bounds established in Lemma~\ref{lemma: bounds on KL divergence}, we can control this mutual information, thereby guaranteeing a strictly positive lower bound on the estimation error.

\begin{theorem}[Generalized Fano's Method]\label{theorem: fano inequality}
    Let $K \subseteq \mathbb{R}^n$ be a parameter space, and let $\left\{\theta _{1},\dots,\theta _{m}\right\} \subset K$ be a finite packing set strictly separated by $\epsilon > 0$ in the Euclidean norm, meaning $\left\lVert \theta _{i}-\theta _{j}\right\rVert_{2}^{}\ge \epsilon $ for all $i \neq j$. Suppose $J$ is a random variable uniformly distributed over the index set $[m]$, and given $J = j$, the data is generated as $Y \sim P_{\theta^j}$. Then, the minimax estimation risk $\mathcal{R}_{n}(f,K)$ is lower bounded by
    \begin{equation*}
        \mathcal{R}_{n}(f,K) \ge \frac{\epsilon^2}{4} \bigg(1 - \frac{I(Y; J) + \log 2}{\log m}\bigg),
    \end{equation*}
    where $I(Y; J)$ denotes the mutual information between the observation $Y$ and the parameter index $J$.
\end{theorem}

To apply Theorem \ref{theorem: fano inequality}, we must control the mutual information $I(Y; J)$. By standard information-theoretic identities and the convexity of the Kullback-Leibler divergence (see, e.g., Wainwright \cite{wainwright2019}, Eq. 15.52), the mutual information can be bounded by the average KL divergence to any arbitrary distribution $P_\theta$:
\begin{equation*}
    I(Y; J) \le \frac{1}{m} \sum_{j=1}^{m} \operatorname{KL}(P _{\theta _{j}}\left\lVert\right. P _{\theta}), \quad \text{for any } \theta \in \mathbb{R}^n.
\end{equation*}
Utilizing the bounds established in Lemma \ref{lemma: bounds on KL divergence}, we can further upper bound the KL divergence by the Euclidean distance, yielding:
\begin{equation*}
    I(Y; J) \le \frac{1}{m} \sum_{j=1}^{m} C(M) \left\lVert \theta _{j}-\theta \right\rVert_{2}^{2} \le \max_{j \in [m]} C(M) \left\lVert \theta _{j}-\theta \right\rVert_{2}^{2}.
\end{equation*}

By combining Lemma \ref{lemma: bounds on KL divergence} with Theorem \ref{theorem: fano inequality}, we can now establish a localized minimax lower bound. This is achieved by balancing the local metric entropy of the parameter space against the information capacity of the observations. The detailed proof is placed in Appendix \ref{subsubsection: proof of minimax lower bound}.
\begin{lemma}\label{lemma: minimax lower bound}
    Let $N^{\text{loc}}(\epsilon, c)$ denote the local packing number as specified in Definition \ref{def_local_entropy}, where $c > 0$ is a sufficiently large, fixed constant. For any $\epsilon > 0$ that satisfies the entropy condition
    \begin{equation*}
        \log N^{\text{loc}}(\epsilon, c) > 4 \big(C(M)\epsilon^2 \vee \log 2 \big),
    \end{equation*}
    the minimax estimation risk over the parameter space $K$ is lower bounded by
    \begin{equation*}
        \mathcal{R}_{n}(f,K) \ge \frac{\epsilon^2}{8c^2}.
    \end{equation*}
\end{lemma}

\begin{remark}
\label{remark: minimax lower bound to minimax rate}
We shall prove later in Theorem \ref{theorem: main theorem for theoretical algorithm} that Lemma \ref{lemma: minimax lower bound} and the condition that $K$ is a star-shaped constraint yield a lower bound on the minimax rate of the $\ell _{2}$ error as $\epsilon ^{*2}\wedge \text{diam}(K) ^{2}$, where we recall that $\epsilon ^{*}$ is defined as 
\begin{equation*}
   \epsilon^*=\sup\left\{\epsilon: \kappa(M)\epsilon^2\le\log N^{\mathrm{loc}}(\epsilon,c)\right\}, 
\end{equation*}
and $\kappa (M)$ is a constant depending only on $M$.
\end{remark}

\subsection{Upper Bound}\label{subsection: upper bound}

In this section, we focus on deriving an upper bound for the minimax rate defined in \eqref{eq: minimax risk of the estimation}. We take a constructive approach: we introduce Algorithm \ref{algorithm: theoretical algorithm} to generate an estimator and subsequently bound its $\ell _{2}$ error. Later, in Section \ref{SecMinimax}, we demonstrate that this upper bound matches our lower bound, thereby characterizing the minimax rate up to constants depending only on $M$.

Before introducing Algorithm \ref{algorithm: theoretical algorithm}, we establish several key definitions. These concepts not only facilitate the presentation of the algorithm but also play a crucial role in our subsequent proofs.

\begin{definition}\label{def_dominating_point}
Let $\theta^\alpha$ and $\theta^\beta$ be two valid candidates in Algorithm \ref{algorithm: theoretical algorithm}. If
\begin{equation*}
    \sum_{i=1}^{n} \left[ \left(\theta^\alpha_i - \theta^\beta_i \right) T(Y _{i}) + A(\theta^\beta_i) - A(\theta^\alpha_i) \right] \le 0,
\end{equation*}
we say that $\theta^\beta$ dominates $\theta^\alpha$, denoted by $\theta^\beta \succ \theta^\alpha$. Equivalently, we say that $\theta^\alpha$ is dominated by $\theta^\beta$ ($\theta^\alpha \prec \theta^\beta$).
\end{definition}

Definition \ref{def_dominating_point} formalizes the comparison of two candidate points based on their relative likelihoods. Specifically, it identifies the parameter that yields a higher likelihood under the specified model, establishing a principled criterion for selecting a superior candidate.

\begin{definition}\label{def_conditional_dominating_point}
Suppose we are given a candidate set $S = \{\theta^1, \dots, \theta^N\}$. For any $\alpha \in [N]$ and a specified constant $C > 0$, we define the discrepancy function $H(\delta, \theta^\alpha, S)$ as 
\begin{equation*}
    H(\delta, \theta^\alpha, S) := \max \Big\{ \|\theta^\alpha - \theta^\beta\|_{2} \;\Big|\; \beta \in [N], \: \theta^\beta \succ \theta^\alpha, \text{ and } \|\theta^\beta - \theta^\alpha\|_{2} \ge C \delta \Big\},
\end{equation*}
with the convention that the maximum over an empty set evaluates to $0$.
\end{definition}

Definition \ref{def_conditional_dominating_point} quantifies the maximum distance between a given candidate $\theta^\alpha$ and any other point in $S$ that strictly dominates it, provided this distance exceeds the margin $C\delta$. Intuitively, $H(\delta, \theta^\alpha, S)$ acts as a spatial screening metric: a strictly positive value indicates that $\theta^\alpha$ is significantly outperformed by a distant candidate, while a value of zero implies that $\theta^\alpha$ is either globally non-dominated or that all dominating points lie within a tight $\mathcal{O}(\delta)$ neighborhood. Consequently, this function serves as a crucial regularizing criterion in the iterative refinement process of Algorithm \ref{algorithm: theoretical algorithm}.

\subsubsection{Tree-Based Algorithm}\label{subsubsection: theoretical algorithm}

In this section, inspired by \cite{10.1214/25-AOS2576}, we design a tree-based algorithm to generate an effective estimator, thereby deriving an upper bound for the minimax rate \eqref{eq: minimax risk of the estimation}. The algorithm constructs a sequence that converges to a desired optimal point. In essence, we first build a directed tree of points in $K$, the construction of which is detailed in Section \ref{subsubsection: algorithm explain} and Algorithm \ref{algorithm: construction of the pruned tree}. Subsequently, we utilize the observations to iteratively traverse this tree, selecting the limiting point (or any point sufficiently far ahead in the sequence) as our final estimator, as detailed in Algorithm \ref{algorithm: theoretical algorithm}.

Our approach differs from Algorithm 1 in \cite{neykov2023minimax} in a crucial way: instead of selecting a vector that minimizes the $\ell_2$ norm over a specific discrete offspring set $\mathcal{O}(\Upsilon_k)$, our procedure selects a candidate $\theta \in \mathcal{O}(\Upsilon_k)$ that minimizes the discrepancy function $H\left(\frac{\text{diam}(K)}{2^{k}(C+1)}, \theta, \mathcal{O}(\Upsilon_k)\right)$ at each step $k$. Here, $\{\Upsilon_k\}$ for $k \in \NN$ denotes the sequence of iterates traversing the tree, and $\mathcal{O}(\Upsilon_k)$ denotes the set of all offspring of the node $\Upsilon_k$.

Furthermore, drawing upon insights from \cite{10.1214/25-AOS2576}, we introduce a pruning step to the directed tree to resolve a technical oversight in the original analysis by \cite{neykov2023minimax}. For completeness, we conceptually outline this tree-based algorithm and its associated pruning mechanism below.

\subsubsection{Algorithm Illustration and Pruning Procedure}\label{subsubsection: algorithm explain}

In this section, we provide a high-level overview of Algorithm \ref{algorithm: construction of the pruned tree}. Before proceeding, we first establish the necessary notation.

For a directed tree $G$, we denote the set of parents of any node $u \in G$ by $\mathcal{P}(u)$ (we allow multiple parents for a node), such that there is a directed edge $\mathcal{P}(u) \to u$. Conversely, we denote the offspring set of a node $v$ by $\mathcal{O}(v) := \{u \in G \mid \mathcal{P}(u) = v\}$. Each node $u$ is assigned a depth or level $k$, where the root node is situated at level 1, and all offspring of a level-$k$ node belong to level $k+1$ We let $\mathcal{L}(k)$ be the nodes in the level $k$.

The construction initializes by arbitrarily selecting a root node $\Upsilon_1 \in K$, which constitutes level 1. Level 2 is formed by taking a maximal $(d/c)$-packing of the set $B(\Upsilon_1, d) \cap K = K$, with directed edges added from $\Upsilon_1$ to each packing point. No pruning is applied at level 2. For any level $k \ge 3$, assuming level $k-1$ has been constructed and pruned, we proceed as follows: for each node $u$ at level $k-1$, we compute a maximal $d/(2^{k-1}c)$-packing of $B(u, d/2^{k-2}) \cap K$. These packing points serve as preliminary candidates for level $k$, and directed edges are initially drawn from $u$ to each point in its respective packing set.

Next, we lexicographically order these preliminary points at level $k$, denoting them as $u^k_1, \ldots, u^k_{N_k}$, and initialize a list of unprocessed nodes $\mathcal{U}_k = [u^k_1, \ldots, u^k_{N_k}]$.

While $\mathcal{U}_k$ remains non-empty, we select its first available element, say $u^k_i$. We then define the neighborhood set
\begin{equation*}
    \mathcal{T}_k(u^k_i) := \left\{ u^k_j \in \mathcal{U}_k \;\middle|\; \|u^k_j - u^k_i\| \le \frac{d}{2^{k-1}c}, \; j \neq i \right\}.
\end{equation*}

\begin{figure}[ht]
    \includegraphics[width=\textwidth]{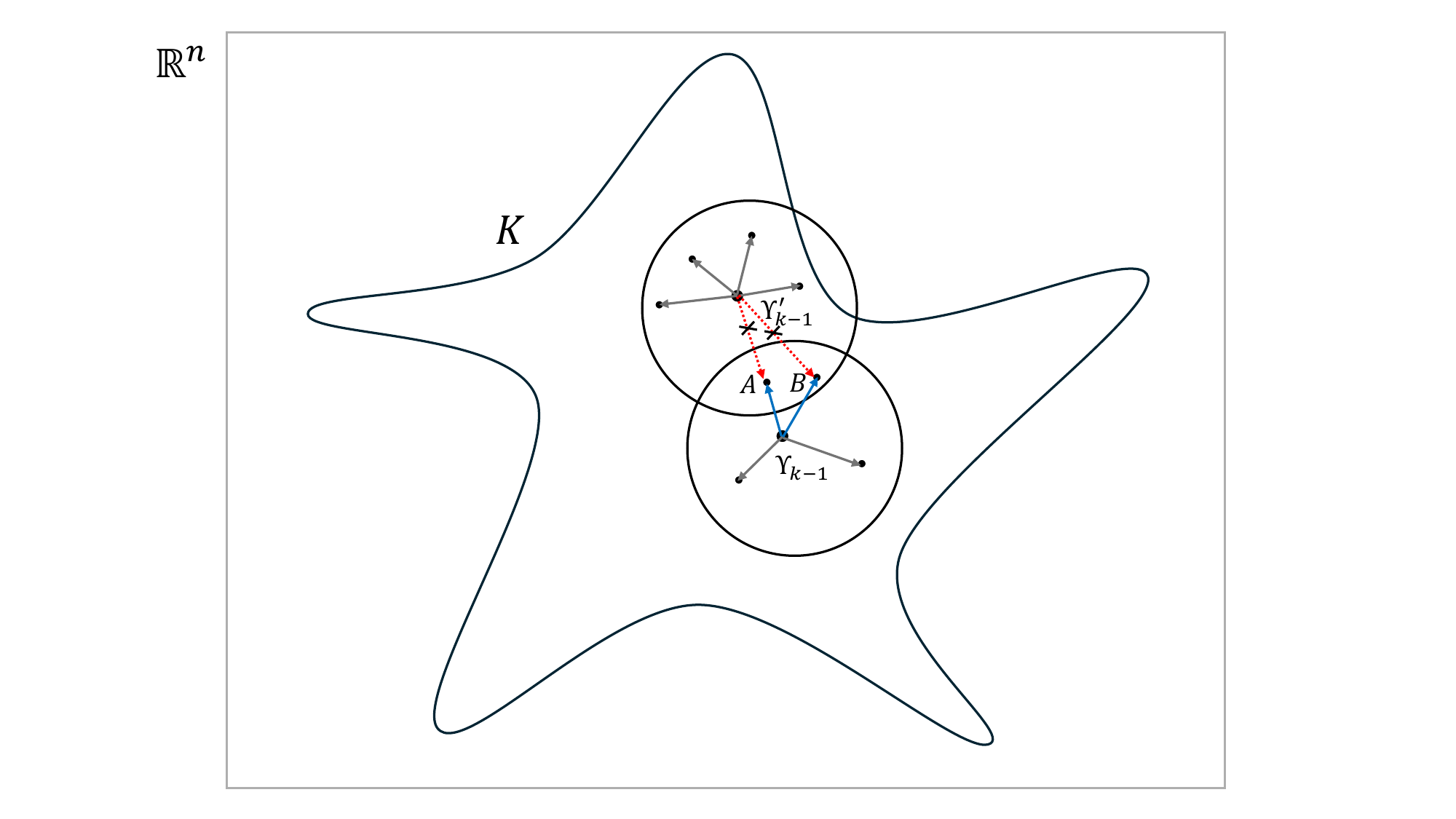}
    \caption{the illustration of the pruning procedure. For the currently first point $\Upsilon _{k-1}$ in $\mathcal{U}_{k-1}$ and the points $A,B$ that are originally $\Upsilon _{k-1}^{\prime }$'s offsprings, we remove the edges from $\Upsilon _{k-1}^{\prime }$ to $A,B$ and add edges from $\Upsilon _{k-1}$ to $A,B$ if $\left\lVert A-\Upsilon _{k-1}\right\rVert_{2}^{}\le d/(2 ^{k-1}c)$ and $\left\lVert B-\Upsilon _{k-1}\right\rVert_{2}^{}\le d/(2 ^{k-1}c)$.}\label{figure: illustration of the pruning procedure}
\end{figure}

By definition, $\mathcal{T}_k(u^k_i) \subseteq \mathcal{U}_k$, ensuring this set consists exclusively of unprocessed nodes. For each node $u^k_j \in \mathcal{T}_k(u^k_i)$ and its parent $\mathcal{P}(u^k_j)$ (note that the parent is unique for any unprocessed point in $\mathcal{U}_{k}$), we delete the existing edge $\mathcal{P}(u^k_j) \to u^k_j$ and replace it with a new directed edge $\mathcal{P}(u^k_j) \to u^k_i$. See Figure \ref{figure: illustration of the pruning procedure} for an illustration.

Following this reassignment, we remove the set $\{u^k_i\} \cup \mathcal{T}_k(u^k_i)$ from $\mathcal{U}_k$. This pruning loop continues until $\mathcal{U}_k$ is empty. The resulting offspring nodes of level $k-1$ constitute the finalized, pruned level $k$. Then we work with the pruned level $k$ and construct the level $k+1$. This procedure is repeated for all $k \ge 3$ to construct the complete pruned graph, as shown in Algorithm \ref{algorithm: construction of the pruned tree}. We note that while this directed tree architecture originates from \cite{10.1214/25-AOS2576} while the subsequent Algorithm \ref{algorithm: theoretical algorithm} adapts the construction but employs a distinct selection criterion.

\begin{algorithm}[ht]
    \caption{Construction of the pruned tree}\label{algorithm: construction of the pruned tree}

    \KwInput{Constraint $K$, root node $\Upsilon _{1}$}

    \KwOutput{The pruned tree $G$}

    Construct maximal $d/c$-packing of $B(\Upsilon _{1},d)\cap K$ as $\mathcal{L}(2)$

    Add edges from $\Upsilon _{1}$ to each node in $\mathcal{L}(2)$

    $k\gets 3$

    \While{TRUE}{
        $\mathcal{U}_{k}\gets\left\{\right\}$

        \For{$u$ in $\mathcal{L}(k-1)$}{

            Construct $d/(2 ^{k-1}c)$ packing of $B _{2}(u,d/2 ^{k-2})\cap K$

            Add the points from the packing to $\mathcal{U} _{k}$

            Add edges from $u$ to each point from the packing

        }

        Sort $\mathcal{U}_{k}$ by the lexicographical order

        \While{$\mathcal{U}_{k}\neq \varnothing $}{
            Select first element of $\mathcal{U}_{k}$, denoted as $\tilde{u}$

            $\mathcal{T}_{k}(\tilde{u})\gets\left\{u \in \mathcal{U}_k \middle| \|u - \tilde{u}\| \le \frac{d}{2^{k-1}c}, u \neq \tilde{u}\right\}.$

            \For{$u ^{\prime }\in \mathcal{T}_{k}(\tilde{u})$}{
                Remove the edge from $\mathcal{P}(u ^{\prime })$ to $u ^{\prime }$

                Add edge from $\mathcal{P}(u ^{\prime })$ to $\tilde{u}$
            }
            Remove $\left\{\tilde{u}\right\}\cup \mathcal{T}_{k}(\tilde{u})$ from $\mathcal{U}_{k}$
        }

        $\mathcal{L}(k)\gets \bigcup\limits _{u \in \mathcal{L}(k-1)}^{}\mathcal{O}(u)$

        $k\gets k+1$
    }

    \Return $G$
    
\end{algorithm}

The iterative construction and pruning procedures described above yield the following crucial lemmas from \cite{10.1214/25-AOS2576}, which establishes the covering and packing properties at each level of the directed tree $G$. Specifically, the nodes at each level $k$ constitute a valid covering and packing of the constraint set $K$. Moreover, for any node $\Upsilon _{k-1}$ at level $k-1$, its offspring set $\mathcal{O}(\Upsilon _{k-1})$ forms an effective cover for the localized neighborhood $\mathcal{B}(\Upsilon _{k-1}, d/2^{k-2}) \cap K$. These geometric structural properties play a foundational role in our subsequent analysis of the estimator's performance.

\begin{lemma}[Lemma 3.1 of \cite{10.1214/25-AOS2576}]\label{lemma: covering and packing properties of the constructed tree}
Define \(\mathcal{L}(J)\) to be the set of nodes forming the pruned graph at level \(k\), and let $G$ be the pruned graph from above and assume $c > 2$. Then for any $J \ge 3$, $\mathcal{L}(J)$ forms a $d/(2^{J-2}c)$-covering of $K$ and a $d/(2^{J-1}c)$-packing of $K$. In addition, for each parent node $\Upsilon_{J-1}$ at level $J-1$, its offspring $\mathcal{O}(\Upsilon_{J-1})$ form a $d/(2^{J-2}c)$-covering of the set $\mathcal{B}(\Upsilon_{J-1}, d/2^{J-2}) \cap K$. Furthermore, the cardinality of $\mathcal{O}(\Upsilon_{J-1})$ is upper bounded by $N^{\mathrm{loc}}(d/2^{J-2}, 2c)$ for $J \ge 2$.
\end{lemma}

\begin{lemma}[Lemma 3.2 of \cite{10.1214/25-AOS2576}]\label{lemma: bound on the cardinality}
For any $\theta \in K$ and $J \ge 2$, $\mathcal{L}(J - 1) \cap B _{2}(\theta, d/2^{J-2})$ has a cardinality upper bounded by $M ^{\text{loc}}(d/2 ^{J-2},c)\le N^{\mathrm{loc}}(d/2^{J-2}, 2c)$.
\end{lemma}

\begin{lemma}[Lemma 3.3 of \cite{10.1214/25-AOS2576}]\label{lemma: bound on the distance between the parent and child}
Let $[\Upsilon_1, \Upsilon_2, \dots]$ be the nodes of an infinite path in the graph $G$, i.e., $\Upsilon_{J+1} \in \mathcal{O}(\Upsilon_J)$ for any $J \in \mathbb{N}$. Then for any integers $1 \le J ^{\prime }\le J$, we have $\|\Upsilon_{J'} - \Upsilon_J\| \le \frac{d(2 + 4c)}{c 2^{J'}}$.
\end{lemma}

Given the pruned tree $G$, the main algorithm runs by starting from the root node $\Upsilon _{1}$ and iteratively select the point from $\mathcal{O}(\Upsilon _{k})$ that minimizes $H \left(\frac{d}{2 ^{k}(C+1)}\theta ,\mathcal{O}(\Upsilon _{k})\right)$ for $\theta \in \mathcal{O}(\Upsilon _{k})$. Such point forms $\Upsilon _{k+1}$. The main algorithm outputs a sequence $\left\{\Upsilon _{1},\Upsilon _{2},\dots\right\}$, as detailed in Algorithm \ref{algorithm: theoretical algorithm}.

\begin{algorithm}[ht]
\caption{Minimax Optimal Algorithm}\label{algorithm: theoretical algorithm}
\SetKwComment{Comment}{/* }{ */}
\KwInput{Constraint $K$, observations $y _{1},\dots,y _{n}$, and the pruned tree $G$ from Algorithm \ref{algorithm: construction of the pruned tree}}
\KwOutput{The sequence $\Upsilon $}

Set $\Upsilon_1 \in K$ as an arbitrary point in $K$

$k \gets 1$

$\Upsilon \gets [\Upsilon_1]$

\While{TRUE}{
    $\Upsilon_{k+1} \gets \argmin\limits _{\theta \in \mathcal{O}(\Upsilon_k)} H\left(\frac{d}{2^{k}(C+1)}, \theta, \mathcal{O}(\Upsilon_k)\right)$ \Comment*[r]{Lexicographically sort points in the minimization to resolve ties by picking the smallest point}

    $\Upsilon$.append$(\Upsilon_{k+1})$

    $k \gets k + 1$
}

\Return $\Upsilon = [\Upsilon_1, \Upsilon_2, \ldots]$\Comment*[r]{Note: The sequence $\Upsilon$ forms a Cauchy sequence. For a proof of this fact, see Lemma 3.3 of \cite{10.1214/25-AOS2576}.}
\end{algorithm}

\subsubsection{Preliminary Lemmas}\label{subsubsection: preliminary lemmas}

In this section, we provide a theoretical analysis of Algorithm \ref{algorithm: theoretical algorithm}. We begin with Lemma \ref{lemma: sub-exponential of T}, which establishes that the sufficient statistic $T(Y _{i})$ is a sub-exponential random variable under our assumptions regarding the exponential family and the constraint $K$. This result serves as an important basis for subsequent Lemma \ref{lemma: high probability for close distance between theta prime and theta prime prime}. The proof of Lemma \ref{lemma: sub-exponential of T} relies on a direct evaluation of the moment generating function of $Y$ and classical properties of exponential families. We defer the detailed proof to Appendix \ref{subsubsection: proof of sub-exponential of T}.

\begin{lemma}\label{lemma: sub-exponential of T}
    The random variable $T(Y)$ is sub-exponential. Namely, there exists a constant $C'(M)$ such that $\log \mathbb{E} \exp [\lambda\cdot[T(Y) - \mathbb{E} T(Y)]] \le  \frac{\lambda^2 C'(M)}{2}$ for all $\lambda \in [-M, M]$.
\end{lemma}

Lemma \ref{lemma: high probability for close distance between theta prime and theta prime prime} parallels Lemma 2.5 in \cite{neykov2023minimax} and serves as a motivation for Definition \ref{def_conditional_dominating_point} and a theoretical guarantee for our selection criterion in Algorithm \ref{algorithm: theoretical algorithm}. Intuitively, this result establishes that if two candidate parameters are separated by a significant margin, the candidate closer to the true parameter $\theta $ is exponentially unlikely to be dominated by the further candidate (in the sense of Definition \ref{def_dominating_point}). By providing a sharp exponential tail bound on this misclassification probability, the lemma guarantees the theoretical soundness of our algorithm when pruning sub-optimal trajectories during each iteration.

\begin{lemma}\label{lemma: high probability for close distance between theta prime and theta prime prime}
    For any two vectors $\theta', \theta'' \in K$, define the indicator function
    \begin{equation*}
        \psi(Y) := \bm{1}_{\left\{ \sum_{i=1}^{n} \left[ (\theta_i' - \theta_i'')T(Y_i) + A(\theta_i'') - A(\theta_i') \right] < 0 \right\}}.
    \end{equation*}
    Assume that $\|\theta' - \theta''\| \ge C \delta$ for some large constant $C \ge \sqrt{\frac{4(C(M)+8c(M))}{c(M)}}$. Then, there exists a constant $\kappa(M) > 0$, also depending only on $M$, such that
    \begin{equation*}
        \sup_{\theta' : \|\theta' - \theta \| \le \delta} \mathbb{P}_{\theta }(\psi(Y) = 1) \le e^{-\kappa(M) \delta^2}.
    \end{equation*}
\end{lemma}

Based on Lemma \ref{lemma: high probability for close distance between theta prime and theta prime prime}, we present Lemma \ref{lemma: guarantee of the selection criterion}, which establishes a high-probability bound on the accuracy of the parameter estimation process within a localized neighborhood of the true parameter $\theta _{0}$. Specifically, this lemma demonstrates that the candidate selected by Algorithm \ref{algorithm: theoretical algorithm} is exponentially unlikely to deviate significantly from the true parameter, thereby guaranteeing the statistical consistency of our algorithmic selection.
\begin{lemma}[]\label{lemma: guarantee of the selection criterion}
    Suppose we are given a set of $N$ points $N_k := \{\theta^1,\ldots, \theta^N\} \subset K' \subseteq K$ that forms a $\delta$-cover of $K'$, and assume that the true parameter satisfies $\theta _{0}\in K'$. If we define the estimator as
    \begin{equation*}
        \overline{\theta} = \argmin_{\theta_\alpha \in N_k} H(\delta, \theta_\alpha, N_k), 
    \end{equation*}
    then for any fixed constant $C$ satisfying the condition in Lemma \ref{lemma: high probability for close distance between theta prime and theta prime prime}, $\overline{\theta}$ satisfies
    \begin{equation*}
        \mathbb{P} \left(\left\lVert \bar{\theta }-\theta \right\rVert_{2}^{}\ge (C+1)\delta \right) \le N e^{-\kappa(M) \delta^2},
    \end{equation*}
    where $\kappa(M)>0$ is the constant from Lemma \ref{lemma: high probability for close distance between theta prime and theta prime prime}.
\end{lemma}

We place the proofs of Lemma \ref{lemma: high probability for close distance between theta prime and theta prime prime} and Lemma \ref{lemma: guarantee of the selection criterion} to Appendix \ref{subsubsection: proof of high probability for close distance between theta prime and theta prime prime} and Appendix \ref{subsubsection: proof of guarantee of the selection criterion} for the reader's reference.

\subsubsection{Minimax Upper Bound}

We formally derive the theoretical guarantee of the performance of Algorithm \ref{algorithm: theoretical algorithm} in this section. Namely, Theorem \ref{theorem: performance of the theoretical algorithm} confirms that we should obtain an estimator $\theta ^{**}$ with its $\ell _{2}$ estimation error uniformly bounded by $\epsilon _{J ^{*}}\wedge d ^{2}$ (see definitions below), given the sufficiently large number of iterations we derive an upper bound of the minimax estimation for the estimator generated by Algorithm \ref{algorithm: theoretical algorithm}. This important result is established on the following Lemma \ref{lemma: property 1 of the estimators in the algorithm} and Lemma \ref{lemma: property 2 of the estimators in the algorithm} motivated from Lemma B.6 and Lemma B.7 of \cite{10.1214/25-AOS2576}. Specifically, given Condition \ref{eq: definition of J} detailed below, Lemma \ref{lemma: property 1 of the estimators in the algorithm} control the tail probability that $\Upsilon _{J}$ (the $J$-th iterate of Algorithm \ref{algorithm: theoretical algorithm}) deviates from the true parameter and Lemma \ref{lemma: property 2 of the estimators in the algorithm} controls the $\ell _{2}$ estimation error based on Lemma \ref{lemma: property 1 of the estimators in the algorithm}, even when Condition \eqref{eq: definition of J} fails for overly large $J$.

\begin{lemma} \label{lemma: property 1 of the estimators in the algorithm}
Let $\epsilon_J$ be defined as $\epsilon _{J}:=\frac{d}{2 ^{J-1}(C+1)}$ and $c \ge 2(C+1)$ possibly depending on $M$. Let $\tilde{J}$ be any positive integer such that 
\begin{equation}\label{eq: definition of J}
\kappa(M)\epsilon_J^2 > 2 \log [ N^{\text{loc}}(c\epsilon_{J}, 2c)]^2 \lor \log 2, 
\end{equation}
Then for $\tilde{J}$, we have
\begin{equation*}
    \mathbb{P}\left( \|\Upsilon _{\tilde{J}} - \theta\| > \frac{d}{2^{\tilde{J}-1}} \right) \le 2 \cdot \mathbf{1}_{\left\{\tilde{J}>1\right\}}\exp\left\{-\frac{\kappa (M)\epsilon _{\tilde{J}}^{2}}{2}\right\}.
\end{equation*}
\end{lemma}

\begin{lemma}\label{lemma: property 2 of the estimators in the algorithm}
Let $\epsilon_J$ be defined in the same way as in Lemma \ref{lemma: property 1 of the estimators in the algorithm} and $c \ge 2(C+1)$. Suppose $\tilde{J}$ is a positive integer such that (\ref{eq: definition of J}) holds and let $J ^{*}$ represent the maximal such $\tilde{J}$. Then if $\theta^{**}$ denotes the output after at least $J^*$ iterations, we have
\begin{equation*}
    \mathbb{E}_{Y}\|\theta - \theta^{**}(Y)\|^2 \le \frac{(16C+19)^{2}}{4}\epsilon_{J ^{*}}^2 + \mathbf{1}_{\left\{J ^{*}>1\right\}} \cdot (16C+19)^{2} \cdot \frac{1}{\kappa(M)}\exp\left\{-\frac{\kappa(M)\epsilon_{J ^{*}}^2}{2}\right\}.
\end{equation*}
\end{lemma}

The proofs of Lemma \ref{lemma: property 1 of the estimators in the algorithm} and Lemma \ref{lemma: property 2 of the estimators in the algorithm} are deferred to Appendix \ref{subsubsection: proof of the property 1} and Appendix \ref{subsubsection: proof of the property 2}, respectively.

We are now equipped to present the main theoretical guarantee for the estimator generated by Algorithm \ref{algorithm: theoretical algorithm}. In Theorem \ref{theorem: performance of the theoretical algorithm}, we establish a rigorous upper bound for its performance and consequently an upper bound for the minimax rate \eqref{eq: minimax risk of the estimation}. For conciseness, we place its proof in Appendix \ref{subsubsection: proof of the performance of the theoretical algorithm}.

\begin{theorem}\label{theorem: performance of the theoretical algorithm}
For any $\epsilon_J = \frac{d}{2^{J-1}(C+1)}, J \ge 1$ and the constant $c$ in Definition \ref{def_local_entropy} satisfying $c \ge 2(C + 1)$, let $J^*$ be the maximal $J$ such that \ref{eq: definition of J} holds (set $J^* = 1$ if this condition never occurs). Then if $\theta^{**}(Y_1, \dots, Y _{n})$ is the output of Algorithm \ref{algorithm: theoretical algorithm} after at least $J^*$ iterations, we have
\[
\mathbb{E}_Y \|\theta^{**}(Y _{1},\dots,Y _{n}) - \theta\|_{2}^2 \lesssim \epsilon_{J^*}^2 \wedge d^2,
\]
where ``$\lesssim $'' omits a constant only depending on $M$.
\end{theorem}

\begin{remark}
\label{RevisedRemark3.13}

By definition, the local metric entropy $N^{\text{loc}}(\epsilon, c)$ is computed using an $(\epsilon / c)$-covering of balls within $K$. In the lower bound established in Lemma \ref{lemma: minimax lower bound}, utilizing $\tilde{c} = 2c$ instead of $c$ leaves the theoretical rate unchanged up to an absolute constant. Consequently, the local metric entropy parameter in the lower bound can be chosen to perfectly match the $2c$ appearing in \eqref{eq: definition of J} of Theorem \ref{theorem: performance of the theoretical algorithm}. Without loss of generality, we can therefore assume that the same sufficiently large constant $c$ appears in both the upper and lower bounds, allowing us to replace \eqref{eq: definition of J} with the following criterion: let $J^*$ be the maximal integer such that
\begin{equation}\label{eq: alternative definition of J}
    \kappa(M)\epsilon_J^2 > 2 \log [ N^{\text{loc}}(c\epsilon_{J},c)]^2 \lor \log 2
\end{equation}
with $J^* = 1$ if this never occurs. We refer the interested reader to Lemma 1.4 of \cite{prasadan2024} for a similar argument and the proof.

\end{remark}

\subsection{Minimax Rate}
\label{SecMinimax}

In this section, we establish the minimax optimality—up to constants only depending on $M$—of the estimator produced by Algorithm \ref{algorithm: theoretical algorithm}. This primary result is formalized in Theorem \ref{theorem: main theorem for theoretical algorithm} below.

\begin{theorem}\label{theorem: main theorem for theoretical algorithm} For the natural parameter estimation problem specified in Section \ref{subsection: problem formulation} and $\epsilon^*$ defined as
    \begin{equation*}
        \epsilon ^{*}:=\sup \{\epsilon: \epsilon >0, \epsilon^2 \kappa(M) \le  \log N^{\text{loc}}(\epsilon,c)\},
    \end{equation*}
    where $c \ge 2(C+1)$ is a sufficiently large constant allowed to depend on $M$ in the definition of \hyperref[def_local_entropy]{local entropy} and $\kappa (M)$ is defined as 
    \begin{equation*}
        \begin{aligned}
            &\kappa (M):=\frac{c(M)^2 C^2}{32 C'(M)}\wedge \frac{|C(M) - c(M)(C-1)^2|}{4};\\ 
            &c(M):=\min\limits _{-M \le \theta \le M}A ^{\prime \prime }(\theta ),C(M):=\max\limits _{-M \le \theta \le M}A ^{\prime \prime }(\theta ),C ^{\prime }(M):=\max\limits _{-2M \le \theta \le 2M}A ^{\prime \prime }(\theta );\\ 
            &C \ge 2\sqrt{\frac{C(M)}{c(M)}+8} \vee \sqrt{\frac{2C(M)}{c(M)}}+1.
        \end{aligned}
    \end{equation*}
    Then the minimax rate \eqref{eq: minimax risk of the estimation} of the $\ell _{2}$ error for estimating $\theta _{0}$ is 
\begin{equation}\label{eq: minimax rate}
    \mathcal{R}_{n}(f,K)\asymp _{M}\epsilon ^{*2}\wedge \text{diam}(K) ^{2},
\end{equation}
where the constant depending only on $M$ is omitted. Moreover, the minimax rate \eqref{eq: minimax rate} is achieved by Algorithm \ref{algorithm: theoretical algorithm}.
\end{theorem}

\begin{remark*}[]
    Theorem \ref{theorem: main theorem for theoretical algorithm} establishes a minimax rate that depends on $M$, the geometric parameter bounding the constraint set $K$. While this dependence is generally unavoidable within our overarching framework, it is instructive to consider the standard Gaussian sequence model ($\sigma = 1$) as a special case. Under this distribution, the log-partition function is $A(\theta) = \theta^2/2$, which yields the absolute constants $c(M) = C(M) = C'(M) = 1$. Consequently, the minimax rate becomes entirely independent of $M$. This seamlessly recovers the optimal rate previously established by \cite{neykov2023minimax}.
\end{remark*}

The proof of Theorem \ref{theorem: main theorem for theoretical algorithm} crucially leverages the foundations established in Lemma \ref{lemma: minimax lower bound} and Theorem \ref{theorem: performance of the theoretical algorithm}. For clarity of presentation, the detailed argument is deferred to Appendix \ref{subsubsection: proof of main theorem for theoretical algorithm}.

\section{Efficient Estimator}\label{section: efficient estimator}

Although Algorithm \ref{algorithm: theoretical algorithm} achieves the optimal minimax rate as established in Theorem \ref{theorem: main theorem for theoretical algorithm}, it remains primarily of theoretical interest. Specifically, for an arbitrary boxed, star-shaped constraint $K$, constructing the packing set for $B_{2}\left(u,\frac{d}{2^{k-2}}\right)\cap K$ at level $k$ ($k \ge 2$) is computationally nontrivial. Furthermore, even if this construction were feasible, the cardinality of such a packing set grows exponentially with the ambient dimension $n$.

To address this computational bottleneck, we propose an efficient algorithm in this section that achieves near-minimax optimal performance. The proposed method operates in polynomial time with respect to the dimension $n$ and the geometric parameters of the constraint $K$, while inflating the minimax $\ell_{2}$ error of the resulting estimator by at most a poly-logarithmic factor. Before detailing the algorithm and establishing its theoretical guarantees, we first introduce additional structural conditions on $K$ and review necessary preliminaries. Similarly, the proofs of some lemmas, theorems, and corollaries are deferred to Appendix \ref{subsection: proofs of the efficient estimator}.

\subsection{Preliminary Properties and Concepts}\label{subsection: preliminary properties and concepts}

First, we assume that the log-partition function $A(\cdot)$ is defined on the entire real line $\mathbb{R}$ and is twice continuously differentiable with a uniformly bounded second derivative. As established in Lemma \ref{lemma: sub-Gaussian of T} below, whose result is stronger than Lemma \ref{lemma: sub-exponential of T}, this structural condition guarantees that the sufficient statistic $T(Y_i)$ is a sub-Gaussian (instead of sub-exponential) random variable for any $i \in [n]$.

\begin{lemma}[Sub-Gaussian property of $T$]\label{lemma: sub-Gaussian of T}
    Under the previous assumptions on the log-partition function $A(\cdot)$, the sufficient statistic $T(Y_i)$ is a sub-Gaussian random variable for all $i \in [n]$.
\end{lemma}

\begin{proof}
    Since the argument holds identically for all $i \in [n]$, we drop the index $i$ for simplicity. Following a similar approach to the proof of Lemma \ref{lemma: sub-exponential of T}, we evaluate the centered log-moment generating function of $T(Y)$. By Taylor's theorem, we have
    \begin{equation*}
        \log \left(\mathbb{E}\exp\left\{\lambda \left[T(Y)-\mathbb{E}T(Y)\right]\right\}\right)=A(\theta +\lambda )-A(\theta )-\lambda \mathbb{E}\left[T(Y)\right]=\frac{1}{2}\lambda^{2}A^{\prime\prime}(\tilde{\theta}),
    \end{equation*}
    for some $\tilde{\theta}$ lying on the line segment between $\theta$ and $\theta +\lambda$. By assumption, $A^{\prime\prime}$ is uniformly bounded on $\mathbb{R}$. Therefore, the centered log-moment generating function is bounded by $\frac{\lambda ^{2}C _{M}^{2}}{2}$ for some universal constant $C _{M}>0$, which formally establishes that $T(Y)$ is sub-Gaussian.
\end{proof}

\begin{remark*}[]
    We note that multiple distributions in the exponential family satisfy such condition, such as Gaussian distribution, (continuous) Bernoulli distribution, multinomial distribution.
\end{remark*}

Next, we impose stronger structural assumptions on the constraint set $K$ beyond the previous star-shaped condition $K \subset [-M,M]^{n}$. Specifically, we assume that $K$ is origin-symmetric (i.e., $\theta \in K \implies -\theta \in K$), convex, and well-balanced. Furthermore, we require its Minkowski gauge to be sign-invariant and exactly 2-convex with a controlled type-2 constant. We formally introduce these geometric concepts below.

\begin{definition}[Well-balanced set]\label{def_well_balanced_set}
  For a subset $K \subset \mathcal{X}=\mathbb{R}^{n}$, $K$ is well-balanced if there exist $0<r \le R \le \infty $ such that $B _{2}(0,r)\subset K \subset B _{2}(0,R)$.
\end{definition}

\begin{remark*}[]
    Under our current settings, we know $R \le M \sqrt{n}<\infty $ since $K \subset [-M,M]^{n}$. We shall also assume in the following that $r,R$ are known constants.
\end{remark*}

\begin{definition}[Minkowski gauge]\label{def_minkowski_gauge}
  Let $\mathcal{X}$ be a vector space and $K \subset \mathcal{X}$ be a fixed subset. The Minkowski gauge (or Minkowski functional) of $K$ is a function $\rho_{K}: \mathcal{X} \rightarrow [0,\infty]$ defined as
  \begin{equation}
    \rho_{K}(\theta) = \inf\limits _{}\left\{r \left\lvert\right. r>0, \theta \in rK\right\}.
  \end{equation}
  By convention, if no such $r$ exists for a given $\theta$, we set $\rho_{K}(\theta) = \infty$.
\end{definition}

\begin{remark*}
  Under our current settings for $K$ (bounded, convex, and origin-symmetric), $\rho _{K}$ is a norm, and $K$ could be equivalently defined as $K=\left\{\theta \left\lvert\right. \theta \in \mathcal{X},\rho _{K}(\theta )\le 1\right\}$. In the subsequent analysis, we assume access to an evaluation oracle for the Minkowski gauge $\rho_K$. Given that $K$ is star-shaped containing the origin, this is essentially equivalent to having a membership oracle for the constraint $K$ (up to arbitrary precision via binary search), which is a standard and mild assumption in the algorithmic literature.
\end{remark*}

\begin{definition}[Sign-invariant norm]\label{def_sign-invariant_norm}
  A norm $\lVert \cdot \rVert$ on $\mathbb{R}^n$ is said to be \textit{sign-invariant} if for any $\theta \in \mathbb{R}^{n}$ and any sign vector $\gamma \in \{-1,1\}^n$, it holds that $\lVert \theta \rVert = \lVert \gamma \odot \theta \rVert$, where $\odot$ denotes the Hadamard (element-wise) product.
\end{definition}

\begin{definition}[Exactly 2-convexity]\label{def_exactly_2_convex_norm}
The Minkowski gauge $\rho_{K}(\cdot)$ of a set $K$ is said to be 2-convex if there exists a constant $c(K)$ depending only on $K$ such that for any integer $m \ge 1$ and any vectors $\theta_{1}, \dots, \theta_{m} \in \mathbb{R}^{n}$, we have
\begin{equation}\label{eq: def_exactly_2_convex_norm}
\rho_{K}\left(\left( \sum_{i=1}^{m} \theta_{i}^{2}\right)^{1/2} \right) \le c(K) \left(\sum_{i=1}^{m} \rho_{K}^{2}(\theta_{i}) \right)^{1/2},
\end{equation}
where $\theta_{i}^{2} = (\theta_{i1}^{2}, \dots, \theta_{in}^{2})^{\top} \in \mathbb{R}^{n}$ denotes the entrywise square of the vector $\theta_{i}$. When $c(K) = 1$, the gauge $\rho_{K}(\cdot)$ is said to be \textit{exactly 2-convex}.
\end{definition}

\begin{remark*}[]
    As an example, all $\ell _{p}$ norm with $2 \le p \le \infty $ is exactly $2$-convex norm.
\end{remark*}

\begin{definition}[Type-2 condition]\label{def_type_2_condition}
  Let $\left\{\epsilon _{i}\right\}_{i=1}^{m}$ be i.i.d. Rademacher random variables for some $m \in \mathbb{N}^{+}$. For $K \subset \mathcal{X}$ above, $K$ is type-$2$ with constant $T _{2}(K)$ if for any $m \in \mathbb{N}^{+}$ and any $\theta _{i}\in \mathbb{R}^{n}$, we have 
  \begin{equation}\label{eq: required condition in type 2 definition}
    \mathbb{E}_{\epsilon }\rho _{K}^{2}\left(\sum\limits_{i=1}^{m}\epsilon _{i}\theta _{i}\right)\le T _{2}^{2}(K)\sum\limits_{i=1}^{m}\rho _{K}^{2}(\theta _{i}).
  \end{equation}
  $T _{2}(K)$ is referred as the type-$2$ constant of $K$. $T _{2}(K)$ is said to be controlled if $T _{2}(K)=\text{polylog}(n)$.
\end{definition}

\begin{remark*}
  By selecting $m=1$, it is straightforward to see that $T _{2}(K)\ge 1, \forall K \subset \mathcal{X}$. We note that the well-known quadratically convex orthosymmetric (QCO) sets, extensively studied in the literature \cite{10.1214/aos/1176347758,li2026robustsignaldetectionquadratically}, satisfy all the aforementioned conditions. This class encompasses a wide range of geometries, including axis-aligned hyperrectangles, ellipsoids, and sets of the form $\left\{ \theta \in \mathbb{R}^n \mid \sum_{i=1}^n a_i \psi(\theta_i^2) \le 1 \right\}$ for some convex function $\psi$ and constants $a_i > 0$. We refer interested readers to Section 6 of \cite{neykov2026polynomialtimenearoptimalestimationcertain} for a formal proof of these geometric properties.
\end{remark*}

\subsection{Quadratic Form Maximization and Kolmogorov Width Approximation}\label{subsection: quadratic form maximization and kolmogorov width approximation}

The geometric properties introduced in the last section facilitate the existence of an approximate separation oracle \cite{grotschel1993geometric,Dadush2012IntegerPL} for the upper covariance region of $K$. Formally, this region is defined as
\begin{equation*}
  \mathcal{U}(K) := \left\{ W \in \mathcal{S}_{+}^{d} \mid \mathbb{E}\left[\rho_{K}^{2}\left(W^{1/2}\bm{g}\right)\right] \le 1 \right\},
\end{equation*}
where $\mathcal{S}_{+}^{n}$ denotes the cone of $n \times n$ positive semi-definite matrices, and $\bm{g} \sim \mathcal{N}(0, \mathbf{I}_{n})$ is a standard Gaussian vector. The existence of such an oracle subsequently enables an approximate solution to the quadratic form maximization problem, a framework established by \cite{10.1145/3406325.3451128} which we briefly review below.

\begin{theorem}[Proposition 4.2 and Theorem 7.6 of \cite{10.1145/3406325.3451128}]\label{theorem: main theorem of the cited paper}
  For any symmetric matrix $X \in \mathbb{R}^{n \times n}$, consider the quadratic maximization problem:
  \begin{equation*}
    \begin{aligned}
      & \text{\rm maximize} && \theta^{\top} X \theta \\ 
      & \text{\rm subject to} && \theta \in K.
    \end{aligned}  
  \end{equation*}
  Suppose $K$ is origin-symmetric, convex, and well-balanced with parameters $r$ and $R$. Further, assume its Minkowski gauge $\rho_{K}$ is sign-invariant and exactly 2-convex with a controlled type-2 constant. Then, there exists a randomized algorithm $\mathcal{A}(X, R, r)$ that, for any such input $X$, outputs a vector $\mathcal{O}_{K}(X) \in K$ in polynomial time in $n$, $\log \left(\frac{R}{r}\right)$, the bit complexity of $X$, and $\log \left(\frac{1}{q}\right)$. Furthermore, with probability at least $1-q$, the output satisfies
  \begin{equation*}
    \mathcal{O}_{K}(X)^{\top} X \mathcal{O}_{K}(X) \gtrsim \frac{1}{\tau _{K}} \max_{\theta \in K} \theta^{\top} X \theta
  \end{equation*}
  for some $\tau _{K}=\operatorname{polylog}(n)$ and can depend on $K$.
\end{theorem}

The specific choice of the probability parameter $q$ is deferred to subsequent discussions. Before stating a key result from \cite{neykov2026polynomialtimenearoptimalestimationcertain}, we first review the fundamental concept of the (approximate) Kolmogorov $k$-width \citep{2ad9b75f-e5d2-3a8c-8211-9d37aa110eb2,neykov2026polynomialtimenearoptimalestimationcertain, li2026efficientrobustconstrainedsignal}. As a crucial geometric preliminary, this notion quantifies the minimax error of approximating a given set by $k$-dimensional linear subspaces under a specified norm.

\begin{definition}[Kolmogorov $k$-width]\label{def_kolmogorov_width_raw}
    Let $\mathcal{X}$ be a Banach space equipped with the norm $\left\lVert \cdot \right\rVert_{}^{}$, and $K \subset \mathcal{X}$ is any subset. The $k$-dimensional Kolmogorov width of $K$ is defined as 
    \begin{equation}\label{eq: raw definition of Kolmogorov width}
        D _{k}(K) = \inf_{V \in \mathcal{V}_{k}} \sup_{\theta \in K} \inf\limits _{\tilde{\theta }\in V}\left\lVert \theta -\tilde{\theta } \right\rVert_{}^{},
    \end{equation}
    where $\mathcal{V}_{k}$ is the set of all $k$-dimensional subspaces of $\mathcal{X}$.
\end{definition}

\begin{remark*}
  Note that the Kolmogorov $k$-width is non-increasing in $k$. Furthermore, when the ambient space $\mathcal{X}$ is $n$-dimensional, we trivially have $D_{0}(K) = \sup\limits _{\theta \in K} \lVert \theta \rVert$ and $D_{n}(K) = 0$. Throughout this paper, unless otherwise specified, we restrict our attention to the Euclidean norm and set $\mathcal{X} = \mathbb{R}^{n}$. Under this setting, the best approximation of $\theta$ within a subspace $V$ is exactly its orthogonal projection $P_{V}\theta$. Consequently, the original definition in \eqref{eq: raw definition of Kolmogorov width} can be equivalently formulated as
  \begin{equation}\label{eq: definition of Kolmogorov width}
    D_{k}(K) = \inf_{P \in \mathcal{P}_{k}} \sup_{\theta \in K} \lVert \theta - P\theta \rVert,
  \end{equation}
  where $\mathcal{P}_{k}$ denotes the set of all rank-$k$ orthogonal projection matrices in $\mathbb{R}^{n \times n}$.
\end{remark*}

The Kolmogorov $k$-width is a fundamental concept in approximation theory with extensive applications; we refer readers to \cite{pinkus2012n} for a comprehensive overview. Furthermore, this concept has been successfully deployed across diverse fields, including matrix theory, neural networks, optimal transport, and partial differential equations (PDEs) \cite{Floater_2021,EVANS20091726,10.1093/imanum/dru066,GREIF2019216,PAPAPICCO2022114687,Arbes2025}.

Despite its theoretical significance, calculating the exact Kolmogorov $k$-width is computationally intractable, even for highly structured sets. For instance, \cite{brieden2002geometric} rigorously establishes the computational hardness of this problem for convex polytopes. To circumvent this computational bottleneck, \cite{neykov2026polynomialtimenearoptimalestimationcertain, li2026efficientrobustconstrainedsignal} introduced the following convex minimax relaxation to approximate the Kolmogorov $k$-width.

\begin{definition}[Approximate Kolmogorov $k$-Width]\label{def_approximate_kolmogorov_width}
Let the ambient space be $\mathbb{R}^{n}$ equipped with the Euclidean norm. For any subset $K \subset \mathbb{R}^{n}$ and any integer $0 \le k \le n$, we define the approximate Kolmogorov $k$-width, denoted by $\tilde{D}_{k}^{2}(K)$, as the optimal value of the following convex minimax problem:
\begin{equation}\label{eq: SDP problem}
    \begin{aligned}
        & \min_{X \in \mathcal{S}^{n}} \max_{\theta \in K} \theta^{\top} X \theta \\
        & \text{\rm subject to} \quad \operatorname{tr}(X) = n-k, \quad \bm{0} \preceq X \preceq \mathbf{I}_{n},
\end{aligned}
\end{equation}
where $\mathcal{S}^{n}$ denotes the set of $n \times n$ symmetric matrices, and $\mathbf{I}_{n}$ is the identity matrix. We denote a corresponding optimal solution to \eqref{eq: SDP problem} by $X^{\star}_{k}$.
\end{definition}

The formulation in \eqref{eq: SDP problem} is indeed a valid relaxation of the Kolmogorov $k$-width, in the sense that $\tilde{D}_{k}^{2}(K) \le D_{k}^{2}(K)$ for all $0 \le k \le n$. This holds because the matrix $\mathbf{I}_{n} - P_{V}$ naturally satisfies the constraints of \eqref{eq: SDP problem} for any projection operator $P_{V}$ onto a $k$-dimensional subspace $V \in \mathcal{V}_{k}$. Furthermore, by definition, we trivially have $\theta^{\top} X_{k}^{\star} \theta \le \tilde{D}_{k}^{2}(K)$ for any $\theta \in K$.

We note that solving \eqref{eq: SDP problem} exactly remains computationally intractable in general, meaning the exact optimal solution $X_{k}^{\star}$ is practically inaccessible. However, for the inner maximization problem over $\theta \in K$, a direct application of Theorem \ref{theorem: main theorem of the cited paper} provides an approximate oracle $\mathcal{O}_{K}(X) \in K$. Specifically, for any feasible $X$, it guarantees that $\mathcal{O}_{K}(X)^{\top} X \mathcal{O}_{K}(X) \ge \frac{1}{\tau _{K}} \max_{\theta \in K} \theta^{\top} X \theta$ with high probability. Building upon this approximate oracle, we can leverage either subgradient descent (Algorithm 1 of \cite{neykov2026polynomialtimenearoptimalestimationcertain}) or the modified ellipsoid method (Algorithm 3 of \cite{li2026efficientrobustconstrainedsignal}) to efficiently compute an approximate solution to the overall minimax problem \eqref{eq: SDP problem}. As established by Theorems 2.5 and 3.1 of \cite{neykov2026polynomialtimenearoptimalestimationcertain}, this procedure yields a solution accurate up to a factor of $\tau _{K}$ with probability at least $1-q$ with $q=\frac{1}{1+\tilde{d} ^{3}},\tilde{d}=2R$.

We are now ready to state a crucial result from \cite{neykov2026polynomialtimenearoptimalestimationcertain}. This lemma establishes an upper bound on the Kolmogorov $k$-width in terms of $\tilde{d}$, $\tau _{K}$, and $c$ for a specifically chosen dimension $k$. Consequently, it guarantees that we can bound the approximation error of the minimax problem \eqref{eq: SDP problem} using a subspace whose dimension is bounded by $\mathcal{O}(\tilde{d}^{2}/C^{2})$.

\begin{lemma}[Lemma 2.1 of \cite{neykov2026polynomialtimenearoptimalestimationcertain}]\label{lemma: kolmogorov width control}
Let $K$ be a type-2 convex set with type-2 constant $T_2(K)$, such that $B _{2}(0,r) \subseteq K \subseteq B _{2}(0,R)$ for some known $r,R$. Let $\tau _{K}$ be a given scalar which can depend on $K$. Assume that $\frac{\tilde{d}}{\sqrt{\tau _{K}}T _{2}(K)}\ge 2\epsilon ^{*}$ and $\lceil 4\lceil \log_c(\sqrt{\tau _{K}}T _2(K))\rceil \epsilon ^{*2}\rceil \le \frac{\tilde{d}^{2}}{C ^{2}}$ for $\tilde d = 2R \ge d$ and some sufficiently large universal constant $C$. Then for $\tilde{k} = \lceil \frac{\tilde{d}^{2}}{C ^{2}} \rceil \wedge n$ we have
\begin{equation}
    D _{\tilde{k}}(K)\lesssim \frac{\tilde{d}}{\sqrt{\tau _{K}}c}.
\end{equation}
\end{lemma}

The following theorem is directly from Theorem \ref{theorem: main theorem of the cited paper}, Lemma \ref{lemma: kolmogorov width control}, and the discussions above.

\begin{theorem}[]\label{theorem: polynomial time approximation to the SDP problem}
    Given a symmetric, convex, well-balanced set $K$ with its Minkowski gauge $\rho _{K}$ being a sign-invariant and exact 2-convex with controlled type-$2$ constant, suppose that $\frac{\tilde{d}}{\sqrt{\tau _{K}}T _{2}(K)}\ge 2\epsilon ^{*}$ and $\lceil 4\lceil \log_c(\sqrt{\tau _{K}}T _2(K))\rceil \epsilon ^{*}\rceil \le \frac{\tilde{d}^{2}}{C ^{2}}$ for $\tilde d = 2R \ge d$, then we are able to find a $X ^{\dagger }$ that is feasible to \eqref{eq: SDP problem} with $k=\tilde{k}=\lceil \frac{\tilde{d}^{2}}{C ^{2}} \rceil \wedge n$ such that 
    \begin{equation}\label{eq: property of the approximation oracle}
        \sup\limits _{\theta \in K}\theta ^{\top } X ^{\dagger }\theta \lesssim \frac{\tilde{d}^{2}}{c ^{2}}
    \end{equation}
    with probability greater than $1-\frac{1}{1+\tilde{d}^{3}}$. In addition, all the procedures finish with in polynomial time of $n, \log \left(\frac{R}{r}\right), \log \tilde{d}$.
\end{theorem}

\begin{remark*}
  We note that the structural assumptions on $K$ in Theorem \ref{theorem: polynomial time approximation to the SDP problem} can actually be significantly relaxed. Specifically, one only requires that $K$ is well-balanced and can be represented in the form $K = \{x \in \mathbb{R}^{n} \mid \lVert Ax \rVert \le 1\}$, where $A \in \mathbb{R}^{m \times n}, m \ge n$ is a known matrix and $\lVert \cdot \rVert$ is an exactly $2$-convex, sign-invariant norm in $\mathbb{R}^{m}$ with a controlled type-$2$ constant. This generalized result is formally established in Theorems 3.1 and 3.3 of \cite{neykov2026polynomialtimenearoptimalestimationcertain}.
\end{remark*}

\subsection{Weak Projection and the Properties}\label{subsection: weak projection and the properties}

The final algorithmic building block addresses the natural constraints of the problem, which enables us to weakly project an unconstrained estimator back into the constraint set $K$ with arbitrary accuracy. Specifically, we review the following projection results established by \cite{Dadush2012IntegerPL}.

\begin{lemma}[Theorem 2.5.9 of \cite{Dadush2012IntegerPL}]\label{lemma: raw approximation from dadush}
    Suppose that $K \subset \mathbb{R}^{n}$ such that $B _{2}(\theta ,r)\subset K \subset B _{2}(\theta ,R)$ for some $\theta \in \mathbb{R}^{n}$ and $r,R>0$ and there is a weak membership oracle $O _{K}$ for $K$. If $f:\mathbb{R}^{n}\rightarrow \mathbb{R}$ is an $L$-Lipschitz convex function that can be (efficiently) evaluated, then for any $\epsilon >0$, a rational number $\omega \in \mathbb{Q}$ and a vector $\theta _{O _{K}}\in K$ can be computed using $O _{K}$ in polynomial time such that 
    \begin{equation*}
        \omega -\epsilon \le \min\limits _{\theta \in K}f(\theta )\le f(\theta _{O _{K}})\le \omega .
    \end{equation*}
\end{lemma}

Taking $f=\left\lVert \cdot \right\rVert_{2}^{}$ in Lemma \ref{lemma: raw approximation from dadush} and noticing that $\left\lVert \cdot \right\rVert_{2}^{}$ is $1$-Lipschitz, we obtain the following existence of the weak projection of $K$:

\begin{theorem}[existence of the weak projection]\label{theorem: existence of the weak projection}
    For the constraint $K$ in the current problem and any $\theta \in \mathbb{R}^{n},\epsilon >0$, there is rational number $\omega \in \mathbb{Q}$ and a weak projection $\tilde{P} _{K}(\cdot )$ satisfying 
    \begin{equation}\label{eq: property of the weak projection}
        \omega -\epsilon \le \min\limits _{\tilde{\theta }\in K}\left\lVert \theta -\tilde{\theta }\right\rVert_{2}^{}\le \left\lVert \theta -\tilde{P}_{K}(\theta )\right\rVert_{2}^{}\le \omega .
    \end{equation}
    Furthermore, all the computation finishes within polynomial-time complexity of $n, \log \left(\frac{R}{r}\right), \log \left(\frac{1}{\epsilon }\right)$.
\end{theorem}

The following corollary from \cite{neykov2026polynomialtimenearoptimalestimationcertain} bounds the distance between any arbitrary point $\xi \in K$ and the projected point $\tilde{P}_{K}(\theta)$ in terms of $\epsilon$ and the distance between $\xi$ and the unprojected point $\theta$. The proof follows directly from elementary properties of convex sets and the triangle inequality. We defer the proof to Appendix \ref{subsubsection: proof of the corollary related to the weak projection} for the reader's reference.

\begin{corollary}[]\label{corollary: control of the distance between the weak projection and any point in K}
    For any $\xi \in K$, we have $\left\lVert \tilde{P}_{K}(\theta )-\xi \right\rVert_{2}^{}\le \left\lVert \theta -\xi \right\rVert_{2}^{}+\sqrt{\epsilon ^{2}+2\epsilon \left\lVert \theta -\xi \right\rVert_{2}^{}}$.
\end{corollary}

\begin{remark*}[]
    Though $\epsilon $ is arbitrary, we cannot deduce $\left\lVert \tilde{P}_{K}(\theta )-\xi \right\rVert_{2}^{}\le \left\lVert \theta -\xi \right\rVert_{2}^{}$ simply by letting $\epsilon \rightarrow 0$ in Corollary \ref{corollary: control of the distance between the weak projection and any point in K}, as the determination of $\tilde{P}_{K}(\theta )$ in Theorem \ref{theorem: existence of the weak projection} depends on $\epsilon $.
\end{remark*}

\subsection{Efficient Algorithm Illustration}\label{subsection: efficient algorithm illustration}

We now formally present Algorithm \ref{algorithm: main algorithm}, a computationally efficient procedure that yields an estimator achieving a near-optimal minimax $\ell_2$ error. The complete algorithm is listed in Algorithm \ref{algorithm: main algorithm}, and the detailed description of the algorithm and some relevant motivation is as follows.

The procedure fundamentally depends on the value of $r$. The regime where $r \ge \sqrt{n}$ is trivial by Lemma \ref{lemma: relationship between minimax rate and r}. For the regime where $r \le \sqrt{n}$, we initialize the procedure with the constraint set $K$. In the $j$-th iteration, we define a localized constraint $K_{(j)} := K \cap B_{2}(0, \tilde{d}_{(j)}/2)$. We then compute the target dimension $\tilde{k}_{(j)}$ as specified in Lemma \ref{lemma: kolmogorov width control}, and obtain $X_{(j)}^{\dagger}$ from the convex minimax relaxation problem \eqref{eq: SDP problem} via Theorem \ref{theorem: polynomial time approximation to the SDP problem}. Next, we define the search space $E_{(j)} := \big\{ \hat{\theta}_{(j)} + A_{(j)}^{\dagger}(\theta - \hat{\theta}_{(j)}) \mid \theta \in K \cap (\hat{\theta}_{(j)} + 2K_{(j)}) \big\}$, where $A_{(j)}^{\dagger} := (\mathbf{I} - X_{(j)}^{\dagger})^{1/2}$. Over this set $E_{(j)}$ (which is naturally parameterized by $\theta \in K \cap (\hat{\theta}_{(j)} + 2K_{(j)})$ by definition of $E _{(j)}$), we optimize the log-likelihood $\sum\limits _{i=1}^{n}\theta_{i} T(y_{i}) - A(\theta_{i})$ via the ellipsoid method to obtain $\hat{\theta}_{\text{m}}^{(j+1)}$. For the unconstrained intermediate estimator $\hat{\theta}_{\text{r}}^{(j+1)} := \hat{\theta}_{(j)} + A_{(j)}^{\dagger}\big(\hat{\theta}_{\text{m}}^{(j+1)} - \hat{\theta}_{(j)}\big)$, we project it back onto $K$ using the weak projection oracle from Theorem \ref{theorem: existence of the weak projection}. This yields the next estimator $\hat{\theta}_{(j+1)} = \tilde{P}_{K}(\hat{\theta}_{\text{r}}^{(j+1)})$. We then update the distance bound as $\tilde{d}_{(j+1)} = \sqrt{\frac{8C_{A}}{c_{A}}} \frac{(\sqrt{3}+1)\tilde{d}_{(j)}}{c}$. The loop terminates once $\tilde{d}_{(j)} \le r \vee c ^{2}$, returning the final estimator $\hat{\theta} = \hat{\theta}_{(j)}$.

The core of the algorithm lies in the definition of $E_{(j)}$. To provide some high-level intuition, suppose that at the $j$-th iteration, the $\ell_{2}$ distance between the current estimator $\hat{\theta}_{(j)}$ and the true parameter $\theta_{0}$ is bounded by $\tilde{d}_{(j)}$, i.e., $\lVert \hat{\theta}_{(j)} - \theta_{0} \rVert_{2} \le \tilde{d}_{(j)}$ (which is naturally satisfied by the first iteration). This implies $\theta_{0} - \hat{\theta}_{(j)} \in B_{2}(0, \tilde{d}_{(j)})$. Furthermore, since $\theta_{0}, \hat{\theta}_{(j)} \in K$, we have $\theta_{0} - \hat{\theta}_{(j)} \in 2K$. Consequently, $\theta_{0} \in \hat{\theta}_{(j)} + 2K_{(j)}$ by the definition of $K_{(j)}$. To further narrow down the search space, we apply the linear operator $A_{(j)}^{\dagger}$. This transformation ensures that the ``proxy true parameter'' $\hat{\theta}_{(j)} + A_{(j)}^{\dagger}(\theta_{0} - \hat{\theta}_{(j)})$ is strictly closer to $\theta_{0}$ than the previous iterate $\hat{\theta}_{(j)}$ given the valid $X _{(j)}^{\dagger }$ and $A _{(j)}^{\dagger }$, because its squared $\ell_2$ distance satisfies $\lVert \theta_{0} - \big(\hat{\theta}_{(j)} + A_{(j)}^{\dagger}(\theta_{0} - \hat{\theta}_{(j)})\big) \rVert_{2}^{2} = (\theta_{0} - \hat{\theta}_{(j)})^{\top} X_{(j)}^{\dagger} (\theta_{0} - \hat{\theta}_{(j)}) \le \tilde{d}_{(j)}^{2} / c^{2}$ by Theorem \ref{theorem: polynomial time approximation to the SDP problem}. Meanwhile, as established by Theorem \ref{theorem: exponentially decreasing distance between the estimator and true parameter} below, the distance between $\hat{\theta }_{\text{m}}^{(j+1)}$ and $\hat{\theta }_{(j)}+A _{(j)}^{\dagger }(\theta _{0}-\hat{\theta }_{(j)})$ is also bounded by $\tilde{d}_{(j)}^{2}/c ^{2}$ up to some universal constant. Therefore, after at most a polynomial number of iterations, this exponential shrinkage yields a highly accurate estimator.

\begin{algorithm}[htbp]
    \caption{construction of the estimator $\hat{\theta }$}\label{algorithm: main algorithm}

    \KwInput{$K,r,R,y _{1},\dots,y _{n}$}

    \KwOutput{The estimator $\hat{\theta }$ of $\theta $}

    \If{$r \gtrsim \sqrt{n}$}{
        \Return arbitrary $\hat{\theta }\in K$
    }

    $\tilde{d}_{(0)}\gets R$

    $\hat{\theta }_{(0)}\gets \bm{0},j\gets 0$

    \While{$\tilde{d}_{(j)}\gtrsim r \vee c ^{2}$}{
        $K _{(j)}\gets K \cap B _{2}\left(0,\frac{\tilde{d}_{(j)}}{2}\right)$

        $\tilde{k} _{(j)}\gets \lceil \frac{\tilde{d}_{(j)} ^{2}}{C ^{2}} \rceil \wedge n$ as in Lemma \ref{lemma: kolmogorov width control}

        Obtain $X ^{\dagger }_{(j)}$ for $K=K _{(j)}$ by solving the convex minimax relaxation \eqref{eq: SDP problem} by Theorem \ref{theorem: polynomial time approximation to the SDP problem}

        $A ^{\dagger }_{(j)}\gets \left(\mathbf{I}_{n}-X ^{\dagger }_{(j)}\right)^{\frac{1}{2}}$

        $E _{(j)}\gets \left\{\hat{\theta } _{(j)}+A ^{\dagger }_{(j)}(\theta -\hat{\theta } _{(j)})\left\lvert\right. \theta \in K \cap \left(\hat{\theta } _{(j)}+2K _{(j)}\right)\right\}$

        Obtain $\hat{\theta } _{\text{m}}^{(j+1)}$ by maximizing the log-likelihood function $\sum\limits _{i=1}^{n}\left[\theta _{i}T(y _{i})-A(\theta _{i})\right]$ over $E _{(j)}$ parameterized by $\hat{\theta }_{(j)}+2K _{(j)}$ via the ellipsoid method with accuracy $\sqrt{\frac{2C _{M}}{c _{M}}}\frac{\tilde{d}_{(j)}}{c}$ (see the following remark)

        $\hat{\theta }_{(j+1)}\gets \tilde{P}_{K}\left[\hat{\theta }_{(j)}+A ^{\dagger }_{(j)}\left(\hat{\theta } _{\text{m}}^{(j+1)}-\hat{\theta }_{(j)}\right)\right]$

        $\tilde{d} _{(j+1)}\gets \sqrt{\frac{8C _{M}}{c _{M}}}\frac{\left(\sqrt{3}+1\right)\tilde{d}_{(j)}}{c}$

        $j\gets j+1$
    }

    \Return $\hat{\theta }_{\text{p}}=\hat{\theta }_{(j)}$

    \tcc{
        \begin{remark*}[]
            (1), we have a separation oracle for $K \cap \hat{\theta }_{(j)}+2K _{(j)}$ since we have the access of the Minkowski gauge of the constraint $K$ and trivially $B _{2}\left(0,\frac{\tilde{d}_{(j)}}{2}\right)$. By Theorem 2.5.9 of \cite{Dadush2012IntegerPL}, such ellipsoid method is feasible. (2), we set the accuracy of the ellipsoid method as $\sqrt{\frac{2C _{M}}{c _{M}}}\frac{\tilde{d}_{(j)}}{c}$. Such setting does not impact the efficiency of the algorithm since the number of iteration of the ellipsoid method is $\mathcal{O}\left(\log \left(\frac{\tilde{d} _{(j)}}{c}\right)\right)$. additionally, since $K$ is bounded and consequently $\sum\limits _{i=1}^{n}\theta _{i}T(y _{i})-A(\theta _{i})$ is Lipschitz, we therefore have $\left\lVert \hat{\theta }_{\text{m}}^{(j+1)}-\theta _{\text{m}}^{(j+1)}\right\rVert_{2}^{}\lesssim \frac{\tilde{d}_{(j)}}{c}$ as well, where $\theta _{\text{m}}^{(j+1)}$ is the true maximizer of the log-likelihood in the $j$-th iteration.
        \end{remark*}  
    }
    
\end{algorithm}

\subsection{Performance Analysis}\label{subsection: performance analysis}

This section is devoted to formally proving that the estimator produced by Algorithm \ref{algorithm: main algorithm} achieves the minimax rate of $(\epsilon^*)^2 \wedge \operatorname{diam}(K)^2$, up to poly-logarithmic factors. To establish this result, we first introduce the concept of sub-Gaussian stochastic processes \cite{10.1214/aop/1065725175}, which serves as a natural extension of sub-Gaussian random variables to the setting of stochastic processes.

\begin{definition}[sub-Gaussian tail stochastic process]
    For a stochastic process $\left(R _{t}\right)_{t \in T}$, we say that it has sub-Gaussian tail with respect to some metric $d _{T}(\cdot ,\cdot )$ defined on $T \times T$ if for any $t _{1},t _{2}\in T$ and $u \ge 0$, we have 
    \begin{equation*}
        \mathbb{P}\left(\left|R _{t _{1}}-R _{t _{2}}\right|\ge \sqrt{u}d _{T}(t _{1},t _{2})\right)\le 2e ^{-u}.
    \end{equation*}
\end{definition}

For stochastic processes with sub-Gaussian tails, we leverage a fundamental property that bounds the expected supremum of the process—and consequently its tail probability—using the diameter of the index set $T$ and the well-known $\gamma_2$-functional (\cite{10.1214/aop/1065725175,talagrand2005generic}). Such bounds on the expectation and tail probabilities were originally established by \cite{10.1214/aop/1008956336}. Subsequently, \cite{10.1214/EJP.v20-3760} sharpened these results and provided a simplified proof. We note that sub-Gaussian processes are a special case of a broader class of stochastic processes with mixed tails (\cite{10.1214/aop/1008956336}), though we omit those details here. For a comprehensive treatment of sub-Gaussian and mixed-tail stochastic processes, we refer the reader to \cite{10.1093/acprof:oso/9780199535255.001.0001,10.1214/EJP.v20-3760,talagrand2021upper}.

\begin{lemma}[Theorem 3.5 of \cite{10.1214/EJP.v20-3760}]\label{lemma: moment and tail bounds of sp}
    If the stochastic process $\left\{R _{t}\right\}_{t \in T}$ has a sub-Gaussian tail with respect to $d _{T}$, then there is a positive constant $c _{1}$ such that for any $1 \le p<\infty $, the following control of the expected supreme difference holds:
    \begin{equation}\label{eq: raw sp moment bound}
        \left(\mathbb{E}\sup\limits _{t \in T}\left|R _{t}-R _{t _{0}}\right|^{p}\right)^{\frac{1}{p}}\le c _{1}\gamma _{2}(T,d _{T})+2 \sup\limits _{t \in T}\left(\mathbb{E}\left\lVert R _{t}-X _{t _{0}}\right\rVert_{}^{p}\right)^{\frac{1}{p}},
    \end{equation}
    where $\gamma _{2}(T,d _{T})$ is the $\gamma _{2}$-functional w.r.t. $\left(R _{t}\right)_{t \in T}$ (see \cite{10.1214/aop/1065725175}).

    As a consequence, for any $u \ge 1$, there exist universal constants $c _{1},c _{2}$ such that
    \begin{equation}\label{eq: raw sp tail bound}
        \mathbb{P}\left(\sup\limits _{t \in T}\left|R _{t}-R _{t _{0}}\right|\ge c _{1}\gamma _{2}(T,d _{T})+c _{2}\sqrt{u}\Delta _{d _{T}}(T)\right)\le e ^{-u},
    \end{equation}
    where $\Delta _{d _{T}}(T):=\sup\limits _{t _{1},t _{2}\in T}d _{T}(t _{1},t _{2})$ is the diameter of $T$ under $d _{T}$.
\end{lemma}

The following lemma is a direct application of Lemma \ref{lemma: moment and tail bounds of sp}, where the corresponding stochastic process in our setting is defined as $R _{\theta } ^{(j)}:=\sum\limits _{i=1}^{n}\theta _{i}\left(T(Y _{i})-\mathbb{E}\left[T(Y _{i})\right]\right), \theta \in E _{(j)}$. The proof relies on a combination of Lemma \ref{lemma: moment and tail bounds of sp}, a concrete evaluation of the $\gamma _{2}$-functional of $(R _{\theta }^{(j)})_{\theta \in E _{(j)}}$ and $\Delta _{\ell _{2}}(E _{(j)})$, and a specific choice of $C$ in Lemma \ref{lemma: kolmogorov width control}. We defer the detailed proof to Appendix \ref{subsubsection: proof of the property of our stochastic process}.

\begin{lemma}[]\label{lemma: property of our stochastic process}
    In the $j$-th iteration, for the stochastic process $R _{\theta } ^{(j)}:=\sum\limits _{i=1}^{n}\theta _{i}\left(T(Y _{i})-\mathbb{E}\left[T(Y _{i})\right]\right), \theta \in E _{(j)}$, we have 
    \begin{equation}\label{eq: sp tail bound}
        \mathbb{P}\left(\sup\limits _{\theta \in E _{(j)}}\left|R _{\theta } ^{(j)}-R _{\theta ^{\prime }}^{(j)}\right|\gtrsim \left(\frac{\tilde{d}_{(j)}}{c}\right)^{2}\right)\lesssim \exp\left\{-\frac{\tilde{d}_{(j)}^{2}}{c ^{4}}\right\},
    \end{equation}
    where $\theta ^{\prime }=\hat{\theta }_{(j)}+A _{(j)}^{\dagger }\left(\theta _{0}-\hat{\theta }_{(j)}\right)$ and ``$\lesssim $'' omits only universal constants.
\end{lemma}

Building on Lemma \ref{lemma: high probability for close distance between theta prime and theta prime prime}, Lemma \ref{lemma: moment and tail bounds of sp}, Lemma \ref{lemma: property of our stochastic process}, Theorem \ref{theorem: polynomial time approximation to the SDP problem}, and Corollary \ref{corollary: control of the distance between the weak projection and any point in K}, the following core theorem establishes that the distance between the estimator $\hat{\theta }_{(j+1)}$ and the true parameter $\theta _{0}$ decays exponentially, provided that the previous iterate $\hat{\theta }_{(j)}$ is well-behaved. We kindly defer the proof to Appendix \ref{subsubsection: proof of theoretical guarantee of the efficient algorithm} for the reader's reference.

\begin{theorem}[]\label{theorem: exponentially decreasing distance between the estimator and true parameter}
    Suppose in the $j$-th iteration, we have $\left\lVert \hat{\theta }_{(j)}-\theta _{0}\right\rVert_{2}^{}\le \tilde{d}_{(j)}$ and $\tilde{d}_{(j)}$ satisfies the conditions of Lemma \ref{lemma: kolmogorov width control}, then with probability greater than $1-\exp\left\{-\frac{\tilde{d}_{(j)}^{2}}{c ^{4}}\right\}-\frac{1}{1+\tilde{d}_{(j)}^{3}}$, we have $\left\lVert \hat{\theta }_{(j+1)}-\theta _{0}\right\rVert_{2}^{}\le \sqrt{\frac{8C _{M}}{c _{M}}}\frac{\left(\sqrt{3}+1\right)\tilde{d}_{(j)}}{c}=\tilde{d}_{(j+1)}$.
\end{theorem}

In Algorithm \ref{algorithm: main algorithm}, we directly output an arbitrary point from $K$ when $r \gtrsim \sqrt{n}$, and iterate while $\tilde{d}_{(j)}\gtrsim r$ otherwise. The following lemma justifies this algorithmic design, demonstrating that such a case split does not result in any deviation from the theoretical minimax rate $\epsilon ^{*2}\wedge \text{diam}(K) ^{2}$. We defer the proof to Appendix \ref{subsubsection: proof of relationship between minimax rate and r}.

\begin{lemma}[]\label{lemma: relationship between minimax rate and r}
    If $r \gtrsim \sqrt{n}$, then $\epsilon ^{*}\asymp \sqrt{n}$; if $r \lesssim \sqrt{n}$, then $\epsilon ^{*}\gtrsim r$.
\end{lemma}

The following lemma plays a crucial role in bounding the $\ell _{2}$ error of the estimator $\hat{\theta }$ produced by Algorithm \ref{algorithm: main algorithm}. We will leverage this result subsequently in Lemma \ref{lemma: matching the minimax rate under ideal assumptions} and Theorem \ref{theorem: main theorem for the algorithm}. We place the proof of this lemma in Appendix \ref{subsubsection: proof of relationship between epsilon and diam}.

\begin{lemma}[]\label{lemma: relationship between epsilon and diam}
    Suppose $c \ge e$, then we have $\epsilon ^{*}\gtrsim _{M}1 \wedge \text{diam}(K) $, where ``$\gtrsim _{M}$'' omits the constant that only depends on $M$.
\end{lemma}

An important observation is that one essential condition required by Theorem \ref{theorem: exponentially decreasing distance between the estimator and true parameter} is that the previous iterate $\hat{\theta }_{(j)}$ must be well-behaved. This good behavior heavily depends on the validity of $X _{(j)}^{\dagger }$ and $A _{(j)}^{\dagger }$, as well as the event described in Lemma \ref{lemma: property of our stochastic process}. Therefore, to analyze the performance of the estimator $\hat{\theta }$ generated by Algorithm \ref{algorithm: main algorithm}, a natural approach is to partition the analysis into an ``ideal'' case, where $\hat{\theta }_{(j)}$ remains well-behaved throughout, and a ``non-ideal'' case, where at least one of the aforementioned conditions is violated sometime in the iteration (even though the algorithm itself is oblivious to this failure). We first consider the ideal case, assuming that the approximation oracle $X _{(j)}^{\dagger }$ is always valid and the event in Lemma \ref{lemma: property of our stochastic process} holds for all iterations $j=1,2,\dots$. Under this assumption, the estimator $\hat{\theta }$ achieves the minimax rate, as established by the following lemma, which is an immediate consequence of Lemmas \ref{lemma: relationship between minimax rate and r} and \ref{lemma: relationship between epsilon and diam}. We place its proof in Appendix \ref{subsubsection: proof of matching the minimax rate under ideal assumptions}.

\begin{lemma}[]\label{lemma: matching the minimax rate under ideal assumptions}
    Assume that $X _{(j)}^{\dagger }$ satisfies \eqref{eq: property of the approximation oracle} for $K _{(j)}$ and the event in Lemma \ref{lemma: property of our stochastic process} holds for $j=1,2,\dots$. Then we have
    \begin{equation}\label{eq: matching the minimax rate under ideal assumptions}
        \max\limits _{\theta _{0}\in K}\mathbb{E}_{Y}\left[\left\lVert \hat{\theta }_{\text{p}}(Y)-\theta _{0}\right\rVert_{2}^{2}\mathbf{1}_{\text{ideal}}\right]\asymp \epsilon ^{\star 2}\wedge \text{diam}(K) ^{2}.
    \end{equation}
    Equivalently, such $\hat{\theta }$ matches the minimax rate \eqref{eq: minimax rate}.
\end{lemma}

However, since it is possible that either $X _{(j)}^{\dagger }$ or the stochastic process $R _{\theta } ^{(j)}$ does not behave according to expectation, we must control the error once the algorithm deviates from its ideal trajectory. This control is established by the following lemma, which guarantees that if the algorithm fails to follow its expected path, the maximum $\ell _{2}$ error of the subsequent estimators will not exceed the distance bound $\tilde{d}_{(j)}$ from the last well-behaved iteration, up to a universal constant. The proof heavily depends on the good property of the weak projection of $K$, as established in Corollary \ref{corollary: control of the distance between the weak projection and any point in K}. We defer the detailed proof to Appendix \ref{subsubsection: proof of maximum error of the subsequent estimator}.

\begin{lemma}[]\label{lemma: maximum error of the subsequent estimator}
    Suppose Algorithm \ref{algorithm: main algorithm} runs normally ($X _{j}^{\dagger }$ is valid, and $R _{\theta }^{(j)}$ is well controlled, hence $\left\lVert \hat{\theta }_{(j)}-\theta _{0}\right\rVert_{2}^{}\le \tilde{d}_{(j)}$) until $j _{0}$-iteration. Then for any $j \ge j _{0}+1$, we have
    \begin{equation*}
        \left\lVert \hat{\theta }_{(j)}-\theta _{0}\right\rVert_{2}^{}\le \left(1+\frac{\sqrt{2+\sqrt{2}}}{1-\sqrt{q}}\right)^{2}\tilde{d}_{j _{0}}\lesssim \tilde{d}_{j _{0}},
    \end{equation*}
    where $q=\frac{\tilde{d}_{(j+1)}}{\tilde{d}_{(j)}}=\sqrt{\frac{8C _{M}}{c _{M}}}\frac{\sqrt{3}+1}{c}$ is from Algorithm \ref{algorithm: main algorithm}.
\end{lemma}

We now state and prove the primary theoretical guarantee for Algorithm \ref{algorithm: main algorithm}. Specifically, Theorem \ref{theorem: main theorem for the algorithm} formally establishes that for a broad class of convex constraints, the corresponding nonparametric regression estimator can be computed efficiently with a negligible statistical penalty.

\begin{theorem}[]\label{theorem: main theorem for the algorithm}
    For the exponential family \eqref{eq: exponential family} and the constraint $K$ for the natural parameter $\theta $, suppose that $K$ is define as $K=\left\{x \left\lvert\right. x \in \mathbb{R}^{n},\left\lVert Ax\right\rVert_{}^{}\le 1\right\}$, where $A \in \mathbb{R}^{m \times n}$ is known and $\left\lVert \cdot \right\rVert_{}^{}$ is an exactly $2$-convex, sign-invariant, and type-$2$ norm with small type-$2$ constant, and $K$ is well-balanced, i.e., there are known constants $r,R$ such that $B _{2}(0,r)\subset K \subset B _{2}(0,R)$. Let $c>2(C+1)\vee \sqrt{\frac{8C _{M}}{c _{M}}}\left(\sqrt{3}+1\right)$, then the estimator $\hat{\theta }_{\text{p}}$ produced by Algorithm \ref{algorithm: main algorithm} achieves
    \begin{equation}\label{eq: main inequality for the estimator}
        \mathcal{R}_{n}(f,K)\asymp _{M}\epsilon ^{*2}\wedge \operatorname{diam}(K) ^{2}\lesssim \sup\limits _{\theta _{0}\in K}\mathbb{E}_{Y}\left[\left\lVert \hat{\theta }_{\text{p}}(Y)-\theta _{0}\right\rVert_{2}^{2}\right]\lesssim _{M}\Upsilon (K)\epsilon ^{* 2}\wedge \text{diam}(K) ^{2},
    \end{equation}
    where $\Upsilon (K)$ only depends on $T _{2}(K)$ and $\kappa _{\text{m}}$ polynomially.  Furthermore, the runtime of the computation for $\hat{\theta }$ is polynomial in $n$ and $\log \left(\frac{R}{r}\right)$.
\end{theorem}

At a high level, the proof of Theorem \ref{theorem: main theorem for the algorithm} synthesizes the ideal-case guarantees from Lemma \ref{lemma: matching the minimax rate under ideal assumptions} with a bound on the $\ell_2$ estimation error under ``non-ideal'' conditions. This latter analysis crucially relies on Lemma \ref{lemma: maximum error of the subsequent estimator}. The complete formal proof is deferred to Appendix \ref{subsubsection: proof of main theorem for the algorithm}.

\subsection{Natural Embedding of Low-Dimensional Constraints}\label{subsection: natural embedding}

As a significant extension of Algorithm \ref{algorithm: main algorithm} and Theorem \ref{theorem: main theorem for the algorithm}, and expanding upon the remark following Theorem \ref{theorem: polynomial time approximation to the SDP problem}, we now consider the case where $K$ is a natural embedding of a lower-dimensional constraint. Suppose $\tilde{K} \subset \mathbb{R}^m$ is a well-balanced set represented as $\tilde{K} := \{x \in \mathbb{R}^m \mid \Vert{}Ax\Vert{} \le 1\}$, where $A \in \mathbb{R}^{p \times m}$ ($p \ge m$) is a known matrix and $\Vert{}\cdot\Vert{}$ is an exactly $2$-convex, sign-invariant norm on $\mathbb{R}^p$ with a controlled type-2 constant. Let $K$ be the high-dimensional representation of this constraint, defined as $K = U E_{m,n}(\tilde{K})$. Here, $U \in \mathbb{R}^{n \times n}$ is a known orthogonal transformation and $E_{m,n}: \mathbb{R}^m \to \mathbb{R}^n$ (for $m \le n$) is the standard zero-padding embedding, mapping $E_{m,n}(x) = (x, \mathbf{0})$. Additionally, we assume there exists a constant $M$, irrelevant to $n$, such that $K \subset [-M,M]^{n}, \tilde{K}\subset [-M,M]^{m}$. Under this structural assumption, there exists a computationally efficient estimator $\hat{\theta}_{\text{p}}$ whose $\ell_2$ risk is minimax optimal up to logarithmic factors. We establish this guarantee through the following arguments.

First, the diameter is preserved under the orthogonal embedding, yielding $\operatorname{diam}(K) = \operatorname{diam}(\tilde{K})$. Furthermore, there is a direct correspondence between their respective Kolmogorov $k$-widths. Specifically, the orthogonality of $U$ alongside the \hyperref[def_kolmogorov_width_raw]{definition} of the Kolmogorov width dictates that
\begin{equation*}
    D_k(K) =
    \begin{cases}
        D_k(\tilde{K}), & \text{if } 0 \le k \le m, \\
        0, & \text{if } m+1 \le k \le n.
    \end{cases}
\end{equation*}
Consequently, following Definition \ref{def_approximate_kolmogorov_width}, the approximation of $D_k(K)$ for $0 \le k \le m$ can be reformulated entirely within the lower-dimensional space:
\begin{equation}\label{eq: equivalent approximation}
    \begin{aligned}
        & \min_{X \in \mathcal{S}^m} \max_{\theta \in \tilde{K}} \theta^\top X \theta \\ 
        & \text{subject to} \quad \operatorname{tr}(X) = m-k, \quad \bm{0} \preceq X \preceq \mathbf{I}_{m}.
    \end{aligned}
\end{equation}
Any valid approximate solution to \eqref{eq: equivalent approximation} directly translates to an effective approximation for the original $n$-dimensional formulation. As a result, Lemma \ref{lemma: kolmogorov width control} and Theorem \ref{theorem: polynomial time approximation to the SDP problem} continue to hold when we approximate \eqref{eq: SDP problem} within the corresponding $m$-dimensional subspace. This dimensionality reduction allows us to efficiently compute $X _{(j)}^\dagger$ and $A _{(j)}^\dagger$ as described in Section \ref{subsection: efficient algorithm illustration}. By identical reasoning, setting $f = \Vert{}\cdot\Vert{}_2$ in Lemma \ref{lemma: raw approximation from dadush} and Theorem \ref{theorem: existence of the weak projection} enables us to restrict the computation of the weak projection purely to the $m$-dimensional subspace.

Second, the local metric entropy of $K$ identically matches that of $\tilde{K}$. Specifically, for any subset $K_1 \subset \mathbb{R}^m$ and any $\epsilon > 0$, we claim that $N(\epsilon, K_1) = N(\epsilon, E_{m,n}(K_1)) = N(\epsilon, UE_{m,n}(K_1))$. The second equality follows immediately from the orthogonal invariance of the Euclidean metric. To establish the first equality, observe that any $\epsilon$-covering $\{x_1, \dots, x_N\}$ of $K_1$ maps to an $\epsilon$-covering $\{E_{m,n}(x_1), \dots, E_{m,n}(x_N)\}$ of $E_{m,n}(K_1)$, yielding $N(\epsilon, E_{m,n}(K_1)) \le N(\epsilon, K_1)$. Conversely, let $\{z_1, \dots, z_N\}$ be an $\epsilon$-covering of $E_{m,n}(K_1)$. Define the projection operator $T_{n,m}: \mathbb{R}^n \to \mathbb{R}^m$ that simply drops the last $n-m$ coordinates. For any $y \in K_1$, there exists some $z_i$ such that $\Vert{}z_i - E_{m,n}(y)\Vert{}_2 \le \epsilon$. Because coordinate projection is a contraction mapping, we have $\Vert{}T_{n,m}(z_i) - y\Vert{}_2 \le \Vert{}z_i - E_{m,n}(y)\Vert{}_2 \le \epsilon$. Thus, the set $\{T_{n,m}(z_1), \dots, T_{n,m}(z_N)\}$ forms an $\epsilon$-covering of $K_1$, proving $N(\epsilon, K_1) \le N(\epsilon, E_{m,n}(K_1))$. In light of Theorem \ref{theorem: main theorem for theoretical algorithm} and Definition \ref{def_local_entropy}, this isometric equivalence dictates that the minimax rate in this setting is governed entirely by the $m$-dimensional local entropy of $\tilde{K}$. Consequently, Lemmas \ref{lemma: relationship between minimax rate and r} and \ref{lemma: relationship between epsilon and diam} remain fully applicable, provided we substitute the ambient dimension $n$ with $m$ and evaluate $r$ relative to $\tilde{K}$.

Finally, by the structural definition $K = U E_{m,n}(\tilde{K})$, any parameter vector $\theta \in K$ admits the representation $\theta = U \tilde{\theta}$, where $\tilde{\theta} = (\tilde{\theta}_1, \dots, \tilde{\theta}_m, 0, \dots, 0)^\top \in E_{m,n}(\tilde{K})$. Consequently, the empirical objective in Line 11 of Algorithm \ref{algorithm: main algorithm} can be rewritten as
\begin{equation*}
    \sum_{i=1}^n \left[ \theta_i T(y_i) - A(\theta_i) \right] = \sum_{i=1}^n \left[ \sum_{j=1}^m u_{i,j} \tilde{\theta}_j T(y_i) - A \left( \sum_{j=1}^m u_{i,j} \tilde{\theta}_j \right) \right],
\end{equation*}
where $u_{i,j}$ denotes the $(i,j)$-th entry of the orthogonal matrix $U$. Because the log-partition function $A(\cdot)$ is inherently convex, this reformulated objective is jointly concave with respect to the active components $(\tilde{\theta}_1, \dots, \tilde{\theta}_m)$. Therefore, the maximization step can be efficiently executed using the ellipsoid method directly within the lower-dimensional domain $\tilde{K}$.

In summary, we have established that all prerequisite lemmas and theorems for Theorem \ref{theorem: main theorem for the algorithm} remain fully valid under this embedding. Furthermore, the minimax estimation rate is strictly governed by the local metric entropy of the lower-dimensional set $\tilde{K}$. Consequently, we formally confirm the existence of a computationally efficient estimator in this setting.

\section{Example}\label{section: example}

\subsection{The Parametric Case: One-Dimensional Estimation}

Consider the classical parametric setting where we observe $Y_1, \dots, Y_n$ independently from a single exponential family distribution. In this scenario, the parameter vector is restricted to a one-dimensional subspace, meaning the constraint set $K \subset [-M, M]^n$ is simply a line segment. According to the \hyperref[def_local_entropy]{definition} of $N^{\text{loc}}(\epsilon, c)$, the local covering number satisfies
\begin{equation*}
    N^{\text{loc}}(\epsilon, c) \asymp 
    \begin{cases} 
        c, & \text{if } \epsilon < \text{diam}(K)/2, \\ \frac{c\text{diam}(K)}{2\epsilon}, & \text{if } \epsilon \ge \text{diam}(K)/2.
    \end{cases}
\end{equation*}
The transition between these two regimes occurs at $\epsilon = \text{diam}(K)/2 \asymp M\sqrt{n}$. By the \hyperref[def_critical_radius]{definition} of $\epsilon^*$, the critical radius must satisfy the balance equation $\epsilon^2 \kappa(M) \asymp \log N^{\text{loc}}(\epsilon, c)$. Assuming the solution falls within the first regime, equating these terms yields $\epsilon^* \asymp \sqrt{\frac{\log c}{\kappa(M)}}$. Because $c$ and $\kappa(M)$ are constants that depend solely on $M$, this implies $\epsilon^* \asymp_M 1$. As $n \to \infty$, we clearly have $\epsilon^* \asymp_M 1 \ll M\sqrt{n}$, which confirms that $\epsilon^*$ strictly resides in the first regime. Consequently, the minimax estimation rate is $\epsilon^{*2} \wedge \text{diam}(K)^2 \asymp_M 1$.

To see how this vector estimation rate aligns with classical parametric theory, recall that the constraint set $K$ restricts the parameter vector to a one-dimensional subspace where all components are identical. Consequently, we can express the true parameter vector as $\theta = \theta_0 \mathbf{1}_n$, where $\theta_0 \in [-M, M]$ is the underlying scalar natural parameter and $\mathbf{1}_n$ is the all-ones vector. The squared $\ell_2$ error for the vector $\theta$ is given by $\|\hat{\theta} - \theta\|^2 = n(\hat{\theta}_0 - \theta_0)^2$. Since we have established that the minimax rate for the vector evaluates to $\|\hat{\theta} - \theta\|^2 \asymp_M 1$, the mean squared error for the scalar parameter scales as $(\hat{\theta}_0 - \theta_0)^2 \asymp_M \frac{1}{n}$. 

This scalar rate seamlessly recovers the fundamental limit established by the Cramér-Rao lower bound for any unbiased estimator. For $n$ independent and identically distributed observations from a nonsingular exponential family, the Fisher information for $\theta_0$ is $\mathcal{I}_n(\theta_0) = n A''(\theta_0)$. Because $\theta_0$ is constrained to the compact interval $[-M, M]$ and the family is nonsingular, the unit Fisher information $A''(\theta_0)$ is bounded by some positive constants $c(M),C(M)$, implying $\mathcal{I}_n(\theta_0) \asymp_M n$. The Cramér-Rao lower bound dictates that the mean squared error is bounded below by $\mathcal{I}_n(\theta_0)^{-1} \asymp_M \frac{1}{n}$, which can be achieved by the MLE. Thus, our general nonparametric framework gracefully collapses to the exact minimax optimal parametric rate, rigorously verifying the consistency of our bounds with classical statistical theory. In addition, our results reveal that the general estimators are not significantly statistically better than the asymptotically unbiased estimators in this scenario.

\subsection{Constrained Bernoulli Regression}\label{subsection: constrained bernoulli regression}

In this section, we consider monotone constrained Bernoulli regression as a concrete example to illustrate the theoretical framework developed in Sections \ref{section: theoretical results} and \ref{section: efficient estimator}. Suppose we observe independent random variables $Y_i \sim \operatorname{Ber}(p_i)$ for $i \in [n]$, where the probabilities $p_i$ satisfy the bounded log-odds condition $\max_{i \in [n]} \left| \log \frac{p_i}{1-p_i} \right| \le M$ for some known constant $M$. Equivalently, this condition guarantees the existence of a known constant $\alpha := \alpha(M) > 0$ such that $\alpha \le \min\limits _{i \in [n]}p _{i} \le \max\limits _{i \in [n]} p _{i}\le 1 - \alpha$. Our primary objective is to estimate the natural parameters, namely the log-odds $\theta_i = \log \left( \frac{p_i}{1-p_i} \right)$.

We first note that the Bernoulli distribution belongs to the class of nonsingular exponential families. Under this parameterization, the sufficient statistic is $T(Y_i) = Y_i$, the natural parameter is the log-odds $\theta_i$, and the log-partition function is $A(\theta_i) = \log(1 + e^{\theta_i})$. Furthermore, this log-partition function is twice continuously differentiable, with its second derivative uniformly bounded by $1/4$.

We now introduce a star-shaped constraint on the parameter $\theta$. Our geometric setup draws inspiration from the multivariate isotonic regression framework of \cite{isotonic_general_dimensions}. Let $q \in \NN$ be a fixed dimensional constant independent of the sample size $n$. We define the lattice $\LL_{q,n} \subset \RR^q$ as 
\begin{equation*}
    \LL_{q,n} = \bigotimes \limits_{j=1}^{q}\left\{ \frac{1}{n^{1/q}}, \frac{2}{n^{1/q}}, \dots, 1 \right\},
\end{equation*}
which consists of exactly $(n^{1/q})^q = n$ points. In the univariate case ($q=1$), this structure reduces to $n$ equispaced points on the interval $(0, 1]$, separated by a distance of $n^{-1}$.

Let us define the function class
$$ \mathcal{F}_q^{[a,b]} := \{ f: [0,1]^q \to [a,b] \mid f \text{ is non-decreasing in each variable} \}. $$
Formally, a function $f$ is non-decreasing in each variable if, for any $(z_1, \dots, z_q) \in [0,1]^q$ and any $\delta_i > 0$ such that $z_i + \delta_i \le 1$, we have
\begin{equation*}
    f(z_1, \dots, z_{i-1}, z_i, z_{i+1}, \dots, z_q) \le f(z_1, \dots, z_{i-1}, z_i + \delta_i, z_{i+1}, \dots, z_q)
\end{equation*}
for all $1 \le i \le q$. We are specifically interested in the set of evaluations of these monotone functions on the lattice, defined as
\begin{equation*}
    Q_{-M,M} = \{ (f(l_1), f(l_2), \dots, f(l_n)) \mid f \in \mathcal{F}_q^{[-M,M]} \},
\end{equation*}
where $l_1, \dots, l_n$ denote the distinct elements of $\LL_{q,n}$.

In this example, we assume that the true parameter vector is $\theta = (\theta_1, \dots, \theta_n) \in Q_{-M,M}$. Observe that $Q_{-M,M}$ is a convex set and is therefore inherently star-shaped. Because the logit link function $x \mapsto \log\left(\frac{x}{1-x}\right)$ is strictly monotonically increasing and bijective on $(0,1)$, the condition $\theta \in Q_{-M,M}$ is strictly equivalent to $p = (p_1, \dots, p_n) \in Q_{\alpha, 1-\alpha}$, where
\begin{equation*}
    Q_{\alpha, 1-\alpha} = \{ (f(l_1), f(l_2), \dots, f(l_n)) \mid f \in \mathcal{F}_q^{[\alpha, 1-\alpha]} \}.
\end{equation*}
Under this structural assumption, our objective is to determine the minimax estimation rate for the vector $\theta$ with respect to the squared $\ell_2$ norm.

\subsubsection{One-Dimensional Case (\texorpdfstring{$q=1$}{q=1})}

For the one-dimensional case ($q = 1$), we derive the minimax rate by leveraging the analysis following Lemma 3.2 in \cite{prasadan2024}. Specifically, the local metric entropy satisfies
\begin{equation*}
    \log N_{Q_{-M,M}}^{\text{loc}}(\epsilon, c) \asymp
    \begin{cases} 
        \frac{\sqrt{2nM}}{\epsilon}, & \text{if } \epsilon \gtrsim \frac{2M}{\sqrt{n}}, \\
        n, & \text{otherwise.}
    \end{cases}
\end{equation*}
The critical radius $\epsilon^*$ is determined by solving the fundamental balance equation $\epsilon^2 \asymp \log N_{Q_{-M,M}}^{\text{loc}}(\epsilon, c)$. Substituting the entropy bound into this relation, we obtain
\begin{equation*}
    \epsilon^2 \asymp \frac{\sqrt{2nM}}{\epsilon} \quad \implies \quad \epsilon^3 \asymp \sqrt{n}.
\end{equation*}
Consequently, the minimax estimation rate is $\epsilon^2 \asymp n^{1/3}$. This result closely mirrors Theorem 2.3 in \cite{pmlr-v89-neykov19a}, though their result characterizes the estimation of the probability vector $p$ rather than the natural parameter $\theta$.

\subsubsection{Two-Dimensional Case (\texorpdfstring{$q=2$}{q = 2})}

For the two-dimensional case ($q = 2$), applying Lemmas 3.1 and 3.4 from \cite{prasadan2024} yields the following upper bound on the local metric entropy:
\begin{equation*}
    \log N^{\text{loc}}_{Q_{-M,M}}(\epsilon, c) \lesssim \left(\frac{\epsilon}{2\sqrt{2nM}}\right)^{-2} \left(\log \frac{2\sqrt{2nM}}{\epsilon}\right)^2.
\end{equation*}
To establish an upper bound on the estimation rate via Theorem \ref{theorem: main theorem for theoretical algorithm}, it suffices to identify the minimal order of $\epsilon$ that satisfies the fundamental balance condition $\epsilon^2 \gtrsim \log N^{\text{loc}}_{Q_{-M,M}}(\epsilon, c)$. Substituting the entropy bound into the right-hand side of this inequality, we require the smallest $\epsilon$ such that
\begin{equation*}
    \epsilon^4 \gtrsim 8nM \left(\log\frac{2\sqrt{2nM}}{\epsilon}\right)^2.
\end{equation*}
It is straightforward to verify that $\epsilon^2 \asymp \sqrt{n} \log n$ satisfies this condition. Consequently, the minimax upper bound for the estimation rate is given by $\sqrt{n} \log n$.

\subsubsection{Case with Higher Dimensions (\texorpdfstring{$q \ge 3$}{q >= 3})}

For higher dimensions ($q \geq 3$), we rely on Lemma 3.6 from \cite{prasadan2024} to derive the minimax rate.

\begin{lemma}[Lemma 3.6 of \cite{prasadan2024}]\label{lemma: akshayslemma}
Suppose that $\frac{\epsilon}{2M} \gtrsim \sqrt{\frac{n}{n^{1/q}}}$ and $q \ge 3$. Then, for a sufficiently large constant $c$, the local metric entropy satisfies
\begin{equation*}
    \log N^{\text{loc}}_{Q_{-M,M}}(\epsilon, c) \asymp \left( \frac{\epsilon}{2\sqrt{2nM}} \right)^{-2(q-1)}.
\end{equation*}
\end{lemma}
Consequently, the minimax rate is determined by solving the fundamental balance equation:
\begin{equation*}
    \epsilon^2 \asymp \left( \frac{\epsilon}{2\sqrt{2nM}} \right)^{-2(q-1)}.
\end{equation*}
Solving this relation yields $\epsilon \asymp n ^{(q-1)/2q}$. Note that this scaling validates the initial assumption required by Lemma \ref{lemma: akshayslemma}. Therefore, the minimax estimation rate is established as $\epsilon^2 \asymp n^{1- 1/q}$.

\subsection{Constrained Logistic Regression}

To demonstrate the practical utility of our efficient estimator, this section applies our framework to constrained logistic regression. Suppose we observe independent responses $Y_i \sim \operatorname{Ber}(p_i)$ for $i = 1, \dots, n$, where the success probabilities are parameterized as $p_i = \frac{1}{1 + \exp\left\{-\beta ^{\top } x _{i}\right\}},\beta \in \mathbb{R}^{p},x _{i}\in \mathbb{R}^{p}$, where $p$ is the number of predictors. We operate under a fixed-design setting, assuming that $x_1, \dots, x_n$ are linear independent and deterministic rather than random. We denote $X :=(x _{1},\dots,x _{n})^{\top } \in \mathbb{R}^{n \times p}$ with $\operatorname{rank}\left(X\right)=p$. Fixed-design experiments are prevalent in practice, particularly when experimental limitations—such as those encountered in destructive or sensitivity testing—restrict observations to predetermined predictor levels. For extensive discussions on the necessity and implications of fixed designs, we refer readers to \cite{dette2013efficiency} for clinical trials and toxicology; \cite{nelson2004accelerated} for reliability and stress testing in engineering; \cite{wichmann2001psychometric} for psychophysics; and \cite{10.1214/088342306000000105} for broader applications within generalized linear models (GLMs).

Under our settings, the log-likelihood function of the model can be written as 
\begin{equation*}
    \begin{aligned}
        \log \mathcal{L}(y _{1},\dots,y _{n};\beta )&=\log \left(\prod\limits _{i=1}^{n}p _{i}^{y _{i}}(1-p _{i})^{y _{i}}\right)\\ 
        &=\log \left[\prod\limits _{i=1}^{n}\left(\frac{1}{1+\exp\left\{-\beta ^{\top } x _{i}\right\}}\right)^{y _{i}}\cdot \left(\frac{\exp\left\{-\beta ^{\top } x _{i}\right\}}{1+\exp\left\{-\beta ^{\top }x _{i}\right\}}\right)^{1-y _{i}}\right]\\ 
        &=-\sum\limits _{i=1}^{n}\log \left(\exp\left\{\beta ^{\top }x _{i}\right\}+1\right)+\sum\limits _{i=1}^{n}\beta ^{\top } x _{i}y _{i}.
    \end{aligned}
\end{equation*}
Compared to \eqref{eq: exponential family}, we know that the natural parameter is $\theta =(\theta _{1},\dots,\theta _{n})^{\top }=X \beta $, the sufficient statistic is $T(Y _{i})=Y _{i}$, and $A(\theta _{i})=\log \left(e ^{\theta _{i}}+1\right)$ with $A ^{\prime \prime }(\theta _{i})=e ^{\theta _{i}}/\left(e ^{\theta _{i}}+1\right)^{2}\in (0,\frac{1}{4}]$.

We consider an ellipsoidal constraint for $\beta $ of the form $\tilde{K} = \{\beta \in \mathbb{R}^p \left\lvert\right. \beta^\top \tilde{Q} \beta \le 1\}$ for some positive definite matrix $\tilde{Q} \in \mathbb{R}^{p \times p}$. We denote the lengths of the semi axes of the ellipsoid as $\sqrt{\tilde{a} _{1}},\sqrt{\tilde{a} _{2}},\dots,\sqrt{\tilde{a} _{p}}$. By permuting the coordinate axes, we can assume without loss of generality that $0 <\sqrt{\tilde{a}_{1}} \le \sqrt{\tilde{a}_{2}} \le \dots \le \sqrt{\tilde{a}_{p}}:=\tilde{M}$. Intuitively, such constraint for $\beta $ implies a constraint for $\theta $. If we solve $\beta $ from some $\theta \in \mathcal{C}(X)$, where $\mathcal{C}(\cdot )$ denotes the column space, we have $\beta =\left(X ^{\top } X\right)^{-1}X ^{\top } \theta $ and consequently the constraint for $\theta $ is:
\begin{equation}\label{eq: constraint for theta}
    \begin{aligned}
        &K=\left\{\theta \left\lvert\right. \theta \in \mathbb{R}^{n},\theta ^{\top } Q \theta \le 1,\theta \in \mathcal{C}(X)\right\},\\ 
        &Q :=X \left(X ^{\top } X\right)^{-1}\tilde{Q}\left(X ^{\top } X\right)^{-1}X ^{\top } .
    \end{aligned}
\end{equation}
We know from \eqref{eq: constraint for theta} that $K$ is a degenerate $p$-dimensional ellipsoid in $\mathbb{R}^{n}$. Similarly, we denote the lengths of the semi axes of this ellipsoid as $\sqrt{a _{1}},\dots,\sqrt{a _{n}}$ and assume $0=a _{1}=\dots=a _{n-p}<a _{n-p+1}\le \dots \le a _{n}:=M$. Furthermore, there is a continuous bijection between $(\tilde{a} _{1},\dots,\tilde{a} _{p},\tilde{M})$ and $(a _{n-p+1},\dots,a _{n},M)$.

For the critical radius $\epsilon ^{*}$, we have the following lemma, whose proof is directly from the \hyperref[def_critical_radius]{definition} of $\epsilon ^{*}$, standard properties of the Kolmogorov $k$-width, and entropy number. The proof is deferred to Appendix \ref{subsubsection: proof of critical radius for ellipsoid}.

\begin{lemma}[]\label{lemma: critical radius for ellipsoid}
    Define the critical index $k ^{*}=\min\limits \left\{k \left\lvert\right. k \ge 0, D _{k}^{2}(K)=a _{n-k}\le (k+1)/\kappa (M)\right\}$, where we recall $D _{k}(K)$ is the Kolmogorov $k$-width of the set $K$ (Definition \ref{def_kolmogorov_width_raw}). Then, we have $\epsilon ^{*}\asymp _{M}\sqrt{k ^{*}}$ given $k ^{*}\neq 0$.
\end{lemma}

As a corollary, Lemma \ref{lemma: critical radius for ellipsoid} yields the following minimax rate for the constrained logistic regression. We place its proof in Appendix \ref{subsubsection: proof of minimax rate for constrained logistic regression}.

\begin{corollary}[]\label{corollary: minimax rate for constrained logistic regression}
    The minimax rate of the $\ell _{2}$ estimation error for the natural parameter $\theta _{0}=X\beta _{0}$ in the constrained logistic regression is
    \begin{equation*}
        \mathcal{R}_{n}(f _{\text{logistic}},K)\asymp _{M}(k ^{*}+1)\wedge a _{n}.
    \end{equation*}
\end{corollary}

\begin{remark*}[]
    We briefly remark on the case where the constraint $K$ is a hyperrectangle with semi-axes $a_1, \dots, a_n$. It is straightforward to verify that in this setting, the minimax rate analysis for the $\ell_2$ estimation error closely parallels that of the ellipsoidal case. Specifically, if we alternatively define the critical index as
    \begin{equation*}
        k^* := \min \left\{ k \ge 0 \mid \sum_{i=1}^{k} a_{n-i}^2 \le \frac{k+1}{\kappa(M)} \right\},
    \end{equation*}
    where $M :=\sqrt{\sum\limits _{i=1}^{n}a _{i}^{2}}$, the minimax rate \eqref{eq: minimax risk of the estimation} is analogously given by $(k^*+1) \wedge \operatorname{diam}(K)^2 \asymp (k^*+1) \wedge \sum\limits _{i=1}^{n} a_i^2$.
\end{remark*}

A notable feature of this example is that the efficient estimator described in Algorithm \ref{algorithm: main algorithm} can be derived under relatively mild conditions. For fixed $x_1, \dots, x_n$, the log-partition function $A(\theta )$ naturally satisfies all distributional requirements for constructing the estimator. Next, we verify the geometric conditions on the constraint set $K$. Since $\theta =X \beta $, we have $v ^{\top } \theta =0$ for any $v \in \operatorname{Ker}(X ^{\top } )$. From \eqref{eq: constraint for theta}, we know that $Q \in \mathbb{R}^{n \times n}$ satisfies $\theta \in \operatorname{Ker}(Q)\Leftrightarrow Q \theta =0 \Leftrightarrow X ^{\top } \theta =0 \Leftrightarrow \theta \in \operatorname{Ker}(X ^{\top } )$. Consider the eigenvalue decomposition for $Q$:
\begin{equation*}
    \begin{aligned}
        & Q=U\Lambda U ^{\top },\\ 
        & U=(u _{1},\dots,u _{p},u _{p+1},\dots,u _{n}),\\ 
        & \Lambda =\operatorname{diag}\{1/a _{n-p+1},\dots,1/a _{n},0,\dots,0\}.
    \end{aligned}
\end{equation*}
We have $\operatorname{Ker}(Q)=\operatorname{Ker}(X ^{\top } )=\operatorname{span}\left\{u _{p+1},\dots,u _{n}\right\}$. Let $\tilde{\theta }=U ^{\top } \theta =(\tilde{\theta }_{1},\dots,\tilde{\theta }_{p},0,\dots,0)$, then $K$ can be expressed as 
\begin{equation*}
    K=UE _{p,n}(K ^{\prime }), K ^{\prime }=\left\{\tilde{\theta }_{p} \mid \tilde{\theta }_{p}\in \mathbb{R}^{p},\sum\limits _{i=1}^{p}\frac{\tilde{\theta }_{p,i}^{2}}{a _{n-i+1}}\le 1\right\},
\end{equation*}
where $E _{p,n}$ is the natural embedding mapping from $\mathbb{R}^{p}$ to $\mathbb{R}^{n}$ described in Section \ref{subsection: natural embedding}. By Section \ref{subsection: natural embedding}, it suffices to verify the geometric conditions on $K ^{\prime }$. An axis-aligned ellipsoid is inherently symmetric and well-balanced, satisfying $B_2(0, a _{n-p+1}) \subseteq K ^{\prime }\subseteq B_2(0, a_n)$. Its Minkowski gauge functional takes the explicit form
\begin{equation*}
    \rho_K(\tilde{\theta }_{p}) = \sqrt{\sum_{i=1}^p \frac{\tilde{\theta } _{p,i}^2}{a_{n-i+1}}}, \quad \forall \tilde{\theta }_{p} \in \mathbb{R}^p.
\end{equation*}
Consequently, for any integer $m \ge 1$ and vectors $\tilde{\theta }_{p}^{(1)}, \dots, \tilde{\theta }_{p}^{(m)}\in \mathbb{R}^p$, it is straightforward to verify that both sides of \eqref{eq: def_exactly_2_convex_norm} evaluate to
\begin{equation*}
    \sqrt{\sum_{i=1}^m \sum_{j=1}^p \frac{\tilde{\theta }_{p,j}^{(i)2}}{a_{n-j+1}}}.
\end{equation*}
Thus, the Minkowski gauge is exactly 2-convex. It remains only to bound the type-2 constant $T_2(K)$, ensuring it grows sufficiently slowly relative to $p \le n$. Specifically, we establish that $T_2(K ^{\prime }) \lesssim \log p$. This bound follows from standard probabilistic techniques, namely by bounding Rademacher complexities via standard Gaussian complexities and subsequently applying a maximal inequality.

Consequently, Algorithm \ref{algorithm: main algorithm} guarantees the efficient construction of the estimator $\hat{\theta }_{\text{p}}$. This procedure achieves polynomial-time complexity provided that $\log \frac{R}{r}$ is bounded by a poly-logarithmic factor in $n$. Given the box constraint $R = a_n = M$, this requirement simplifies to the condition that $|\log r| = |\log a_1|$ grows at most poly-logarithmically with $n$. Geometrically, this ensures that the embedded constraint ellipsoid $K ^{\prime }$ is not severely ill-conditioned.

As established in Theorem \ref{theorem: main theorem for the algorithm}, the $\ell_2$ risk of the estimator $\hat{\theta }_{\text{p}}$ is bounded by
\begin{equation*}
    \max_{\theta \in K} \mathbb{E}_{Y_1,\dots,Y_n} \|\hat{\theta }_{\text{p}}(Y) - \theta \|_2^2 \lesssim_M \Upsilon(K ^{\prime })(k^* + 1) \wedge a_n,
\end{equation*}
where the geometric factor $\Upsilon(K ^{\prime })$ is explicitly given by
\begin{equation*}
    \Upsilon(K) = 2 \sqrt{\tau_{K ^{\prime }_{\text{m}}} T_2(K ^{\prime })} \vee c^2 \sqrt{\lceil 4 \lceil \log_c(\sqrt{\tau_{K ^{\prime }_{\text{m}}}} T_2(K ^{\prime })) \rceil \rceil} \lesssim_M \log p \le \log n.
\end{equation*}
This bound demonstrates that the performance of $\hat{\theta }_{\text{p}}$ is minimax optimal up to a logarithmic factor. As a byproduct, we can further naturally estimate $\beta $ via $\hat{\beta }_{\text{p}}=\left(X ^{\top } X\right)^{-1}X ^{\top }\hat{\theta }_{\text{p}}$. Given fixed $X$, such estimation of $\beta _{0}$ is also minimax optimal up to a logarithmic factor, as there is a bi-Lipschitz homeomorphism between $\mathbb{R}^{p}$ and $\mathcal{C}(X)$. (Note that $\hat{\theta }_{\text{p}}\in \mathcal{C}(X)$ is guaranteed.)

\section{Discussion}\label{section: discussion}

In this paper, we investigated the minimax $\ell_2$ estimation rate for nonparametric exponential family regression under star-shaped constraints. By establishing matching lower and upper bounds (Sections \ref{subsection: lower bound} and \ref{subsection: upper bound}), Theorem \ref{theorem: main theorem for theoretical algorithm} formally characterizes this exact minimax rate. Furthermore, we proposed Algorithm \ref{algorithm: theoretical algorithm}, which is theoretically guaranteed to attain this optimal rate for any sufficiently smooth, nonsingular distribution within the exponential family. Our work significantly generalizes the results of \cite{neykov2023minimax}, which focused exclusively on the Gaussian sequence model under strictly convex constraints.

Beyond statistical optimality, we also address the computational feasibility of our estimators. While Algorithm \ref{algorithm: theoretical algorithm} establishes the fundamental theoretical limits, it is primarily of theoretical interest; explicitly constructing the sequence of packing sets in Algorithm \ref{algorithm: construction of the pruned tree} is computationally prohibitive in practice. To bridge this gap, Section \ref{section: efficient estimator} introduces Algorithm \ref{algorithm: main algorithm}, a computationally efficient procedure that operates in polynomial time under mildly stronger geometric assumptions on the constraint set $K$. Crucially, this efficient estimator retains minimax optimality, matching the theoretical rate up to poly-logarithmic factors in the ambient dimension $n$ and the geometric parameters of $K$.

To illustrate the broad utility of our framework, Section \ref{section: example} applies these general results to explicitly derive the minimax rate for a various of examples and the efficient estimator is also constructed when the geometric conditions are satisfied.

\subsection{Future Work}\label{subsection: future work}
This work opens several avenues for future research. First, as discussed in our problem formulation, the constraint set $K$ is assumed to be contained within a bounded box $[-M,M]^{n}$ for a known $M$. While this boundedness assumption is standard in the literature \cite{neykov2023minimax,prasadan2024characterizing,10.1214/25-AOS2576}, eliminating the dependency on $M$ in the minimax rate would represent a significant theoretical advance. Recently, \cite{10.1214/25-AOS2576} extended their analysis to unbounded star-shaped constraints under Gaussian and sub-Gaussian noise models by utilizing a sequence of nested bounded sets. However, it remains an open question whether this technique can be adapted to general distributions within the exponential family. Furthermore, as acknowledged by the authors of \cite{10.1214/25-AOS2576}, their proposed estimator is computationally intractable. Consequently, developing a computationally efficient estimator for unbounded star-shaped constraints poses a formidable and highly impactful challenge for future work.

Another compelling direction is to investigate whether the efficient estimator constructed in Algorithm \ref{algorithm: main algorithm} can be generalized to broader classes of constraints. In Section \ref{section: efficient estimator}, we imposed several geometric and structural assumptions on $K$ to guarantee computational efficiency. While not overly restrictive, these conditions currently exclude some highly relevant settings, such as $s$-sparse constraints. Consequently, a critical open question is whether efficient estimation remains possible under these more general constraints, and if so, what the minimal structural requirements might be. Conversely, it would be highly instructive to establish formal computational lower bounds—demonstrating a fundamental statistical-computational gap (e.g., \cite{brennan2018reducibility, zhang2014lower})—under which polynomial-time estimation achieving the minimax rate becomes provably intractable.

Finally, a natural but statistically fundamental extension of our work is the corresponding hypothesis testing problem. Specifically, consider testing the null hypothesis $H_0: \beta = \beta_0$ against the alternative $H_1: \beta \neq \beta_0$, where the natural parameter $\beta$ is known a priori to belong to a bounded star-shaped constraint $K$. The fundamental objective in this setting is to determine the minimax separation rate between $H_0$ and $H_1$ under a specified distance metric. For the Gaussian sequence model subject to general quadratically convex orthosymmetric (QCO) constraints, this problem has been resolved by \cite{li2026robustsignaldetectionquadratically,doi:10.1137/24M1636435}, even in the presence of data contamination. However, for general distributions within the regular exponential family, this domain remains largely unexplored. While it is a well-established statistical principle that testing is not harder than estimation—typically yielding faster minimax rates—precisely characterizing how these rates improve under our generalized framework remains an intriguing open question.

\bibliographystyle{abbrvnat}
\bibliography{citations}

\begin{thebibliography}{57}
\providecommand{\natexlab}[1]{#1}
\providecommand{\url}[1]{\texttt{#1}}
\expandafter\ifx\csname urlstyle\endcsname\relax
  \providecommand{\doi}[1]{doi: #1}\else
  \providecommand{\doi}{doi: \begingroup \urlstyle{rm}\Url}\fi

\bibitem[Bhattiprolu et~al.(2021)Bhattiprolu, Lee, and
  Naor]{10.1145/3406325.3451128}
V.~Bhattiprolu, E.~Lee, and A.~Naor.
\newblock A framework for quadratic form maximization over convex sets through
  nonconvex relaxations.
\newblock In \emph{Proceedings of the 53rd Annual ACM SIGACT Symposium on
  Theory of Computing}, STOC 2021, page 870–881, New York, NY, USA, 2021.
  Association for Computing Machinery.
\newblock ISBN 9781450380539.
\newblock \doi{10.1145/3406325.3451128}.
\newblock URL \url{https://doi.org/10.1145/3406325.3451128}.

\bibitem[Bickel and Doksum(2015)]{Bickel2015}
P.~J. Bickel and K.~A. Doksum.
\newblock \emph{Mathematical Statistics: Basic Ideas and Selected Topics}.
\newblock CRC Press, 2nd edition, 2015.

\bibitem[Boucheron et~al.(2013)Boucheron, Lugosi, and
  Massart]{10.1093/acprof:oso/9780199535255.001.0001}
S.~Boucheron, G.~Lugosi, and P.~Massart.
\newblock \emph{Concentration Inequalities: A Nonasymptotic Theory of
  Independence}.
\newblock Oxford University Press, 02 2013.
\newblock ISBN 9780199535255.
\newblock \doi{10.1093/acprof:oso/9780199535255.001.0001}.
\newblock URL \url{https://doi.org/10.1093/acprof:oso/9780199535255.001.0001}.

\bibitem[Brennan et~al.(2018)Brennan, Bresler, and
  Huleihel]{brennan2018reducibility}
M.~Brennan, G.~Bresler, and W.~Huleihel.
\newblock Reducibility and computational lower bounds for problems with planted
  sparse structure.
\newblock In \emph{Conference On Learning Theory}, pages 48--166. PMLR, 2018.

\bibitem[Brieden(2002)]{brieden2002geometric}
Brieden.
\newblock Geometric optimization problems likely not contained in.
\newblock \emph{Discrete \& Computational Geometry}, 28\penalty0 (2):\penalty0
  201--209, 2002.

\bibitem[Brown et~al.(2010)Brown, Cai, and Zhou]{brown2010nonparametric}
L.~D. Brown, T.~T. Cai, and H.~H. Zhou.
\newblock Nonparametric regression in exponential families.
\newblock \emph{The Annals of Statistics}, 38\penalty0 (4):\penalty0
  2005--2046, 2010.

\bibitem[Cai(2012)]{Cai_2012}
T.~T. Cai.
\newblock Minimax and adaptive inference in nonparametric function estimation.
\newblock \emph{Statistical Science}, 27\penalty0 (1), Feb. 2012.
\newblock ISSN 0883-4237.
\newblock \doi{10.1214/11-sts355}.
\newblock URL \url{http://dx.doi.org/10.1214/11-STS355}.

\bibitem[Cai and Zhou(2010)]{cai2010nonparametric}
T.~T. Cai and H.~H. Zhou.
\newblock Nonparametric regression in natural exponential families.
\newblock In \emph{Borrowing Strength: Theory Powering Applications–A
  Festschrift for Lawrence D. Brown}, volume~6, pages 199--216. Institute of
  Mathematical Statistics, 2010.

\bibitem[Canonne et~al.(2026)Canonne, Hopkins, Li, Liu, and
  Narayanan]{doi:10.1137/24M1636435}
C.~Canonne, S.~B. Hopkins, J.~Li, A.~Liu, and S.~Narayanan.
\newblock The full landscape of robust mean testing: Sharp separations between
  oblivious and adaptive contamination.
\newblock \emph{SIAM Journal on Computing}, 55\penalty0 (3):\penalty0
  FOCS23--192--FOCS23--267, 2026.
\newblock \doi{10.1137/24M1636435}.
\newblock URL \url{https://doi.org/10.1137/24M1636435}.

\bibitem[Casella and Berger(2002)]{Casella2002}
G.~Casella and R.~L. Berger.
\newblock \emph{Statistical Inference}.
\newblock Thomson Learning, 2nd edition, 2002.

\bibitem[Chatterjee et~al.(2015)Chatterjee, Guntuboyina, and
  Sen]{chatterjee2015risk}
S.~Chatterjee, A.~Guntuboyina, and B.~Sen.
\newblock On risk bounds in isotonic and other shape restricted regression
  problems.
\newblock \emph{The Annals of Statistics}, 43\penalty0 (4):\penalty0
  1774--1800, 2015.
\newblock \doi{10.1214/15-AOS1324}.

\bibitem[Cohen and DeVore(2015)]{10.1093/imanum/dru066}
A.~Cohen and R.~DeVore.
\newblock Kolmogorov widths under holomorphic mappings.
\newblock \emph{IMA Journal of Numerical Analysis}, 36\penalty0 (1):\penalty0
  1--12, 03 2015.
\newblock ISSN 0272-4979.
\newblock \doi{10.1093/imanum/dru066}.
\newblock URL \url{https://doi.org/10.1093/imanum/dru066}.

\bibitem[Dadush(2012)]{Dadush2012IntegerPL}
D.~Dadush.
\newblock Integer programming, lattice algorithms, and deterministic volume
  estimation.
\newblock 2012.
\newblock URL \url{https://api.semanticscholar.org/CorpusID:126363382}.

\bibitem[DasGupta(2011)]{expfamily}
A.~DasGupta.
\newblock The exponential family and statistical applications.
\newblock In \emph{Probability for Statistics and Machine Learning:
  Fundamentals and Advanced Topics}, chapter~18, pages 583--612. Springer New
  York, 2011.
\newblock ISBN 978-1-4419-9633-6.
\newblock \doi{10.1007/978-1-4419-9634-3_18}.

\bibitem[Dette et~al.(2013)Dette, Bornkamp, and Bretz]{dette2013efficiency}
H.~Dette, B.~Bornkamp, and F.~Bretz.
\newblock On the efficiency of two-stage response-adaptive designs.
\newblock \emph{Statistics in Medicine}, 32\penalty0 (10):\penalty0 1646--1660,
  2013.

\bibitem[Dirksen(2015)]{10.1214/EJP.v20-3760}
S.~Dirksen.
\newblock {Tail bounds via generic chaining}.
\newblock \emph{Electronic Journal of Probability}, 20\penalty0
  (none):\penalty0 1 -- 29, 2015.
\newblock \doi{10.1214/EJP.v20-3760}.
\newblock URL \url{https://doi.org/10.1214/EJP.v20-3760}.

\bibitem[Donoho and Liu(1991)]{donoho1991geometrizing}
D.~L. Donoho and R.~C. Liu.
\newblock Geometrizing rates of convergence, iii.
\newblock \emph{The Annals of Statistics}, 19\penalty0 (2):\penalty0 668--701,
  1991.
\newblock \doi{10.1214/aos/1176348115}.

\bibitem[Donoho et~al.(1990)Donoho, Liu, and MacGibbon]{10.1214/aos/1176347758}
D.~L. Donoho, R.~C. Liu, and B.~MacGibbon.
\newblock {Minimax Risk Over Hyperrectangles, and Implications}.
\newblock \emph{The Annals of Statistics}, 18\penalty0 (3):\penalty0 1416 --
  1437, 1990.
\newblock \doi{10.1214/aos/1176347758}.
\newblock URL \url{https://doi.org/10.1214/aos/1176347758}.

\bibitem[Emery et~al.(2007)]{Emery2007}
M.~Emery et~al.
\newblock \emph{Lectures on Probability Theory and Statistics: Ecole d'Été de
  Probabilités de Saint-Flour XXVIII - 1998}, volume 1738.
\newblock Springer Berlin / Heidelberg, 1st edition, 2007.
\newblock \doi{10.1007/BFb0106703}.
\newblock URL \url{https://doi.org/10.1007/BFb0106703}.

\bibitem[Evans et~al.(2009)Evans, Bazilevs, Babuska, and Hughes]{EVANS20091726}
J.~A. Evans, Y.~Bazilevs, I.~Babuska, and T.~J. Hughes.
\newblock n-widths, sup–infs, and optimality ratios for the k-version of the
  isogeometric finite element method.
\newblock \emph{Computer Methods in Applied Mechanics and Engineering},
  198\penalty0 (21):\penalty0 1726--1741, 2009.
\newblock ISSN 0045-7825.
\newblock \doi{https://doi.org/10.1016/j.cma.2009.01.021}.
\newblock URL
  \url{https://www.sciencedirect.com/science/article/pii/S0045782509000280}.
\newblock Advances in Simulation-Based Engineering Sciences – Honoring J.
  Tinsley Oden.

\bibitem[Fan(1993)]{fan1993local}
J.~Fan.
\newblock Local linear regression smoothers and their minimax efficiencies.
\newblock \emph{Annals of Statistics}, 21\penalty0 (1):\penalty0 196--216,
  1993.
\newblock \doi{10.1214/aos/1176349022}.

\bibitem[Floater et~al.(2021)Floater, Manni, Sande, and Speleers]{Floater_2021}
M.~S. Floater, C.~Manni, E.~Sande, and H.~Speleers.
\newblock Best low-rank approximations and kolmogorov $n$-widths.
\newblock \emph{SIAM Journal on Matrix Analysis and Applications}, 42\penalty0
  (1):\penalty0 330–350, Jan 2021.
\newblock ISSN 1095-7162.
\newblock \doi{10.1137/20m1355720}.
\newblock URL \url{http://dx.doi.org/10.1137/20M1355720}.

\bibitem[Florian~Arbes and Urban(2025)]{Arbes2025}
C.~G. Florian~Arbes and K.~Urban.
\newblock The kolmogorov n-width for linear transport: exact representation and
  the influence of the data.
\newblock \emph{Adv Comput Math}, 51\penalty0 (13), 2025.
\newblock \doi{10.1007/s10444-025-10224-0}.
\newblock URL \url{https://doi.org/10.1007/s10444-025-10224-0}.

\bibitem[Gao et~al.(2020)Gao, Han, and Zhang]{10.1214/18-AOS1792}
C.~Gao, F.~Han, and C.-H. Zhang.
\newblock {On estimation of isotonic piecewise constant signals}.
\newblock \emph{The Annals of Statistics}, 48\penalty0 (2):\penalty0 629 --
  654, 2020.
\newblock \doi{10.1214/18-AOS1792}.
\newblock URL \url{https://doi.org/10.1214/18-AOS1792}.

\bibitem[Greif and Urban(2019)]{GREIF2019216}
C.~Greif and K.~Urban.
\newblock Decay of the kolmogorov n-width for wave problems.
\newblock \emph{Applied Mathematics Letters}, 96:\penalty0 216--222, 2019.
\newblock ISSN 0893-9659.
\newblock \doi{https://doi.org/10.1016/j.aml.2019.05.013}.
\newblock URL
  \url{https://www.sciencedirect.com/science/article/pii/S0893965919301983}.

\bibitem[Gr{\"o}tschel et~al.(1993)Gr{\"o}tschel, Lov{\'a}sz, and
  Schrijver]{grotschel1993geometric}
M.~Gr{\"o}tschel, L.~Lov{\'a}sz, and A.~Schrijver.
\newblock \emph{Geometric Algorithms and Combinatorial Optimization}, volume~2
  of \emph{Algorithms and Combinatorics}.
\newblock Springer-Verlag, Berlin, Heidelberg, 2nd edition, 1993.

\bibitem[Han et~al.(2019)Han, Wang, Chatterjee, and
  Samworth]{isotonic_general_dimensions}
Q.~Han, T.~Wang, S.~Chatterjee, and R.~J. Samworth.
\newblock {Isotonic regression in general dimensions}.
\newblock \emph{The Annals of Statistics}, 47\penalty0 (5):\penalty0 2440 --
  2471, 2019.
\newblock \doi{10.1214/18-AOS1753}.

\bibitem[Jordan(2009)]{Jordan2009}
M.~I. Jordan.
\newblock The exponential family: Basics, 2009.
\newblock URL
  \url{https://people.eecs.berkeley.edu/~jordan/courses/260-spring10/other-readings/chapter8.pdf}.

\bibitem[Khuri et~al.(2006)Khuri, Mukherjee, Sinha, and
  Ghosh]{10.1214/088342306000000105}
A.~I. Khuri, B.~Mukherjee, B.~K. Sinha, and M.~Ghosh.
\newblock {Design Issues for Generalized Linear Models: A Review}.
\newblock \emph{Statistical Science}, 21\penalty0 (3):\penalty0 376 -- 399,
  2006.
\newblock \doi{10.1214/088342306000000105}.
\newblock URL \url{https://doi.org/10.1214/088342306000000105}.

\bibitem[Kolmogoroff(1936)]{2ad9b75f-e5d2-3a8c-8211-9d37aa110eb2}
A.~Kolmogoroff.
\newblock Über die beste annäherung von funktionen einer gegebenen
  funktionenklasse.
\newblock \emph{Annals of Mathematics}, 37\penalty0 (1):\penalty0 107--110,
  1936.
\newblock ISSN 0003486X, 19398980.
\newblock URL \url{http://www.jstor.org/stable/1968691}.

\bibitem[Li and
  Neykov(2026{\natexlab{a}})]{li2026efficientrobustconstrainedsignal}
Y.~Li and M.~Neykov.
\newblock Efficient robust constrained signal detection via kolmogorov width
  approximations, 2026{\natexlab{a}}.
\newblock URL \url{https://arxiv.org/abs/2605.11238}.

\bibitem[Li and
  Neykov(2026{\natexlab{b}})]{li2026robustsignaldetectionquadratically}
Y.~Li and M.~Neykov.
\newblock Robust signal detection with quadratically convex orthosymmetric
  constraints, 2026{\natexlab{b}}.
\newblock URL \url{https://arxiv.org/abs/2308.13036}.

\bibitem[McCullagh and Nelder(1989)]{mccullagh1989generalized}
P.~McCullagh and J.~A. Nelder.
\newblock \emph{Generalized Linear Models}.
\newblock Chapman and Hall/CRC Monographs on Statistics and Applied Probability
  Series. Chapman \& Hall, 1989.
\newblock ISBN 9780412317606.
\newblock URL \url{http://books.google.com/books?id=h9kFH2_FfBkC}.

\bibitem[Nelson(2004)]{nelson2004accelerated}
W.~B. Nelson.
\newblock \emph{Accelerated Testing: Statistical Models, Test Plans, and Data
  Analysis}.
\newblock Wiley Series in Probability and Statistics. John Wiley \& Sons, 2004.

\bibitem[Neykov(2019)]{pmlr-v89-neykov19a}
M.~Neykov.
\newblock Tossing coins under monotonicity.
\newblock In K.~Chaudhuri and M.~Sugiyama, editors, \emph{Proceedings of the
  Twenty-Second International Conference on Artificial Intelligence and
  Statistics}, volume~89 of \emph{Proceedings of Machine Learning Research},
  pages 21--30. PMLR, 16--18 Apr 2019.
\newblock URL \url{https://proceedings.mlr.press/v89/neykov19a.html}.

\bibitem[Neykov(2023)]{neykov2023minimax}
M.~Neykov.
\newblock On the minimax rate of the gaussian sequence model under bounded
  convex constraints.
\newblock \emph{IEEE Transactions on Information Theory}, 69\penalty0
  (2):\penalty0 1244--1260, 2023.
\newblock \doi{10.1109/TIT.2022.3213141}.
\newblock URL \url{https://doi.org/10.1109/TIT.2022.3213141}.

\bibitem[Neykov(2026)]{neykov2026polynomialtimenearoptimalestimationcertain}
M.~Neykov.
\newblock Polynomial-time near-optimal estimation over certain type-2 convex
  bodies, 2026.
\newblock URL \url{https://arxiv.org/abs/2512.22714}.

\bibitem[Nussbaum(1985)]{10.1214/aos/1176349651}
M.~Nussbaum.
\newblock {Spline Smoothing in Regression Models and Asymptotic Efficiency in
  $L_2$}.
\newblock \emph{The Annals of Statistics}, 13\penalty0 (3):\penalty0 984 --
  997, 1985.
\newblock \doi{10.1214/aos/1176349651}.
\newblock URL \url{https://doi.org/10.1214/aos/1176349651}.

\bibitem[Papapicco et~al.(2022)Papapicco, Demo, Girfoglio, Stabile, and
  Rozza]{PAPAPICCO2022114687}
D.~Papapicco, N.~Demo, M.~Girfoglio, G.~Stabile, and G.~Rozza.
\newblock The neural network shifted-proper orthogonal decomposition: A machine
  learning approach for non-linear reduction of hyperbolic equations.
\newblock \emph{Computer Methods in Applied Mechanics and Engineering},
  392:\penalty0 114687, 2022.
\newblock ISSN 0045-7825.
\newblock \doi{https://doi.org/10.1016/j.cma.2022.114687}.
\newblock URL
  \url{https://www.sciencedirect.com/science/article/pii/S004578252200069X}.

\bibitem[Pinkus(2012)]{pinkus2012n}
A.~Pinkus.
\newblock \emph{N-widths in Approximation Theory}, volume~7.
\newblock Springer Science \& Business Media, 2012.

\bibitem[Prasadan and Neykov(2024)]{prasadan2024characterizing}
A.~Prasadan and M.~Neykov.
\newblock Characterizing the minimax rate of nonparametric regression under
  bounded star-shaped constraints.
\newblock \emph{Electronic Journal of Statistics}, 2024.
\newblock \doi{10.1214/25-EJS2419}.

\bibitem[Prasadan and Neykov(2025)]{prasadan2024}
A.~Prasadan and M.~Neykov.
\newblock Some facts about the optimality of the lse in the gaussian sequence
  model with convex constraint.
\newblock \emph{arXiv}, 2025.
\newblock \doi{10.48550/arxiv.2406.05911}.
\newblock URL \url{https://doi.org/10.48550/arxiv.2406.05911}.

\bibitem[Prasadan and Neykov(2026)]{10.1214/25-AOS2576}
A.~Prasadan and M.~Neykov.
\newblock {Information theoretic limits of robust sub-Gaussian mean estimation
  under star-shaped constraints}.
\newblock \emph{The Annals of Statistics}, 54\penalty0 (1):\penalty0 490 --
  515, 2026.
\newblock \doi{10.1214/25-AOS2576}.
\newblock URL \url{https://doi.org/10.1214/25-AOS2576}.

\bibitem[Raskutti et~al.(2009)Raskutti, Wainwright, and Yu]{raskutti2009lower}
G.~Raskutti, M.~J. Wainwright, and B.~Yu.
\newblock Lower bounds on minimax rates for nonparametric regression with
  additive sparsity and smoothness.
\newblock pages 1563--1570, 2009.
\newblock Also available at: https://proceedings.neurips.cc/paper/2009.

\bibitem[Raskutti et~al.(2012)Raskutti, Wainwright, and
  Yu]{JMLR:v13:raskutti12a}
G.~Raskutti, M.~J. Wainwright, and B.~Yu.
\newblock Minimax-optimal rates for sparse additive models over kernel classes
  via convex programming.
\newblock \emph{Journal of Machine Learning Research}, 13\penalty0
  (13):\penalty0 389--427, 2012.
\newblock URL \url{http://jmlr.org/papers/v13/raskutti12a.html}.

\bibitem[Speckman(1985)]{10.1214/aos/1176349650}
P.~Speckman.
\newblock {Spline Smoothing and Optimal Rates of Convergence in Nonparametric
  Regression Models}.
\newblock \emph{The Annals of Statistics}, 13\penalty0 (3):\penalty0 970 --
  983, 1985.
\newblock \doi{10.1214/aos/1176349650}.
\newblock URL \url{https://doi.org/10.1214/aos/1176349650}.

\bibitem[Stone(1980)]{stone1980optimal}
C.~J. Stone.
\newblock Optimal rates of convergence for nonparametric estimators.
\newblock \emph{The annals of Statistics}, pages 1348--1360, 1980.

\bibitem[Talagrand(1996)]{10.1214/aop/1065725175}
M.~Talagrand.
\newblock {Majorizing measures: the generic chaining}.
\newblock \emph{The Annals of Probability}, 24\penalty0 (3):\penalty0 1049 --
  1103, 1996.
\newblock \doi{10.1214/aop/1065725175}.
\newblock URL \url{https://doi.org/10.1214/aop/1065725175}.

\bibitem[Talagrand(2001)]{10.1214/aop/1008956336}
M.~Talagrand.
\newblock {Majorizing measures without measures}.
\newblock \emph{The Annals of Probability}, 29\penalty0 (1):\penalty0 411 --
  417, 2001.
\newblock \doi{10.1214/aop/1008956336}.
\newblock URL \url{https://doi.org/10.1214/aop/1008956336}.

\bibitem[Talagrand(2005)]{talagrand2005generic}
M.~Talagrand.
\newblock \emph{The Generic Chaining: Upper and Lower Bounds of Stochastic
  Processes}.
\newblock Springer Monographs in Mathematics. Springer, Berlin, Heidelberg,
  2005.

\bibitem[Talagrand(2021)]{talagrand2021upper}
M.~Talagrand.
\newblock \emph{Upper and Lower Bounds for Stochastic Processes: Modern Methods
  and Classical Problems}.
\newblock Springer, Cham, Switzerland, 2nd edition, 2021.
\newblock ISBN 978-3-030-82595-9.
\newblock \doi{10.1007/978-3-030-82595-9}.
\newblock URL \url{https://link.springer.com/book/10.1007/978-3-030-82595-9}.

\bibitem[Tsybakov(2008)]{minimaxbook}
A.~B. Tsybakov.
\newblock \emph{Introduction to Nonparametric Estimation}.
\newblock Springer Publishing Company, Incorporated, 1st edition, 2008.
\newblock ISBN 0387790519.

\bibitem[Wainwright(2019)]{wainwright2019}
M.~J. Wainwright.
\newblock \emph{High-Dimensional Statistics: A Non-Asymptotic Viewpoint}.
\newblock Cambridge University Press, 2019.

\bibitem[Wichmann and Hill(2001)]{wichmann2001psychometric}
F.~A. Wichmann and N.~J. Hill.
\newblock The psychometric function: I. fitting, sampling, and goodness of fit.
\newblock \emph{Perception \& Psychophysics}, 63\penalty0 (8):\penalty0
  1293--1313, 2001.

\bibitem[Yang and Barron(1999)]{yang1999information}
Y.~Yang and A.~Barron.
\newblock Information-theoretic determination of minimax rates of convergence.
\newblock \emph{The Annals of Statistics}, 27\penalty0 (5):\penalty0
  1564--1599, 1999.
\newblock \doi{10.1214/aos/1017939142}.

\bibitem[Yang and Tokdar(2015)]{yang2015minimax}
Y.~Yang and S.~T. Tokdar.
\newblock Minimax-optimal nonparametric regression in high dimensions.
\newblock \emph{Annals of Statistics}, 43\penalty0 (2):\penalty0 652--674,
  2015.
\newblock \doi{10.1214/14-AOS1289}.

\bibitem[Zhang et~al.(2014)Zhang, Wainwright, and Jordan]{zhang2014lower}
Y.~Zhang, M.~J. Wainwright, and M.~I. Jordan.
\newblock Lower bounds on the performance of polynomial-time algorithms for
  sparse linear regression.
\newblock In \emph{Conference on Learning Theory}, pages 921--948. PMLR, 2014.

\end{thebibliography}
\newpage

\appendix

\section{Deferred Proofs}\label{section: deferred proofs}

\subsection{Proofs of Section \ref{section: theoretical results}}\label{subsection: proofs of theoretical results}

\subsubsection{Proof of Lemma \ref{lemma: bounds on KL divergence}}\label{subsubsection: proof of bounds on KL divergence}

\begin{proof}
The KL divergence between $f(y _{i};\theta _{i})$ and $f(y _{i};\theta _{i}^{\prime })$ is given by:
\begin{equation*}
\operatorname{KL}(\theta _{i}\left\lVert\right. \theta ^{\prime }_{i})=\mathbb{E}_{\theta _{i}}\left[\log \frac{f(y _{i};\theta _{i})}{f(y _{i};\theta ^{\prime }_{i})}\right],
\end{equation*}
where
\begin{equation*}
    \begin{aligned}
        \log \frac{f(y _{i};\theta _{i})}{f(y _{i};\theta ^{\prime }_{i})} &= \log \frac{h(y_i) \exp(\theta_i T(y_i) - A(\theta_i))}{h(y_i) \exp(\theta_i' T(y_i) - A(\theta_i'))} \\
        &= (\theta_i T(y_i) - A(\theta_i)) - (\theta_i' T(y_i) - A(\theta_i')) \\
        &= (\theta_i - \theta_i') T(y_i) + A(\theta_i') - A(\theta_i)
    \end{aligned}
\end{equation*}
Therefore,
\begin{equation*}
    \begin{aligned}
        \operatorname{KL}(\theta _{i}\left\lVert\right. \theta _{i}^{\prime })&= \mathbb{E}_{\theta_i} \left[ (\theta_i - \theta_i') T(y_i) + A(\theta_i') - A(\theta_i) \right]\\
        &= (\theta_i - \theta_i') A'(\theta_i) + A(\theta_i') - A(\theta_i),
    \end{aligned}
\end{equation*}
where we leverage the well-known property of the exponential family that $\mathbb{E}_{\theta _{i}}T(Y _{i})=A ^{\prime }(\theta _{i})$(\cite{Jordan2009}). By Taylor's theorem, we further have 
\begin{equation*}
    \operatorname{KL}(\theta _{i}\left\lVert\right. \theta _{i}^{\prime }) = (\theta_i - \theta_i') A'(\theta_i) + A(\theta_i') - A(\theta_i)= \frac{1}{2} (\theta_i - \theta_i')^2 A''(\tilde{\theta}),
\end{equation*}
where \(\tilde{\theta} = \alpha\theta_i + (1 - \alpha)\theta_i'\) for some \(\alpha \in [0, 1]\), assuming that \(A\) is twice continuously differentiable. By our assumption that our distribution from the exponential family is not degenerate, we know $A''(\tilde{\theta})=\text{Var}_{\tilde{\theta }}\left(T(Y _{i})\right)>0$ for all $\tilde{\theta }\in K$. Since we are also assuming that $\theta \in [-M,M]^{n}$, we know that there exist positive constants $c(M)$ and $C(M)$ such that $c(M)\le A _{i}^{\prime \prime }(\theta _{i})\le C(M)$ for any $-M \le \theta _{i}\le M$. Therefore, we come to the bound on $\operatorname{KL}(\theta _{i}\left\lVert\right. \theta _{i}^{\prime })$:
\begin{equation*}
    c(M) (\theta_i' - \theta_i)^2 \le \operatorname{KL}(\theta _{i}\left\lVert\right. \theta ^{\prime }_{i}) \le C(M) (\theta_i' - \theta_i)^2.
\end{equation*}
Consequently, we have:
\begin{equation}\label{eq: bounds on KL divergence}
    c(M) \left\lVert \theta ^{\prime }-\theta \right\rVert_{2}^{2} \le \operatorname{KL}(\theta \left\lVert\right. \theta ^{\prime })\le C(M)\left\lVert \theta ^{\prime }-\theta \right\rVert_{2}^{2}.
\end{equation}
This completes the proof.
\end{proof}

\subsubsection{Proof of Lemma \ref{lemma: minimax lower bound}}\label{subsubsection: proof of minimax lower bound}

\begin{proof}
From Fano's inequality (see Theorem \ref{theorem: fano inequality}), we have
\begin{equation*}
    \inf_{\hat{\theta }} \sup_{\theta _{0}\in K} \mathbb{E} \|\hat{\theta }(Y) - \theta _{0}\|^2 \ge \frac{\epsilon^2}{4} \left(1 - \frac{I(Y; J) + \log 2}{\log m}\right).
\end{equation*}
Let $\theta ^{*}\in K$ achieves $N(K \cap B _{2}(\theta ^{*},\epsilon ),\epsilon /c)=N ^{\text{loc}}(\epsilon ,c)$. By the property of the $\epsilon $-packing and $\epsilon $-covering, there exist $m=N ^{\text{loc}}(\epsilon ,c)$ points $\theta _{1},\dots,\theta _{m}$ in $K \cap B _{2}(\theta ^{*},\epsilon )$ forming an $\epsilon /c$ packing. Combining with \eqref{eq: bounds on KL divergence}, we know that
\begin{equation*}
    I(Y; J)=\frac{1}{m} \sum_{j=1}^{m} \operatorname{KL}(\theta _{j}\left\lVert\right. \frac{1}{m}\sum\limits _{i=1}^{m}\theta _{j})\le \frac{1}{m} \sum_{j=1}^{m} C(M) \left\lVert \theta _{j}-\theta ^{*}\right\rVert_{2}^{2}\le \max_j C(M) \|\theta _{j} - \theta ^{*}\|^2 \le C(M)\epsilon ^{2},
\end{equation*}
Therefore, by Fano's inequality we have
\begin{equation*}
    \inf_{\hat{\theta }} \sup_{\theta _{0}\in K} \mathbb{E} \|\hat{\theta }(Y) - \theta _{0}\|^2 \ge \frac{\epsilon^2}{4c^2} \left(1 - \frac{C(M) \epsilon^2 + \log 2}{\log N^{\text{loc}}(\epsilon,c)}\right).
\end{equation*}
By our assumption on $\log ^{\text{loc}}(\epsilon ,c)$, we have
\begin{equation*}
    \log N^{\text{loc}}(\epsilon,c) \ge 4 (\epsilon^2 C(M) \vee \log 2) \ge 2(C(M) \epsilon^2 + \log 2).
\end{equation*}
Therefore, we come to the conclusion:
\begin{equation*}
    \inf_{\hat{\theta }} \sup_{\theta _{0}\in K} \mathbb{E} \|\hat{\theta }(Y) - \theta _{0}\|^2 \ge \frac{\epsilon^2}{4 c^2} \left(1 - \frac{I(Y; J) + \log 2}{\log N ^{\text{loc}}(\epsilon ,c)}\right)\ge \frac{\epsilon^2}{4  c^2} \left(1 - \frac{1}{2}\right) =\frac{\epsilon^2}{8c^2}.
\end{equation*}
The proof is completed.
\end{proof}

\subsubsection{Proof of Lemma \ref{lemma: sub-exponential of T}}\label{subsubsection: proof of sub-exponential of T}

\begin{proof}
    Let us evaluate the moment generating function of $T(Y)$. By the properties of the exponential family, for any $\theta \in K$, we have $\mathbb{E}_{\theta }\left[\exp\left\{\lambda T(Y)\right\}\right]=\exp\left\{A(\theta +\lambda )-A(\theta )\right\}$ (see Theorem 18.3 of \cite{expfamily}). Therefore, we have 
    \begin{equation*}
        \begin{aligned}
            \log \left(\mathbb{E} _{\theta }\exp\left\{\lambda\cdot[T(Y) - \mathbb{E}_{\theta } T(Y)]\right\}\right) & = A(\theta + \lambda) - A(\theta ) - \lambda \cdot \mathbb{E}_{\theta } T(Y) \\
            & = A(\theta +\lambda )-A(\theta )-\lambda A ^{\prime }(\theta )\\ 
            & = \frac{1}{2}\lambda^2 A ^{\prime \prime }(\tilde{\theta }),
        \end{aligned}
    \end{equation*}
    where $\tilde \theta = \theta + (1-\alpha) \lambda \in [-2M,2M]$ for some $\alpha \in [0,1]$. Define $C ^{\prime }(M):=\sup\limits _{-2M \le \tilde{\theta }\le 2M}A ^{\prime \prime }(\tilde{\theta })$. Such $C ^{\prime }(M)$ is finite by the continuity of $A ^{\prime \prime }(\cdot )$.  Hence, whenever $\lambda \in [-M,M]$, we have
     \begin{equation*}
        \log \left(\mathbb{E} _{\theta }\exp\left\{\lambda\cdot[T(Y) - \mathbb{E}_{\theta } T(Y)]\right\}\right)\le \frac{1}{2}\lambda ^{2}C ^{\prime }(M).
    \end{equation*}
    This is equivalently to saying $T(Y)$ is sub-exponential. The proof is completed.
\end{proof}

\subsubsection{Proof of Lemma \ref{lemma: high probability for close distance between theta prime and theta prime prime}}\label{subsubsection: proof of high probability for close distance between theta prime and theta prime prime}

\begin{proof}   
    We start by observing the following identity from Appendix \ref{subsubsection: proof of bounds on KL divergence}:
   \begin{equation*}
    \sum\limits _{i=1}^{n} A(\theta_i'') - A(\theta_i') = \sum\limits _{i=1}^{n}\left[ (\theta_i'' - \theta_i') A'(\theta_i) + \operatorname{KL}(\theta_i \| \theta_i'') - \operatorname{KL}(\theta_i \| \theta_i') \right].
    \end{equation*}
    Therefore, we have
    \begin{equation*}
        \begin{aligned}
            \sum\limits _{i=1}^{n}(\theta_i' - \theta_i'')T(Y_i) + A(\theta_i'') - A(\theta_i')&=\sum\limits _{i=1}^{n}(\theta_i' - \theta_i'')T(Y_i) + \sum\limits _{i=1}^{n} \left[ A(\theta_i'') - A(\theta_i') \right] \\
            &=\sum\limits _{i=1}^{n}(\theta_i' - \theta_i'')T(Y_i) + \sum\limits _{i=1}^{n} \left[ (\theta_i'' - \theta_i') A'(\theta_i) + \operatorname{KL}(\theta_i \| \theta_i'') - \operatorname{KL}(\theta_i \| \theta_i') \right]\\
            &= \sum\limits _{i=1}^{n} (\theta_i' - \theta_i'') \left[T(Y_i) - \mathbb{E}T(Y_i) \right] + \operatorname{KL}(\theta \| \theta'') - \operatorname{KL}(\theta \| \theta'),
        \end{aligned}
    \end{equation*}
    where we recall that $\operatorname{KL}(\theta \left\lVert\right. \theta ^{\prime }):=\sum\limits _{i=1}^{n}\operatorname{KL}(\theta _{i}\left\lVert\right. \theta ^{\prime }_{i})$. By rearranging the terms, we have
    \begin{equation*}
        \mathbb{P}_\theta(\psi(Y) = 1) \le \mathbb{P}_\theta \left(\sum\limits _{i=1}^{n}(\theta_i' - \theta_i'')[T(Y_i) - \mathbb{E}T(Y_i)] \le \operatorname{KL}(\theta \| \theta') - \operatorname{KL}(\theta \| \theta'')\right)
    \end{equation*}
    Next, for any $\theta ^{\prime }: \left\lVert \theta ^{\prime }-\theta \right\rVert_{2}^{}\le \delta $ and $\theta ^{\prime \prime }: \left\lVert \theta ^{\prime \prime }-\theta ^{\prime }\right\rVert_{2}^{}\ge C\delta $, we have by Lemma \ref{lemma: bounds on KL divergence}: 
    \begin{equation*}
        \begin{aligned}
        \operatorname{KL}(\theta \| \theta') - \operatorname{KL}(\theta \| \theta'') & \le C(M) \|\theta - \theta'\|_{2}^2 - c(M)\|\theta - \theta''\|_{2}^2  \\
        & \le C(M) \|\theta - \theta'\|_{2}^2 - c(M) (\|\theta - \theta'\|_{2} - \|\theta'' - \theta'\|_{2})^2\\
        & \le C(M) \delta^2 - c(M) (\delta - (C\delta))^2\\ 
        &= C(M) \delta^2 - c(M) \cdot (C- 1)^2 \delta^2 < 0
        \end{aligned}
    \end{equation*}
    whenever $C>\sqrt{\frac{C(M)}{c(M)}}+1$. Since $T(Y)-\mathbb{E} T(Y)$ is sub-exponential from Lemma \ref{lemma: sub-exponential of T}, by Bernstein's inequality we have
    \begin{equation*}
        \mathbb{P}_\theta(\psi(Y) = 1) \le \exp\left\{-\min\limits \left\{\frac{(\operatorname{KL}(\theta \| \theta') - \operatorname{KL}(\theta \| \theta''))^2}{2C'(M)\|\theta' - \theta''\|^2},\frac{M|\operatorname{KL}(\theta \| \theta')-\operatorname{KL}(\theta \| \theta'')|}{2\|\theta' - \theta''\|_\infty}\right\}\right\}.
    \end{equation*}
    We know by assumption that $\|\theta' - \theta''\|_\infty \le 2 M$. Hence,
    \begin{equation*}
        \frac{M|\operatorname{KL}(\theta \| \theta') - \operatorname{KL}(\theta \| \theta'')|}{2\|\theta' - \theta''\|_\infty} \ge \frac{|C(M) - c(M)(C-1)^2|\delta^2}{4}
    \end{equation*}
    Next, using our bound from above
    \begin{equation*}
        \begin{aligned}
            \operatorname{KL}(\theta \| \theta') - \operatorname{KL}(\theta \| \theta'') & \le C(M) \|\theta - \theta'\|_{2}^2 - c(M) (\|\theta - \theta'\|_{2} - \|\theta'' - \theta'\|_{2})^2\\
            & \le C(M) \|\theta - \theta'\|_{2}^2 + 2c(M) \|\theta - \theta'\|_{2} \|\theta' - \theta''\|_{2} - c(M) \|\theta' - \theta''\|_{2}^2 \\
            & \le C(M)\delta^2 + 2c(M) \delta \|\theta' - \theta''\|_{2} - c(M)\|\theta' - \theta''\|_{2}^2\\
            & \overset{\text{(i)}}{\le }\left[C(M)+8c(M)\right]\delta ^{2}+\frac{1}{2}c(M)^{2}\left\lVert \theta ^{\prime }-\theta ^{\prime \prime }\right\rVert_{2}^{2}-c(M)\left\lVert \theta ^{\prime }-\theta ^{\prime \prime }\right\rVert_{2}^{2}\\ 
            & \overset{\text{(ii)}}{\le } - \frac{1}{4}c(M)\|\theta' - \theta''\|_{2}^2/4,
        \end{aligned}
    \end{equation*}
    where (i) is from the inequality $2ab \le a ^{2}+b ^{2}$ with $a=2 \sqrt{2c(M)}\delta $ and $b=\frac{\sqrt{2}}{2c(M)}\left\lVert \theta ^{\prime }-\theta ^{\prime \prime }\right\rVert_{2}^{}$, and (ii) is from the condition that $C \ge \sqrt{\frac{4(C(M)+8c(M))}{c(M)}}$. Hence, we conclude that
    \begin{equation*}
        \frac{(\operatorname{KL}(\theta \| \theta') - \operatorname{KL}(\theta \| \theta''))^2}{2C'(M)\|\theta' - \theta''\|^2} \ge \frac{c(M)^2\|\theta'  -\theta''\|^2}{32 C'(M)} \ge \frac{c(M)^2C^2 \delta^2}{32 C'(M)}.
    \end{equation*}
    This consequently leads to
    \begin{equation*}
        \begin{aligned}
            \mathbb{P}_\theta(\psi(Y) = 1) &\le \exp\left\{-\min\limits \left\{\frac{(\operatorname{KL}(\theta \| \theta') - \operatorname{KL}(\theta \| \theta''))^2}{2C'(M)\|\theta' - \theta''\|^2}, \frac{M|\operatorname{KL}(\theta \| \theta') - \operatorname{KL}(\theta \| \theta'')|}{2\|\theta' - \theta''\|_\infty}\right\}\right\}\\
            &\le \exp\left\{-\min\limits \left\{\frac{c(M)^2 C^2}{32 C'(M)}, \frac{|C(M) - c(M)(C-1)^2|}{4}\right\}\delta ^{2}\right\}.
        \end{aligned}
    \end{equation*}
    The proof is completed by setting $\kappa (M)=\min\limits \left\{\frac{c(M)^2 C^2}{32 C'(M)}, \frac{|C(M) - c(M)(C-1)^2|}{4}\right\}$.
\end{proof}

\begin{remark}\label{remark: comparison between kappa and C(M)}
If additionally we have $C \ge 1+\sqrt{\frac{2C(M)}{c(M)}}$, it is not hard to see that $\kappa (M)\ge \zeta C(M)$ with $\zeta =\frac{1}{8C ^{\prime }(M)c(M)}\wedge \frac{1}{4}$ with some simple algebra.
\end{remark}

\subsubsection{Proof of Lemma \ref{lemma: guarantee of the selection criterion}}\label{subsubsection: proof of guarantee of the selection criterion}

\begin{proof}
    Recall that $H(\delta ,\theta ^{\alpha },N _{k})$ is defined as, 
    \begin{equation*}
    H(\delta, \theta^\alpha, N _{k}) := \max \Big\{ \|\theta^\alpha - \theta^\beta\|_{2} \;\Big|\; \beta \in [N], \: \theta^\beta \succ \theta^\alpha, \text{ and } \|\theta^\beta - \theta^\alpha\|_{2} \ge C \delta \Big\},
    \end{equation*}
    where we additionally define $\max\limits \left\{\varnothing \right\}=0$. Without loss of generality, we assume that $\theta^\alpha \in N_k$ is the closest point to $\theta$ among all the candidates inside $N_k$, and thus $\|\theta^\alpha - \theta\| \le \delta$ since $N _{k}$ forms a $\delta $-covering of $K ^{\prime }$. According to the triangle inequality, we have 
    \begin{equation*}
        \| \overline{\theta} - \theta \|_{2} \le \| \overline{\theta} - \theta^\alpha \|_{2} + \| 
        \theta^\alpha - \theta  \|_{2}
    \end{equation*}
    Therefore, the fact that $\|\theta^\alpha - \theta\|_{2} \le \delta$, implies $\|\theta^\alpha - \overline{\theta}\|_{2} \ge C\delta$, and hence
    \begin{equation*}
        \mathbb{P} \left( \| \overline{\theta} - \theta \|_{2} > (C+1)\delta \right) \le \mathbb{P} \left( \overline{\theta} \in \{ \theta^\beta : \|\theta^\alpha - \theta^\beta\|_{2} \ge C\delta \} \right) \le \mathbb{P}(H(\delta, \theta^\alpha, N _{k}) > 0)
    \end{equation*}
    From Lemma \ref{lemma: high probability for close distance between theta prime and theta prime prime}, we know that
    \begin{equation*}
        \begin{aligned}
            \mathbb{P}(H(\delta, \theta^\alpha, S) > 0) & = \mathbb{P} \left(\exists \beta : \|\theta^\alpha - \theta^\beta\|_{2} \ge C\delta \text{ and } \sum\limits _{i=1}^{n}(\theta_i^\alpha - \theta_i^\beta)T(Y_i) + A(\theta_i^\beta) - A(\theta_i^\alpha) \le 0 \right) \\
            & \overset{\text{(i)}}{\le }N e^{-\kappa(M) \delta^2},
        \end{aligned}
    \end{equation*}
    where (i) is simply from the union bound. Thus, we conclude that
    \begin{equation*}
        \mathbb{P} \left( \| \overline{\theta} - \theta \| \ge (C+1)\delta \right) \le N e^{-\kappa(M) \delta^2}.
    \end{equation*}
    The proof is completed
\end{proof}

\subsubsection{Proof of Lemma \ref{lemma: property 1 of the estimators in the algorithm}}\label{subsubsection: proof of the property 1}

Before presenting the main proof, we introduce the following probabilistic lemma from \cite{10.1214/25-AOS2576} as a preliminary.

\begin{lemma}[Lemma B.5 from \cite{10.1214/25-AOS2576}]\label{lemma: probability lemma}
Let $J \ge 2$ be an integer and let $A_1, A_2, \ldots, A_J$ be a sequence of events. Then
\begin{equation*}
    \mathbb{P}(A_J) \le \mathbb{P}(A_1) + \sum_{j=2}^J \mathbb{P}(A_j \cap A_{j-1}^c).
\end{equation*}
\end{lemma}

\begin{proof}[Proof of Lemma~\ref{lemma: property 1 of the estimators in the algorithm}]
We first note that if $\tilde{J}$ satisfies \eqref{eq: definition of J}, then so does any $1 \le j \le \tilde{J}$ by the monotonicity of $\epsilon \mapsto N ^{\text{loc}}(\epsilon ,c)$. For $3 \le j \le \tilde{J}$, if $\|\Upsilon_{j-1} - \theta\| \le \frac{d}{2^{j-2}}$, then $\Upsilon_{j-1} = u$ for some $u \in \mathcal{L}(j-1) \cap B _{2}(\theta, d/2 ^{j-2})$ by definition. Therefore, by setting $\delta = \frac{d}{2^{j-1}(C+1)}$ in Lemma \ref{lemma: guarantee of the selection criterion} and applying the union bound, for $\Upsilon _{j}$ we have

\begin{equation*}
    \begin{aligned}
        &\mathbb{P}\left( \|\Upsilon_j - \theta\| > \frac{d}{2^{j-1}}, \|\Upsilon_{j-1} - \theta\| \le \frac{d}{2^{j-2}} \right) \\
        \le &\sum_{u \in \mathcal{L}(j-1) \cap B _{2}\left(\theta, d/2 ^{j-2}\right)} \mathbb{P}\left( \|\Upsilon_j - \theta\| > \frac{d}{2^{j-1}}, \Upsilon_{j-1} = u \right) \\
        = &\sum_{u \in \mathcal{L}(j-1) \cap B _{2}\left(\theta, d/2 ^{j-2}\right)} \mathbb{P} \left( \left\| \argmin_{\theta^\alpha \in \mathcal{O}(u)} H(\delta, \theta^\alpha, \mathcal{O}(u)) - \theta \right\| > (C+1)\delta, \Upsilon_{j-1} = u \right) \\
        \le &\sum_{u \in \mathcal{L}(j-1) \cap B _{2}\left(\theta, d/2 ^{j-2}\right)} \mathbb{P} \left( \left\| \argmin_{\theta^\alpha \in \mathcal{O}(u)} \left( H(\delta, \theta^\alpha, \mathcal{O}(u)) - \theta \right) \right\| > (C+1)\delta \right).
    \end{aligned}
\end{equation*}
The cardinality of  $\mathcal{L}(j-1) \cap \mathcal{B}(\theta, \frac{d}{2^{j-2}})$ is upper bounded by $N ^{\text{loc}}(d / 2^{j-2}, 2c)$ by Lemma \ref{lemma: bound on the cardinality}. Set $K' = B(u, d / 2^{j-2}) \cap K \subseteq K$ in Lemma \ref{lemma: guarantee of the selection criterion}, and recall from Lemma \ref{lemma: covering and packing properties of the constructed tree} that $\mathcal{O}(u)$ forms a $d / (2^{j-2}c) = d / (2^{j-1}(C+1))$-covering of $K'$ with cardinality bounded by $N ^{\text{loc}}(d / 2^{j-2}, 2c)$. Applying Lemma \ref{lemma: guarantee of the selection criterion} and noticing that the summands is constant in $u$, we can bound the probability term from above as:
\begin{equation*}
    \begin{aligned}
        \mathbb{P}\left( \|\Upsilon_j - \theta\| > \frac{d}{2^{j-1}}, \|\Upsilon_{j-1} - \theta\| \le \frac{d}{2^{j-2}} \right) &\le \left[ N^{\text{loc}} \left( \frac{d}{2^{j-2}}, 2c \right) \right]^2 \cdot e^{-\kappa(M) \epsilon_j^2}\\ 
        &= \left[ N^{\text{loc}} \left( 2(C+1)\epsilon_j, 2c \right) \right]^2 \cdot e^{-\kappa(M) \epsilon_j^2}
    \end{aligned}
\end{equation*}
where $\epsilon_j$ is also defined as $\epsilon_j \coloneqq \frac{d}{2^{j-1}(C+1)}.$

Define the event $A_j = \left\{\|\Upsilon_j - \theta\| > d / 2^{j-1}\right\}$. We have already bounded $\mathbb{P}(A_j \cap A_{j-1}^c)$ for $3 \le j \le \tilde{J}$ above. Note that $\mathbb{P}(A_1) = 0$ since $\Upsilon_1$ and $\theta$ both belong to a set of diameter $d$. Recall that for the level $2$, we construct a maximal $d/c$-packing of $B(\Upsilon_1,d) \cap K$, which is therefore also a $d/c$-covering (without pruning). So we may apply Lemma \ref{lemma: guarantee of the selection criterion} and conclude
\begin{equation*}
    \begin{aligned}
        \mathbb{P}(A_2 \cap A_1^c) = \mathbb{P}(A_2) &\le N_{\text{loc}}(d,c) \cdot e^{-\kappa(M) \epsilon_2^2} \\
        &\le \left[ N^{\text{loc}} \left( 2(C+1)\epsilon_2, 2c \right) \right]^2 \cdot e^{-\kappa(M) \epsilon_2^2}.
    \end{aligned}
\end{equation*}
Since $\epsilon _j$ is decreasing in $j$ and $N^{\text{loc}}(\epsilon ,c)$ is non-increasing in $\epsilon $, we may bound $N^{\text{loc}} \left(2(C+1)\epsilon_j, 2c \right)$ with $N^{\text{loc}} \left(2(C+1)\epsilon_J, 2c \right)$. Therefore, according to Lemma \ref{lemma: probability lemma}, for any $1 \le J \le \tilde{J}$ we have
\begin{equation*}
    \begin{aligned}
        \mathbb{P}(A_J) = \mathbb{P}\left( \|\Upsilon_J - \theta\| > \frac{d}{2^{J-1}} \right)  
        &\le \mathbb{P}(A_1) + \sum_{j=2}^J \mathbb{P}(A_j \cap A_{j-1}^c)\\
        &\le \sum_{j=2}^J \mathbb{P}(A_j \cap A_{j-1}^c)\\
        &\le  \left[ N^{\text{loc}} \left( 2(C+1)\epsilon_J, 2c \right) \right]^2 \sum_{j=2}^{J} e^{-\kappa(M) \epsilon_j^2} \\
        &\le \mathbf{1}_{\left\{J>1\right\}}\cdot N^{\text{loc}} \left[(2 (C+1) \epsilon_J, 2c \right)]^2 \cdot \frac{a_J}{1 - a_J},
    \end{aligned}
\end{equation*}
where we set $a_J = e^{-\kappa(M) \epsilon_J^2}$. Now suppose that (\ref{eq: definition of J}) holds. The fact that \( \kappa(M)\epsilon^2_J > \log 2 \) implies \( a_J < 1/2 \). Therefore, we conclude
\begin{equation*}
    \begin{aligned}
        \mathbb{P}\left( \|\Upsilon_J - \theta\| > \frac{d}{2^{J-1}} \right) 
        &\le \mathbf{1}_{\left\{J>1\right\}} \exp\left(\frac{\kappa(M) \epsilon^2_J}{2}\right) \cdot a_J \cdot \frac{1}{1 - a_J} \\
        &\le 2 \cdot \mathbf{1}_{\left\{J>1\right\}} \exp\left(\frac{\kappa(M) \epsilon^2_J}{2}\right) \cdot \exp\left(-\kappa(M) \epsilon_J^2\right) \\
        &= 2 \cdot \mathbf{1}_{\left\{J>1\right\}} \exp\left(-\frac{\kappa(M) \epsilon^2_J}{2}\right),
    \end{aligned}
\end{equation*}
which is exactly the desired inequality. The proof is completed.
\end{proof}

\subsubsection{Proof of Lemma \ref{lemma: property 2 of the estimators in the algorithm}}\label{subsubsection: proof of the property 2}

\begin{proof}
By Lemma \ref{lemma: property 1 of the estimators in the algorithm}, we have for $1 \le J \le J ^{*}$ that
\begin{equation*}
    \mathbbm{P}(\|\Upsilon_J - \theta\| > d/2^{J-1}) \le 2 \cdot \mathbf{1}_{\left\{J>1\right\}} \exp\left\{-\frac{\kappa(M)\epsilon^2_J}{2}\right\}.
\end{equation*}
Now recall our definition of $J^*$ and let $\theta^* = \Upsilon_{J^*}$ be the output after $J^* - 1$ steps. Define $B_J$ to be the event that $\|\theta - \theta^*\| > \omega \epsilon_J$ where $\omega :=\frac{7+6C}{2}.$ We upper bound $\mathbb{P}(B_J)$ under the condition \eqref{eq: definition of J}. By Lemma \ref{lemma: bound on the distance between the parent and child}, for any $1 \le J \le J ^{*}$ (which also implies $J \le J^*$) we have
\begin{equation*}
    \|\Upsilon_J - \theta^*\|_{2} \le \frac{d(2 + 4c)}{c2^J} \le \left( \frac{1}{C+1} + 4 \right) \cdot \frac{d}{2^J},
\end{equation*}
where we recall $c \ge 2(C+1)$. If $A _{J}$ fails to hold, i.e. $\|\Upsilon_J - \theta\|_{2} \le d/2^{J-1}$, then by the triangle inequality and the definitions of $\omega$ and $\epsilon_J$, we have
\begin{equation*}
    \|\theta^* - \theta\|_{2} \le \|\theta^* - \Upsilon_J\| + \|\Upsilon_J - \theta\| \le \left( \frac{1}{C+1} + 4 \right) \cdot \frac{d}{2^J} + \frac{d}{2^{J-1}} = \omega \epsilon_J.
\end{equation*}
This implies that $A _{J}^{\complement }\subset B _{J}^{\complement }$. Therefore, for $1 \le J \le \tilde{J}$,
\begin{equation*}
    \underbrace{\mathbb{P}\left( \|\theta - \theta^* \| > \omega \epsilon_J \right)}_{= \mathbb{P}(B_J)} \le \mathbb{P}(A_J) \le 2 \cdot \mathbf{1}_{\left\{J>1\right\}} \exp\left\{-\frac{\kappa(M)\epsilon_J^2}{2}\right\}.
\end{equation*}
In fact, this inequality holds for all integers $J \le 0$. For such $J$ we have $\omega \epsilon_J > \frac{d}{2^{J-2}} \ge 4d$ so both sides of the inequality are $0$. Since $\bigcup\limits _{-\infty <J \le \tilde{J}}^{} [\epsilon_J, \epsilon_{J-1}] = [\epsilon_{\tilde J}, \infty)$ and $\epsilon_J = \epsilon_{J-1}/2$, we observe that for any $x \ge \epsilon_{\tilde{J}}$, it must belong to an interval $[\epsilon_J, \epsilon_{J-1}]$ for some $J \le \tilde{J}$ and consequently satisfies $2\omega x \ge 2\omega \epsilon_J = \omega \epsilon_{J-1}$. Hence, for $x \ge \epsilon_{\tilde{J}}$,
\begin{equation*}
    \begin{aligned}
        \mathbb{P}(\|\theta - \theta^*\| > 2\omega x) &\le \mathbb{P}(\|\theta - \theta^*\| > \omega \epsilon_{J-1}) \\
        &\le 2 \cdot \mathbf{1}_{\left\{J>1\right\}} \exp\left\{-\frac{\kappa(M)\epsilon_J^2}{2}\right\}\\ 
        &\le 2 \cdot \mathbf{1}_{\left\{J ^{*}>1\right\}}
        \exp\left\{- \frac{\kappa(M)\epsilon_{J-1}^2}{2}\right\} \\
        &\le 2 \cdot \mathbf{1}_{\left\{J ^{*}>1\right\}}\exp\left\{-\frac{\kappa(M)x^2}{2}\right\},
    \end{aligned}
\end{equation*}
where we use the facts that $J \le J ^{\star }$ and $x \le \epsilon _{J-1}$.

We are interested in performing at least $J^*$ steps, so let $\theta^{**}$ be the output of $J^{**} \ge J^*$ steps. Then by Lemma \ref{lemma: bound on the distance between the parent and child}, for $x \ge \epsilon_{\tilde{J}}$ ($\tilde{J} \le J^*$ implies $\epsilon_{\tilde{J}} \ge \epsilon_{J^*}$), we have
\begin{equation*}
    \|\theta^* - \theta^{**}\|_{2} = \|\Upsilon_{J ^{*}
} - \Upsilon_{J^{**}}\|_{2} \le \frac{d(2 + 4c)}{c2^{J^*}} = \frac{5+4C}{7+6C}\epsilon _{J ^{*}}=\frac{5+4C}{7+6C}\omega \epsilon _{J ^{*}}\le \frac{5+4C}{7+6C}\omega x.
\end{equation*}
The triangle inequality implies
\begin{equation*}
    \|\theta - \theta^{**}\|_{2} \le \|\theta - \theta^*\|_{2} + \|\theta^* - \theta^{**}\|_{2} \le \|\theta - \theta^*\|_{2} + \frac{5 + 4C}{7 + 6C} \omega x.  
\end{equation*}
Set $\omega ^{\prime } = \left(2+\frac{5+4C}{7+6C}\right)\omega $. Then, for $x \ge \epsilon_{\tilde{J}}$, by the inequalities above we have
\begin{equation*}
    \begin{aligned}
        \mathbb{P}(\|\theta - \theta^{**}\|_{2} > \omega' x) 
        &\le \mathbb{P}\left( \|\theta - \theta^*\|_{2} + \frac{5 + 4C}{7 + 6C} \omega x > \omega' x \right) \\
        &= \mathbb{P}(\|\theta - \theta^*\|_{2} > 2 \omega x) \\
        &\le 2 \cdot \mathbf{1}_{\left\{J ^{*}>1\right\}}
        \exp\left\{-\frac{\kappa(M) x^2}{2}\right\}. 
    \end{aligned}
\end{equation*}
Therefore, we conclude that
\begin{equation*}
    \begin{aligned}
        \mathbb{E} _{Y} \|\theta - \theta^{**}(Y)\|_{2}^2 
        &= \int_0^{\infty} \mathbb{P}(\|\theta - \theta^{**}\|_{2}^2 > t) \, dt = \int_0^{\infty} 2t \cdot \mathbb{P}(\|\theta - \theta^{**}\| > t) \, dt \\
        &= 2\omega'^2 \int_0^{\infty} u \cdot \mathbb{P}(\|\theta - \theta^{**}\| > \omega' u) \, du \\
        &\le 2\omega'^2 \int_0^{\epsilon_{\tilde{J}}} u \, du + 2\omega'^2 \int_{\epsilon_{\tilde{J}}}^{\infty} u \cdot \mathbb{P}(\|\theta - \theta^{**}\| > \omega' u) \, du \\
        &\le \omega'^2 \epsilon_{\tilde{J}}^2 + \mathbf{1}_{\left\{\tilde{J}>1\right\}} \cdot 4\omega'^2 \cdot \int_{\epsilon_{\tilde{J}}}^{\infty} u \exp\left\{-\frac{\kappa(M)u^2}{2}\right\} du \\
        &= \omega'^2 \epsilon_{\tilde{J}}^2 + \mathbf{1}_{\left\{\tilde{J}>1\right\}} \cdot 4\omega'^2 \cdot \frac{1}{\kappa(M)}\exp\left\{-\frac{\kappa(M)\epsilon_{\tilde{J}}^2}{2}\right\}\\
        &= \frac{(16C+19)^{2}}{4}\epsilon_{\tilde{J}}^2 + \mathbf{1}_{\left\{\tilde{J}>1\right\}} \cdot (16C+19)^2 \cdot \frac{1}{\kappa(M)} \exp\left\{-\frac{\kappa(M)\epsilon_{\tilde{J}}^2}{2}\right\}.
    \end{aligned}
\end{equation*}
Setting $\tilde{J}=J ^{*}$ completes the proof.
\end{proof}

\subsubsection{Proof of Theorem \ref{theorem: performance of the theoretical algorithm}}\label{subsubsection: proof of the performance of the theoretical algorithm}

\begin{proof}
If $J^* = 1$, then $\epsilon_{J^*} \asymp d$ and $\epsilon_{J^*}^2 \wedge d^2 \asymp d^2$. Since $\theta ^{**},\theta _{0}\in K$, we trivially have $\mathbb{E}_Y \|\theta ^{**}(Y _{1},\dots,Y _{n}) - \theta _{0}\|_{2}^2 \le d^2$. Therefore, we assume $J^* > 1$ and consequently $\epsilon _{J ^{*}}\lesssim d$ in the following proof.

We first note that $J ^{*}$ must be finite since the left-hand side of \eqref{eq: definition of J} is an increasing function of $\epsilon_J$ (thus decreasing in $J$) while the right-hand side is a non-increasing function of $\epsilon_J$ (thus non-decreasing in $J$). By Lemma \ref{lemma: property 2 of the estimators in the algorithm}, we have 
\begin{equation*}
    \mathbb{E}_{Y}\|\theta _{0}- \theta^{**}(Y _{1},\dots,Y _{n})\|^2 \le \frac{(16C+19)^{2}}{4} \epsilon_{J ^{*}}^2 + \mathbf{1}_{\left\{J ^{*}>1\right\}} \cdot (16C+19)^{2}\cdot \frac{1}{\kappa(M)}\exp\left\{-\frac{\kappa(M)\epsilon_{J ^{*}}^2}{2}\right\}.   
\end{equation*}
Since $\kappa(M) \epsilon_{J^*}^2 \ge \log 2$ by definition of $J^*$, we have $\exp\left\{-\frac{\kappa(M)\epsilon_{\tilde{J}}^2}{2}\right\} \le \frac{1}{\sqrt{2}}$ and also $\frac{1}{\kappa(M)} < \frac{\epsilon_{J^*}^2}{\log 2}$. Thus, we further have
\begin{equation*}
    \begin{aligned}
        \mathbb{E}_{Y}\|\theta - \theta^{**}(Y _{1},\dots,Y _{n})\|_{2}^2 &\le \frac{(16C+19)^{2}}{4} \epsilon_{J ^{*}}^2 + \mathbf{1}_{\left\{J ^{*}>1\right\}} \cdot (16C+19)^{2}\cdot \frac{1}{\kappa(M)}\exp\left\{-\frac{\kappa(M)\epsilon_{J ^{*}}^2}{2}\right\} \\
        &\le \left(\frac{1}{4}+\frac{1}{\sqrt{2}\log 2}\right)(16C+19)^{2}\epsilon _{J ^{*}}^{2}.
    \end{aligned}
\end{equation*}
Combining with $\epsilon _{J ^{*}}\lesssim d$, we complete the proof.

\end{proof}

\subsubsection{Proof of Theorem \ref{theorem: main theorem for theoretical algorithm}}\label{subsubsection: proof of main theorem for theoretical algorithm}

\begin{proof}

Our proof is naturally divided into two cases based on the magnitude of $\epsilon^{*2} \kappa(M)$. Before detailing these cases, we first eliminate the degenerate scenario where $\epsilon^* = 0$. In this situation, $N^{\text{loc}}(\epsilon, c) = 1$ for any sufficiently small $\epsilon$, implying that the constraint set $K$ consists of a single point (i.e., $d = 0$). Consequently, our algorithm trivially outputs this exact point, achieving a minimax rate of 0. We may therefore assume without loss of generality that $\epsilon^* > 0$. This strict positivity implies $d > 0$, which guarantees that $N^{\text{loc}}(\epsilon, c)$ can be made arbitrarily large by choosing a sufficiently large constant $c$. With this established, we partition our analysis into the following two cases: $\epsilon^{*2} \kappa(M) > 8\log 2$ and $\epsilon^{*2} \kappa(M) \leq 8\log 2$.

\textbf{Case 1:}  $\epsilon^{*2} \kappa(M) > 8\log 2$. 

By the definition of $\epsilon ^{*}$, we always have $(\epsilon ^{*}-\gamma )^{2} \kappa(M) \le \log N^{\text{loc}}(\epsilon^*-\gamma ,c)$ for any $\gamma >0$. Defining $\delta^* := \epsilon^* [\sqrt{\frac{\kappa(M)}{8C(M)}} \wedge 1/2]\le \epsilon ^{*}/2$ and noticing that $N ^{\text{loc}}(\epsilon ,c)$ is non-increasing in $\epsilon $, we have
\begin{equation*}
    \begin{aligned}
        \log N^{\text{loc}}(\delta^*, c) &\ge \lim_{\gamma \to 0} \log N^{\text{loc}}(\epsilon^* - \gamma, c) \\ 
        &\ge \lim_{\gamma \to 0} [(\epsilon^* - \gamma)^2 \kappa(M)]\\
        &= \frac{\epsilon^{*2} \kappa(M)}{2} + \frac{\epsilon^{*2} \kappa(M)}{2}\\
        &\ge 4\delta^{*2} C(M) + 4\log 2\\
    \end{aligned}   
\end{equation*}
Therefore, the $\delta ^{*}$ specified above satisfies the condition of Lemma \ref{lemma: minimax lower bound}, and we can conclude that the minimax rate \eqref{eq: minimax risk of the estimation} is lower bounded by $\delta^* \asymp \epsilon^*$ up to some constant depending only on $M$.  

On the other hand, Theorem \ref{theorem: performance of the theoretical algorithm} establishes that $ \epsilon _{J ^{*}} \wedge d^2$ is an upper bound on the minimax rate \eqref{eq: minimax risk of the estimation}. If we can find a $\tilde{\epsilon} \asymp \epsilon^*$ such that $\tilde{\epsilon}^2 \ge \epsilon_{J^*}^{2}$, which further implies
\begin{equation*}
    \epsilon_{J^*}^{2} \wedge d^2 \le \tilde{\epsilon}^2 \wedge d^2,
\end{equation*}
then we know that actually $\epsilon_{J^*}^{2} \wedge d^2 \asymp \epsilon ^{*2}\wedge d ^{2}$ and consequently the upper bound and lower bound of the minimax rate \eqref{eq: minimax risk of the estimation} match.

Pick any fixed but universal constant $D > 1$, and let $\tilde{\epsilon} = 2D \epsilon^*$. Using the definition of $\epsilon^*$, we have
\begin{equation*}
    \kappa(M) \tilde{\epsilon}^2 = 4\kappa(M) (D \epsilon^*)^2 \ge 4 \log N^{\text{loc}}(D \epsilon^*,c) \ge 4 \log N^{\text{loc}} \left( Dc 2 \epsilon^*,c \right) =2 \log [N^{\text{loc}} \left( c \tilde{\epsilon},c \right)]^2,
\end{equation*}
where we assume $2c>1$ without loss of generality. By our assumption of the case 1, we also have
\begin{equation*}
    \kappa(M) \tilde{\epsilon}^2 = 4 \kappa(M) D^2 \epsilon^{*2} > \log 2.
\end{equation*}
Therefore, $\tilde{\epsilon}$ satisfies (\ref{eq: alternative definition of J}).

Suppose first that $\epsilon_J$ satisfies (\ref{eq: alternative definition of J}) for some $J \ge 1$. Then by the maximality of $J^*$ and the definition of $\epsilon _{J ^{*}}$, we know that $\epsilon_{J^* + 1} = \epsilon_{J^*}/2 \le \tilde{\epsilon }$. Thus $\epsilon _{J ^{*}}\le 2 \tilde{\epsilon }$. If no such $\epsilon_J$ exists, we must have $J^* = 1$, and trivially we have $\epsilon_{J^*}\le \tilde{\epsilon }$. We then conclude that we always have $\epsilon _{J ^{*}}\le 2\tilde{\epsilon }$. According to our previous argument, we know the minimax rate is characterized by \eqref{eq: minimax rate} in the case 1.

\textbf{Case 2:}  $\epsilon^{*2} \kappa(M) < 8\log 2$ 

In this case, we have $\epsilon^{*} \le \sqrt{\frac{8\log 2}{\kappa(M)}}$. Set $\epsilon =2\epsilon ^{*}>\epsilon ^{*}$. By the definition of $\epsilon^*$ as the supremum of $\{\epsilon: \epsilon^2 \kappa(M) \le \log N^{\text{loc}}(\epsilon,c)\}$, $\epsilon $ satisfies $\log N^{\text{loc}}(2\epsilon^*,c) < (2\epsilon^*)^2 \kappa(M) = 4\epsilon^{*2} \kappa(M) \le 32 \log 2$. We now claim that $d<8\epsilon^*$. Suppose not, i.e., $d\ge 8\epsilon^*$. Since $K$ is star-shaped, there must exist (by the triangle inequality) a line segment contained within $K$ with length at least $d/3 \ge 8\epsilon^*/3 > 2\epsilon^*$. Take $L^*$ as any line segment of length $2\epsilon^*$ that is entirely contained in $K$. Let $m$ be the midpoint of $L^*$. Then $L^* \subseteq B _{2}(m,\epsilon^*) \cap K \subseteq B _{2}(m,2\epsilon ^{*})\cap K$. Then partition \( L^* \)  into \(c\) sub-intervals of length \( 2\epsilon^*/c \). By choosing \( c \) large enough, we can produce a \( 2\epsilon^*/c \)-packing set of \( B _{2}(m,\epsilon ^{*})\cap K \) of cardinality exceeding \( \exp(32\log 2) \), violating our claim that \( \log N^{\text{loc}}(2\epsilon^*,c) \le 32 \log 2 \). Thus we conclude \( d \le 8\epsilon^* \le 8 \sqrt{\frac{8 \log 2}{\kappa(M)} } \), which indicates that the set $K$ is entirely inside a ball of radius $8\sqrt{\frac{8 \log 2}{\kappa(M)}}$ and $\kappa(M)d^2 \le 64\kappa(M)\epsilon^{*2} \le 512\log 2$.

Now take $\widetilde{\epsilon}=l d$ for any fixed $l \le \min\left\{\frac{1}{3},\sqrt{\frac{\zeta}{4}}\right\}$, where $\zeta$ is defined in Remark \ref{remark: comparison between kappa and C(M)}. Then $\kappa(M)\tilde{\epsilon }^{2}/l ^{2}= \kappa(M)d^2 \le 512\log 2$. Since $K$ contains a line segment of length at least $d/3$ and $l\le 1/3$, this line segment contains a subsegment $\widetilde{L}$ of length $\widetilde{\epsilon}=ld$. Let $\widetilde{m}$ be the midpoint of $\widetilde{L}$. Then $\widetilde{L}\subseteq B(\widetilde{m},\widetilde{\epsilon})\cap K$. Partition $\widetilde{L}$ into subintervals of length $\widetilde{\epsilon}/c$, we obtain a $\widetilde{\epsilon}/c$-packing set of $B(\widetilde{m},\widetilde{\epsilon})\cap K$ with cardinality at least $\lfloor c \rfloor+1$. Therefore, for sufficiently large $c$, $\log N^{\text{loc}}(\widetilde{\epsilon},c) \ge \log(\lfloor c \rfloor+1)>512\log 2$. Moreover, since $l^2\le \zeta/4$ and $\kappa(M)\ge \zeta C(M)$ by Remark \ref{remark: comparison between kappa and C(M)}, we have $4\widetilde{\epsilon}^2C(M)=4l^2d^2C(M)\le \zeta d^2C(M)\le \kappa(M)d^2\le512\log2\le \log N ^{\text{loc}}(\tilde{\epsilon },c)$. It follows that $\log N^{\text{loc}}(\widetilde{\epsilon},c)>4(\widetilde{\epsilon}^2C(M)\vee\log2)$. Consequently, $\tilde{\epsilon }$ satisfies the condition of Lemma \ref{lemma: minimax lower bound}, and the minimax rate \eqref{eq: minimax risk of the estimation} is lower bounded by $\widetilde{\epsilon}^2=l^2d^2 \asymp d ^{2}\asymp \epsilon ^{*2}\wedge d ^{2}$ up to some constants. For the upper bound, $d^2$ is trivially an upper bound on the minimax rate \eqref{eq: minimax risk of the estimation} and $\epsilon _{J ^{*}}^{2}\wedge d ^{2}\le d ^{2}$. Therefore, the lower bound given in Lemma \ref{lemma: minimax lower bound} and the upper bound given in Theorem \ref{theorem: performance of the theoretical algorithm} match except for a constant depending on $M$.

In conclusion, the minimax rate is given by $\epsilon^{*2}\wedge d^2$ up to some constant only depending on $M$ in the both cases. The proof is completed.
\end{proof}

\subsection{Proofs of Section \ref{section: efficient estimator}}\label{subsection: proofs of the efficient estimator}

\subsubsection{Proof of Corollary \ref{corollary: control of the distance between the weak projection and any point in K}}\label{subsubsection: proof of the corollary related to the weak projection}

\begin{proof}
    Without loss of generality we assume $P _{K}(\theta )=\mathop{\arg \min}\limits _{\tilde{\theta }\in K}\left\lVert \theta -\tilde{\theta }\right\rVert_{2}^{}$. By the property of the convex set, we know that 
    \begin{equation*}
        \left\lVert \tilde{P}_{K}(\theta )-P _{K}(\theta )\right\rVert_{2}^{2}+\left\lVert \theta -P _{K}(\theta )\right\rVert_{2}^{2}\le \left\lVert \theta -\tilde{P}_{K}(\theta )\right\rVert_{2}^{2}\overset{\text{(i)}}{\le }\left\lVert \theta -P _{K}(\theta )\right\rVert_{2}^{2}+\epsilon ^{2}+2\epsilon \left\lVert \theta -P _{K}(\theta )\right\rVert_{2}^{},
    \end{equation*}
    where (i) is from Theorem \ref{theorem: existence of the weak projection}. Hence, we have 
    \begin{equation*}
        \left\lVert \tilde{P}_{K}(\theta )-P _{K}(\theta )\right\rVert_{2}^{}\le \sqrt{\epsilon ^{2}+2\epsilon \left\lVert \theta -P _{K}(\theta )\right\rVert_{2}^{}}
    \end{equation*}
    By the triangle inequality, we have 
    \begin{equation*}
        \begin{aligned}
            \left\lVert \xi -\tilde{P}_{K}(\theta )\right\rVert_{2}^{}&\le \left\lVert \xi -P _{K}(\theta )\right\rVert_{2}^{}+\left\lVert \tilde{P}_{K}(\theta )-P _{K}(\theta )\right\rVert_{2}^{}\\ 
            &\le \left\lVert \theta -\xi \right\rVert_{2}^{}+\sqrt{\epsilon ^{2}+2\epsilon \left\lVert \theta -P _{K}(\theta )\right\rVert_{2}^{}}\\ 
            &\le \left\lVert \theta -\xi \right\rVert_{2}^{}+\sqrt{\epsilon ^{2}+2\epsilon \left\lVert \theta -\xi \right\rVert_{2}^{}}.
        \end{aligned}
    \end{equation*}
    The proof is completed.
\end{proof}

\subsubsection{Proof of Lemma \ref{lemma: property of our stochastic process}}\label{subsubsection: proof of the property of our stochastic process}

\begin{proof}
    From Lemma \ref{lemma: sub-Gaussian of T}, we know $T(Y _{i})$ is a sub-Gaussian random variable for any $i \in [n]$. Therefore, $R _{\theta } ^{(j)}$ has a sub-Gaussian tail with respect to the canonical $\ell _{2}$ norm on $E _{(j)}$. Next, by the majorizing measure theorem (\cite{talagrand2021upper}), we know that $\gamma _{2}(E _{(j)},\ell _{2})\asymp \mathcal{G}(E _{(j)})$, where the Gaussian width $\mathcal{G}(E _{(j)})$ can be bounded as:
    \begin{equation*}
        \begin{aligned}
            \gamma _{2}(E _{(j)},\ell _{2})\asymp \mathcal{G}(E _{(j)})&=\mathbb{E}\sup\limits _{\theta \in E _{(j)}}\theta ^{\top } \bm{g}\\ 
            &=\sup\limits _{\theta \in K \cap \hat{\theta }_{(j)}+2K _{(j)}}\sum\limits _{i=1}^{n}\left[\hat{\theta }_{(j),i}+A _{(j),i}^{\dagger \top }\left(\theta -\hat{\theta }_{(j)}\right)\right]g _{i}\\ 
            &=\mathbb{E}\sup\limits _{\theta \in K \cap \hat{\theta }_{(j)}+2K _{(j)}}\left(A _{(j)}^{\dagger }\theta \right)^{\top } \bm{g}\\ 
            &\overset{\text{(i)}}{\le }\mathbb{E}\sup\limits _{\theta \in 2K _{(j)}}\left(A _{(j)}^{\dagger }\theta \right)^{\top } \bm{g}\\
            &\overset{\text{(ii)}}{\le }\sup\limits _{\theta \in 2K _{(j)}}\left\lVert \theta \right\rVert_{2}^{}\mathbb{E}\left\lVert A _{(j)}^{\dagger }\bm{g}\right\rVert_{2}^{}\\ 
            &\overset{\text{(iii)}}{\le }\sup\limits _{\theta \in 2K _{(j)}}\left\lVert \theta \right\rVert_{2}^{}\sqrt{\mathbb{E}\left[\bm{g}^{\top } \left(\mathbf{I}_{n}-X _{(j)}^{\dagger }\right)\bm{g}\right]}\\ 
            &\le \tilde{d}_{(j)}\sqrt{\tilde{k} _{(j)}},
        \end{aligned}
    \end{equation*}
    where $\hat{\theta }_{(j),i},A _{(j),i}^{\dagger }$ represent the $i$-th coordinate of $\hat{\theta }_{(j)}$ and $i$-th row of $A _{(j)}^{\dagger }$, (i) is because $K \cap \hat{\theta }_{(j)}+2K _{(j)}\subset \hat{\theta }_{(j)}+2K _{(j)}$ and the translation property of the Gaussian width, (ii) is from the Cauchy--Schwarz's inequality, and (iii) is from the Jensen's inequality and the definition of $A _{(j)}^{\dagger }$. From the selection of $\tilde{k} _{(j)}$, we know that $\gamma _{2}(E _{(j)},\ell _{2})\lesssim \frac{\tilde{d}_{(j)}^{2}}{C}$.

    For $\Delta _{\ell _{2}}(E _{(j)})$, we trivially have 
    \begin{equation*}
        \begin{aligned}
            \Delta _{\ell _{2}}(T)&=\sup\limits _{\theta _{1},\theta _{2}\in E _{(j)}}\left\lVert \theta  _{1}-\theta _{2}\right\rVert_{2}^{}\\ 
            &=\sup\limits _{\theta _{1},\theta _{2}\in K \cap \hat{\theta }_{(j)}+2K _{(j)}}\left\lVert \hat{\theta }_{(j)}+A _{j}^{\dagger }\left(\theta _{1}-\hat{\theta }_{(j)}\right)-\left(\hat{\theta }_{(j)}+A _{j}^{\dagger }\left(\theta _{2}-\hat{\theta }_{(j)}\right)\right)\right\rVert_{2}^{}\\
            &=\sup\limits _{\theta _{1},\theta _{2}\in K \cap \hat{\theta }_{(j)}+2K _{(j)}}\left\lVert A _{(j)}^{\dagger }\left(\theta _{1}-\theta _{2}\right)\right\rVert_{2}^{}\\ 
            &\le \tilde{d}_{(j)}.
        \end{aligned}
    \end{equation*}
    The proof is completed by selecting $C \asymp c ^{2}$ and $u \asymp \frac{\tilde{d}_{(j)}^{2}}{c ^{4}}$ (from the loop condition, we know $u \ge 1$) in Lemma \ref{lemma: moment and tail bounds of sp}.
\end{proof}

\subsubsection{Proof of Theorem \ref{theorem: exponentially decreasing distance between the estimator and true parameter}}\label{subsubsection: proof of theoretical guarantee of the efficient algorithm}

\begin{proof}
    Denote $\theta _{\text{m}}^{(j+1)}$ as the true maximizer of the log-likelihood parameterized by $\hat{\theta }_{(j)}+2K _{(j)}$ in the $j$-th iteration. Let $\theta ^{\prime }=\hat{\theta }_{(j)}+A _{(j)}^{\dagger }\left(\theta _{0}-\hat{\theta }_{(j)}\right)$ and $\theta ^{\prime \prime }=\hat{\theta }_{(j)}+A _{(j)}^{\dagger }\left(\theta _{\text{m}}^{(j+1)}-\hat{\theta }_{(j)}\right)$, where $A _{(j)}^{\dagger }$ from Algorithm \ref{algorithm: main algorithm} is guaranteed to be valid with probability greater than $1-\frac{1}{1+\tilde{d}_{(j)}^{3}}$ by Theorem \ref{theorem: polynomial time approximation to the SDP problem}. Since $\theta _{\text{m}}^{(j+1)}$ is the maximizer, we have 
    \begin{equation*}
        \sum\limits _{i=1}^{n}\theta ^{\prime }_{i}T(y _{i})-A(\theta ^{\prime }_{i})\le \sum\limits _{i=1}^{n}\theta ^{\prime \prime }_{i}T(y _{i})-A(\theta ^{\prime \prime }_{i})
    \end{equation*}
    From the proof of Lemma \ref{lemma: high probability for close distance between theta prime and theta prime prime}, this is equivalent to 
    \begin{equation}\label{eq: equivalent expression for the maximizer}
        \sum\limits _{i=1}^{n}\left(\theta ^{\prime }_{i}-\theta ^{\prime \prime }_{i}\right)\left[T(y _{i})-\mathbb{E}\left[T(y _{i})\right]\right]\le \operatorname{KL}\left(\theta _{0}\left\lVert\right. \theta ^{\prime }\right)-\operatorname{KL}\left(\theta _{0}\left\lVert\right. \theta ^{\prime \prime }\right).
    \end{equation}
    From \eqref{eq: bounds on KL divergence}, we have
    \begin{equation*}
        c _{M}\left\lVert \theta _{0}-\theta ^{\prime }\right\rVert_{2}^{2}\le \operatorname{KL}(\theta _{0}\left\lVert\right. \theta ^{\prime })\le C _{M}\left\lVert \theta _{0}-\theta ^{\prime \prime }\right\rVert_{2}^{2}.
    \end{equation*}
    By our assumption on $\hat{\theta }_{(j)}$ and the property of $A _{(j)}^{\dagger }$, we know 
    \begin{equation*}
        \operatorname{KL}\left(\theta _{0}\left\lVert\right. \theta ^{\prime }\right)\le C _{M}\left\lVert \left(\mathbf{I}_{n}-A _{(j)}^{\dagger }\right)(\theta _{0}-\hat{\theta } _{(j)})\right\rVert_{2}^{2}= C _{M}\left(\theta _{0}-\hat{\theta } _{(j)}\right)^{\top } X _{(j)}^{\dagger }\left(\theta _{0}-\hat{\theta }_{(j)}\right)\overset{\text{(i)}}{\le } C _{M}\left(\frac{\tilde{d} _{(j)}}{c}\right)^{2},
    \end{equation*}
    where (i) is from Theorem \ref{theorem: polynomial time approximation to the SDP problem}.

    In the same way, we have $\operatorname{KL}(\theta _{0}\left\lVert\right. \theta ^{\prime \prime })\ge c _{M}\left\lVert \theta _{0}-\theta ^{\prime \prime }\right\rVert_{2}^{2}$. By Lemma \ref{lemma: property of our stochastic process}, we know (universal constants omitted) with probability greater than $1-\exp\left\{-\frac{\tilde{d}_{(j)}^{2}}{c ^{4}}\right\}$, we have 
    \begin{equation*}
        \sup\limits _{\theta \in E _{(j)}}\left|R _{\theta } ^{(j)}-R _{\theta ^{\prime }}\right|\le C _{M}\left(\frac{\tilde{d}_{(j)}}{c}\right)^{2},
    \end{equation*}
    where we recall that $R _{\theta } ^{(j)}:=\sum\limits _{i=1}^{n}\theta _{i}\left(T(Y _{i})-\mathbb{E}\left[T(Y _{i})\right]\right), \theta \in E _{(j)}$. Conditioning on such event, we have 
    \begin{equation*}
        c _{M}\left\lVert \theta _{0}-\theta ^{\prime \prime }\right\rVert_{2}^{2}\le \operatorname{KL}\left(\theta _{0}\left\lVert\right. \theta ^{\prime \prime }\right)\overset{\text{(i)}}{\le }\operatorname{KL}\left(\theta _{0}\left\lVert\right. \theta ^{\prime }\right)-\left(R _{\theta ^{\prime }}-R _{\theta ^{\prime \prime }}\right)\le 2C _{M}\left(\frac{\tilde{d}_{(j)}}{c}\right)^{2}.
    \end{equation*}
    Therefore, we have $\left\lVert \theta _{0}-\theta ^{\prime \prime }\right\rVert_{2}^{}\le \sqrt{\frac{2C _{M}}{c _{M}}}\frac{\tilde{d}_{(j)}}{c}$. Since the accuracy of the ellipsoid method is $\sqrt{\frac{2C _{M}}{c _{M}}}\frac{\tilde{d}_{(j)}}{c}$ as well, we know that 
    \begin{equation*}
        \left\lVert \theta _{0}-\left(\hat{\theta }_{(j)}+A ^{\dagger }_{(j)}\left(\hat{\theta } _{\text{m}}^{(j+1)}-\hat{\theta }_{(j)}\right)\right)\right\rVert_{2}^{}\le \left\lVert \theta _{0}-\theta ^{\prime \prime }\right\rVert_{2}^{}+\left\lVert A _{(j)}^{\dagger }\left(\theta _{\text{m}}^{(j+1)}-\hat{\theta }_{\text{m}}^{(j+1)}\right)\right\rVert_{2}^{}\le \sqrt{\frac{8C _{M}}{c _{M}}}\frac{\tilde{d}_{(j)}}{c}
    \end{equation*}
    as $\mathbf{0}\preceq A _{(j)}^{\dagger }\preceq \mathbf{I}_{n}$. Setting $\epsilon =\left\lVert \theta _{0}-\left(\hat{\theta }_{(j)}+A ^{\dagger }_{(j)}\left(\hat{\theta } _{\text{m}}^{(j+1)}-\hat{\theta }_{(j)}\right)\right)\right\rVert_{2}^{}$ in Corollary \ref{corollary: control of the distance between the weak projection and any point in K}, we find that 
    \begin{equation*}
        \begin{aligned}
            \left\lVert \theta _{0}-\hat{\theta }_{(j+1)}\right\rVert_{2}^{}&=\left\lVert \theta _{0}-\tilde{P}_{K}\left[\hat{\theta }_{(j)}+A ^{\dagger }_{(j)}\left(\hat{\theta } _{\text{m}}^{(j+1)}-\hat{\theta }_{(j)}\right)\right]\right\rVert_{2}^{}\\ 
            &\le \left(\sqrt{3}+1\right)\left\lVert \theta _{0}-\left(\hat{\theta }_{(j)}+A ^{\dagger }_{(j)}\left(\hat{\theta } _{\text{m}}^{(j+1)}-\hat{\theta }_{(j)}\right)\right)\right\rVert_{2}^{}\\ 
            &\le \sqrt{\frac{8C _{M}}{c _{M}}}\frac{\left(\sqrt{3}+1\right)\tilde{d}_{(j)}}{c}.
        \end{aligned}
    \end{equation*}
    Aggregating the error probability as $\frac{1}{1+\tilde{d}_{(j)}^{3}}+\exp\left\{-\frac{\tilde{d}_{(j)}^{2}}{c ^{4}}\right\}$, the proof is completed.
\end{proof}

\subsubsection{Proof of Lemma \ref{lemma: relationship between minimax rate and r}}\label{subsubsection: proof of relationship between minimax rate and r}

\begin{proof}
    The proof is directly from the definition of $\epsilon ^{*}$. Recall that $\epsilon ^{*}:=\sup\limits _{\epsilon }\left\{\epsilon ^{2}\kappa (M)\le \log N ^{\text{loc}}(\epsilon ,c)\right\}$. If $r \gtrsim \sqrt{n}$, from the packing number of the Euclidean ball, we know that $\log N ^{\text{loc}}(\sqrt{n},c)\asymp n$. Therefore, for $\epsilon \asymp \sqrt{n}$, we have $\epsilon ^{2}\asymp \log N ^{\text{loc}}(\epsilon ,c)$. Since $\epsilon ^{2}\kappa (M)$ is increasing in $\epsilon $ while $\log N ^{\text{loc}}(\epsilon ,c)$ is non-increasing in $\epsilon $, we know that $\epsilon ^{*}\asymp \sqrt{n}$.

    On the other side, if $r \lesssim \sqrt{n}$, similarly, for $\epsilon \asymp r$, we have $\epsilon ^{2}\lesssim n \asymp \log N ^{\text{loc}}(r,c)$. Therefore, we have $\epsilon ^{*}\gtrsim r$ by definition.
\end{proof}

\subsubsection{Proof of Lemma \ref{lemma: relationship between epsilon and diam}}\label{subsubsection: proof of relationship between epsilon and diam}

\begin{proof}
    We separate the cases for $\text{diam}(K) \ge 1$ and $\text{diam}(K) <1$.

    Assume $\text{diam}(K) \ge 1$. By the definition of $\text{diam}(K) $, there exists a line segment $S \subset K$ with length $1 \le \text{diam}(K) $. Therefore, take $\epsilon _{0}=\frac{1}{2} \wedge \frac{1}{\sqrt{\kappa _{M}}}$, we have $N ^{\text{loc}}(\epsilon _{0})\ge \sup\limits _{\nu \in K}N \left(\frac{\epsilon _{0}}{c},B _{2}(\nu ,\epsilon _{0})\cap S\right)\ge \frac{2\epsilon _{0}}{2\epsilon _{0}/c}=c$, and hence $\log N ^{\text{loc}}(\epsilon )\ge 1$. Since $\epsilon _{0}^{2}\kappa _{M}\le 1$, we know that $\epsilon ^{*}\ge \epsilon _{0}\asymp _{M}1$.

    Assume $\text{diam}(K)<1$. Take $\epsilon _{0}=\frac{\text{diam}(K) }{2(\sqrt{\kappa _{M}}\vee 1)}$. By the same argument, we know $\log N ^{\text{loc}}(\epsilon _{0})\ge 1$. Since $\epsilon _{0}^{2}\kappa _{M}=\frac{\text{diam}(K) ^{2}}{4(1 \vee \kappa _{M})}\kappa _{M}\le \frac{\text{diam}(K) ^{2}}{4}<1$, we know $\epsilon ^{*}\ge \epsilon _{0}\asymp _{M}\text{diam}(K) $. The proof is completed.
\end{proof}

\subsubsection{Proof of Lemma \ref{lemma: matching the minimax rate under ideal assumptions}}\label{subsubsection: proof of matching the minimax rate under ideal assumptions}

\begin{proof}
    The proof is directly from Lemmas \ref{lemma: relationship between minimax rate and r} and \ref{lemma: relationship between epsilon and diam}.

    First suppose that $r \gtrsim \sqrt{n}$. By Algorithm \ref{algorithm: main algorithm}, an arbitrary $\hat{\theta }\in K$ is the output. By Lemma \ref{lemma: relationship between minimax rate and r}, we know that $\epsilon ^{*}\asymp \sqrt{n}$. Since $K \in [-M,M]^{n}$, we trivially have $\left\lVert \theta _{0}-\hat{\theta }\right\rVert_{2}^{}\le M \sqrt{n}\lesssim _{M}\sqrt{n}$, which means that $\hat{\theta }$ achieves the minimax rate under this case.

    Now suppose $r \le \sqrt{n}$, which implies that $\epsilon ^{*}\gtrsim r$. We further split this case into two sub-cases. If $r \ge c ^{2}$, Algorithm \ref{algorithm: main algorithm} stops whenever $\tilde{d}_{(j)}\le r$. With the ideal assumptions, for $\hat{\theta }=\hat{\theta }_{(j)}$, by Theorem \ref{theorem: exponentially decreasing distance between the estimator and true parameter} and Lemma \ref{lemma: relationship between minimax rate and r} we have $\left\lVert \theta _{0}-\hat{\theta }\right\rVert_{2}^{}\le \tilde{d}_{(j)}\le r \lesssim \epsilon ^{*}$, which indicates that $\hat{\theta }$ achieves the minimax rate. Otherwise, we have $r \le c ^{2}$ and consequently $\left\lVert \theta _{0}-\hat{\theta }\right\rVert_{2}^{}\lesssim c ^{2}$. By Lemma \ref{lemma: relationship between epsilon and diam}, if $\text{diam}(K)\gtrsim 1$, then $\epsilon ^{*}\gtrsim _{M}1 \gtrsim \left\lVert \theta _{0}-\hat{\theta }\right\rVert_{2}^{}$; if $\text{diam}(K) \lesssim 1$, then $\epsilon ^{*}\gtrsim \text{diam}(K) $. However, we trivially have $\text{diam}(K) \gtrsim \left\lVert \theta _{0}-\hat{\theta }\right\rVert_{2}^{}$ by the definition of $\text{diam}(K) $. The discussions above demonstrates that $\hat{\theta }$ always achieves the minimax rate when $X _{(j)}^{\dagger }$ are always valid and the supreme of the stochastic process is well controlled. The proof is completed.
\end{proof}

\subsubsection{Proof of Lemma \ref{lemma: maximum error of the subsequent estimator}}\label{subsubsection: proof of maximum error of the subsequent estimator}

\begin{proof}
    Since for any $j \ge j _{0}+1$, $X _{(j)}^{\dagger }$ is no longer valid, we should prove the argument from the statement of the algorithm. We start with $\hat{\theta }_{j _{0}}$. Since it is the last step on the normal trajectory, we have $\left\lVert \hat{\theta }_{(j _{0})}-\theta _{0}\right\rVert_{2}^{}\le \tilde{d}_{(j)}$. Consider $\left\lVert \hat{\theta }_{(j+1)}-\theta _{0}\right\rVert_{2}^{}$, by definition we have 
    \begin{equation*}
        \begin{aligned}
            \left\lVert \hat{\theta }_{j+1}-\theta _{0}\right\rVert_{2}^{}&=\left\lVert \Pi _{K}\left[\hat{\theta }_{(j)}+A _{(j)}^{\dagger }\left(\hat{\theta }_{\text{m}}^{(j+1)}-\hat{\theta }_{(j)}\right)\right]-\theta _{0}\right\rVert_{2}^{}\\ 
            &\overset{\text{(i)}}{\le }\left\lVert \hat{\theta }_{(j)}+A _{j}^{\dagger }\left(\hat{\theta }_{\text{m}}^{(j+1)}-\hat{\theta }_{(j)}\right)-\theta _{0}\right\rVert_{2}^{}+\sqrt{\epsilon ^{2}+2\epsilon \left\lVert \hat{\theta }_{(j)}+A _{j}^{\dagger }\left(\hat{\theta }_{\text{m}}^{(j+1)}-\hat{\theta }_{(j)}\right)-\theta _{0}\right\rVert_{2}^{}}\\ 
            &\overset{\text{(ii)}}{\le }\left\lVert \hat{\theta }_{(j)}-\theta _{0}\right\rVert_{2}^{}+\left\lVert \hat{\theta }_{\text{m}}^{(j+1)}-\hat{\theta }_{(j)}\right\rVert_{2}^{}+\sqrt{\epsilon ^{2}+2\epsilon \left[\left\lVert \hat{\theta }_{(j)}-\theta _{0}\right\rVert_{2}^{}+\left\lVert \hat{\theta }_{\text{m}}^{(j+1)}-\hat{\theta }_{(j)}\right\rVert_{2}^{}\right]}\\ 
            &\overset{\text{(iii)}}{\le }\left\lVert \hat{\theta }_{(j)}-\theta _{0}\right\rVert_{}^{}+\tilde{d}_{(j)}+\epsilon +\sqrt{2\epsilon \left(\tilde{d}_{(j)}+\left\lVert \hat{\theta }_{(j)}-\theta _{0}\right\rVert_{2}^{}\right)},
        \end{aligned}
    \end{equation*}
    where (i) is from Corollary \ref{corollary: control of the distance between the weak projection and any point in K}, (ii) is from the triangle inequality and $\mathbf{0}\preceq A _{(j)}^{\dagger }\preceq \mathbf{I}_{n}$ (this always holds no matter whether $X _{(j)}^{\dagger }$ is valid or not), and (iii) is from that $\hat{\theta }_{\text{m}}^{(j+1)}\in \hat{\theta }_{(j)}+2K _{(j)}$, hence $\hat{\theta }_{\text{m}}^{(j+1)}-\hat{\theta }_{(j)}\in 2K _{(j)}=2K \cap B _{2}(0,\tilde{d}_{(j)})$, and consequently $\left\lVert \hat{\theta }_{\text{m}}^{(j+1)}-\hat{\theta }_{(j)}\right\rVert_{2}^{}\le \tilde{d}_{(j)}$. Now taking $\epsilon =\tilde{d}_{(j)}$, we further have 
    \begin{equation*}
        \begin{aligned}
            \left\lVert \hat{\theta }_{(j+1)}-\theta _{0}\right\rVert_{2}^{}&\le \left\lVert \hat{\theta }_{(j)}-\theta _{0}\right\rVert_{2}^{}+2 \tilde{d}_{(j)}+\sqrt{2 \tilde{d}_{(j)}\left(\tilde{d}_{(j)}+\left\lVert \hat{\theta }_{(j)}-\theta _{0}\right\rVert_{2}^{}\right)}\\ 
            &\le \left\lVert \hat{\theta }_{(j)}-\theta _{0}\right\rVert_{2}^{}+\left(2+\sqrt{2}\right)\tilde{d}_{(j)}+\sqrt{2 \tilde{d}_{(j)}\left\lVert \hat{\theta }_{(j)}-\theta _{0}\right\rVert_{2}^{}}.
        \end{aligned}
    \end{equation*}
    Define $b _{j}:=\left\lVert \hat{\theta }_{(j)}-\theta _{0}\right\rVert_{2}^{}$, we thus have 
    \begin{equation*}
        \begin{aligned}
            &b _{j+1}\le b _{j}+\left(2+\sqrt{2}\right)\tilde{d}_{(j)}+\sqrt{2 \tilde{d}_{(j)}b _{j}}\le \left(\sqrt{b _{j}}+\sqrt{2+\sqrt{2}}\sqrt{\tilde{d}_{(j)}}\right)^{2}\\ 
            \Rightarrow &\sqrt{b _{j+1}}\le \sqrt{b _{j}}+\sqrt{2+\sqrt{2}}\sqrt{\tilde{d}_{(j)}}.
        \end{aligned}
    \end{equation*}
    Summating the inequality above from $j _{0}$ to arbitrary $j \ge j _{0}+1$, we have 
    \begin{equation*}
        \sqrt{b _{j}}\le \sqrt{b _{j _{0}}}+\frac{\sqrt{2+\sqrt{2}}}{1-\sqrt{q}}\tilde{d}_{j _0}\overset{\text{(i)}}{\le }\left(1+\frac{\sqrt{2+\sqrt{2}}}{1-\sqrt{q}}\right)\sqrt{\tilde{d}_{j _{0}}}, \forall j \ge j _{0}+1,
    \end{equation*}
    where $q=\frac{\tilde{d}_{(j+1)}}{\tilde{d}_{(j)}}=\sqrt{\frac{8C _{M}}{c _{M}}}\frac{\sqrt{3}+1}{c}$ and (i) is from that $\left\lVert \hat{\theta }_{(j _{0})}-\theta _{0}\right\rVert_{2}^{}\le \tilde{d}_{(j _{0})}$. The proof is completed.
\end{proof}

\subsubsection{Proof of Theorem \ref{theorem: main theorem for the algorithm}}\label{subsubsection: proof of main theorem for the algorithm}

\begin{proof}
    When $\text{diam}(K) \le 1$, by Lemma \ref{lemma: relationship between epsilon and diam} we know $\epsilon ^{*}\gtrsim _{M}\text{diam}(K) $, hence \eqref{eq: main inequality for the estimator} trivially holds. Therefore, we assume $\text{diam}(K) \ge 1$, and consequently $\epsilon ^{*}\gtrsim _{M} 1$ in the following analysis.
    
    In the ideal case where $X _{(j)}^{\dagger }$ is always valid and $R _{\theta }^{(j)}$ is always well controlled, we know by Lemma \ref{lemma: matching the minimax rate under ideal assumptions} that $\hat{\theta }_{\text{p}}(Y)$ exactly achieves the minimax rate.

    There are three possibilities that could break the ideal assumptions. First, the modified ellipsoid method to compute $X _{(j)}^{\dagger }$ by \cite{10.1145/3406325.3451128,neykov2026polynomialtimenearoptimalestimationcertain} could fail, whose probability is smaller than $\frac{1}{1+\tilde{d}_{(j)}^{3}}$ by Theorem \ref{theorem: polynomial time approximation to the SDP problem}. Second, the conditions required by Theorem \ref{theorem: polynomial time approximation to the SDP problem}, i.e., $\tilde{d}_{(j)}\ge 2 \sqrt{\tau _{K _{(j)}}T _{2}(K _{(j)})}\epsilon ^{*}\vee c ^{2} \sqrt{\lceil 4 \lceil \log _{c}(\sqrt{\tau _{K _{(j)}}}T _{2}(K _{(j)})) \rceil \epsilon ^{\star 2} \rceil }$ might not be satisfied after some iterations. Third, the behavior of $R _{\theta }^{(j)}$ might fall into a low probability region where the supreme is not bounded within the desired order. The probability of such region is smaller than $\exp\left\{-\frac{\tilde{d}_{(j)}^{2}}{c ^{4}}\right\}$ in $j$-th iteration by Theorem \ref{theorem: exponentially decreasing distance between the estimator and true parameter}. By the definition of the Minkowski gauge, the type-2 constant and $K _{(j)}$, it is not difficult to show that $T _{2}(K _{(j)})\le \sqrt{2}T _{2}(K), \forall j$. Define $\tau _{\text{m}}:=\max\limits _{j \ge 0}\tau _{K _{(j)}}$. To simplify the analysis, we only consider the $\ell _{2}$ error that is greater than $\tilde{d}_{\text{e}}:=2 \sqrt{\tau _{\text{m}}T _{2}(K _{(j)})}\epsilon ^{*}\vee c ^{2} \sqrt{\lceil 4 \lceil \log _{c}(\sqrt{\tau _{\text{m}}}T _{2}(K _{(j)})) \rceil \rceil }\epsilon ^{*}$. By Lemma \ref{lemma: maximum error of the subsequent estimator}, the case that the $\ell _{2}$ error greater than $\tilde{d}_{\text{e}}$ could only happen where the ideal assumptions are broken for some $\tilde{d}_{(j)}\gtrsim \tilde{d}_{\text{e}}$. By our definition of $\tilde{d}_{\text{e}}$, the second possibility of the failure is excluded, which means that either the first or the third possibilities happens during the iterations when $\tilde{d}_{(j)}\gtrsim \tilde{d}_{\text{e}}$. For any $t \ge \tilde{d}_{\text{e}}$, from the discussions above, we have 
    \begin{equation*}
        \mathbb{P}\left(\left\lVert \hat{\theta }_{\text{p}}(Y)-\theta _{0}\right\rVert_{2}^{}\ge t\right)\le \sum\limits _{j: \tilde{d}_{(j)}\gtrsim t}^{}\exp\left\{-\frac{-\tilde{d}_{(j)}^{2}}{c ^{4}}\right\}+\frac{1}{1+\tilde{d}_{(j)}^{3}}\lesssim \exp\left\{-\frac{t ^{2}}{c ^{4}}\right\}+\frac{1}{t ^{3}},
    \end{equation*}
    where some universal constant is omitted. Therefore, for any $\theta _{0}\in K$, we have 
    \begin{equation*}
        \begin{aligned}
            \mathbb{E}\left[\left\lVert \hat{\theta }_{\text{p}}(Y)-\theta _{0}\right\rVert_{2}^{2}\mathbf{1}_{\text{broken}}\right]&\le \displaystyle\int _{0}^{\infty }\mathbb{P}\left(\left\lVert \hat{\theta }_{\text{p}}(Y)-\theta _{0}\right\rVert_{2}^{}\ge t\right)2t dt\\ 
            &=\tilde{d}_{\text{e}}^{2}+2\displaystyle\int _{\tilde{d}_{\text{e}}}^{\infty }\mathbb{P}\left(\left\lVert \hat{\theta }_{\text{p}}(Y)-\theta _{0}\right\rVert_{2}^{}\ge t\right)tdt\\ 
            &\lesssim \tilde{d}_{\text{e}}^{2}+\frac{1}{\tilde{d}_{\text{e}}}+1.
        \end{aligned}
    \end{equation*}
    Since we assume that $\text{diam}(K) \ge 1$, we have $\epsilon ^{*}\gtrsim _{M}1$ by Lemma \ref{lemma: relationship between epsilon and diam} and consequently $\tilde{d}_{e}\gtrsim _{M}1$. We come to the performance for $\hat{\theta }$ that 
    \begin{equation*}
        \mathbb{E}\left[\left\lVert \hat{\theta }_{\text{p}}(Y)-\theta _{0}\right\rVert_{2}^{2}\mathbf{1}_{\text{broken}}\right]\lesssim \tilde{d}_{\text{e}}^{2}=\Upsilon (K)\epsilon ^{\star 2},
    \end{equation*}
    where $\Upsilon (K)=2 \sqrt{\tau _{K _{\text{m}}}T _{2}(K)}\vee c ^{2} \sqrt{\lceil 4 \lceil \log _{c}(\sqrt{\tau _{K _{\text{m}}}}T _{2}(K)) \rceil \rceil }$. Taking supreme over $\theta _{0}\in K$ and combining with \eqref{eq: matching the minimax rate under ideal assumptions}, \eqref{eq: main inequality for the estimator} is proved.

    Now it remains to prove the efficiency of Algorithm \ref{algorithm: main algorithm}. We first note that $\kappa _{\text{m}}$ is at most a poly-logarithmic factor in $n$, as shown by Theorems 3.1, 3.3 in \cite{neykov2026polynomialtimenearoptimalestimationcertain}. Also, the total amount of iterations of Algorithm \ref{algorithm: main algorithm} is no greater than $\log \left(\frac{R}{r}\right)$ while the time complexity of each iteration is at most polynomial in $n,\log \left(\frac{R}{r}\right)$ by Theorems \ref{theorem: main theorem of the cited paper}, \ref{theorem: polynomial time approximation to the SDP problem}. Therefore, the total time complexity of Algorithm \ref{algorithm: main algorithm} is at most polynomial in $n,\log \left(\frac{R}{r}\right)$. The proof is completed.
\end{proof}

\subsection{Proofs of Section \ref{section: example}}\label{subsection: proofs of section of example}

\subsubsection{Proof of Lemma \ref{lemma: critical radius for ellipsoid}}\label{subsubsection: proof of critical radius for ellipsoid}

\begin{proof}
    Consider $\tilde{\epsilon }=\sqrt{32c ^{2}k ^{*}/\kappa (M)}$, and take $P _{k ^{*}}^{*}$ as the optimal projection corresponding the Kolmogorov $k$-width. For the ellipsoid $K$ in the finite-dimensional space $\mathbb{R}^{n}$, it is not hard to verify that $P ^{*}_{k ^{*}}$ is merely the truncation operator on the last $k$ coordinates. Now for any $\beta \in K$, consider a maximal $\tilde{\epsilon }/c$ packing for $B _{2}(\beta ,\tilde{\epsilon })\cap K$. For any two points $x,y$ in this packing $L$, by definition we have $\left\lVert x-y\right\rVert_{2}^{}\ge \tilde{\epsilon }/c$. From the triangle inequality and the \hyperref[def_kolmogorov_width_raw]{definition} of the Kolmogorov $k$-width, we know 
    \begin{equation*}
        \begin{aligned}
            \tilde{\epsilon }/c &\le \left\lVert x-y\right\rVert_{2}^{}\\ 
            &\le \left\lVert x-P ^{*}_{k ^{*}}x\right\rVert_{2}^{}+\left\lVert y-P ^{*}_{k ^{*}}y\right\rVert_{2}^{}+\left\lVert P ^{*}_{k ^{*}}x-P ^{*}_{k ^{*}}y\right\rVert_{2}^{}\\ 
            &\le 2D _{k ^{*}}(K)+\left\lVert P _{k ^{*}}^{*}x-P _{k ^{*}}^{*}y\right\rVert_{2}^{}\\ 
            &=2 \sqrt{a _{n-k ^{*}}}+\left\lVert P _{k ^{*}}^{*}x-P _{k ^{*}}^{*}y\right\rVert_{2}^{}.
        \end{aligned}
    \end{equation*}
    Therefore, we have $\left\lVert P _{k ^{*}}^{*}x-P _{k ^{*}}^{*}y\right\rVert_{2}^{}\ge \tilde{\epsilon }/c-2 \sqrt{a _{n-k ^{*}}}\ge \tilde{\epsilon }/c-2 \sqrt{\frac{k ^{*}+1}{\kappa (M)}}\ge \frac{\tilde{\epsilon }}{2c}$. This means that $L$ induces a $\frac{\tilde{\epsilon }}{2c}$ packing set with the same cardinality. However, from the property of the projection, we know $\left\lVert P _{k ^{*}}^{*}x-P _{k ^{*}}^{*}\beta \right\rVert_{2}^{}\le \left\lVert x-\beta \right\rVert_{2}^{}\le \tilde{\epsilon }$, which implies that $\left\{P _{k ^{*}}^{*}x:x \in L\right\}$ is contained within a $k ^{*}$-dimensional ball $B _{2}(P _{k ^{*}}^{*}\beta ,\tilde{\epsilon })$. From the standard volumetric argument of the entropy, we know $\log N(\tilde{\epsilon }/c,B _{2}(\beta ,\tilde{\epsilon })\cap K)\le \log M(\tilde{\epsilon }/c,B _{2}(\beta ,\tilde{\epsilon })\cap K)\le \log M \left(\frac{\tilde{\epsilon }}{2c},B _{2}(P _{k ^{*}}^{*}\beta ,\tilde{\epsilon })\right)\asymp k ^{*}=\frac{\tilde{\epsilon }^{2}\kappa (M)}{32c ^{2}}$. Consequently taking supreme over $\beta \in K$, we come to $\epsilon ^{*}\lesssim _{M}\tilde{\epsilon }\asymp _{M}\sqrt{k ^{*}}$.

    On the other hand, the ellipsoid $K$ contains a $k ^{*}$-dimensions ball with radius $\sqrt{\frac{k ^{*}}{\kappa (M)}}$ since $a _{n-k+1}>k ^{*}/\kappa (M)$ for $k=1,2,\dots,k ^{*}$. Take $\tilde{\epsilon }=\sqrt{\frac{k ^{*}}{\kappa (M)}}$, we know that
    \begin{equation*}
        \log N ^{\text{loc}}(\tilde{\epsilon },c)\ge \log N(\tilde{\epsilon }/c,B _{2}(\mathbf{0},\tilde{\epsilon }))\asymp k ^{*}=\tilde{\epsilon }^{2}\kappa (M).
    \end{equation*}
    By the \hyperref[def_critical_radius]{definition} of $\epsilon ^{*}$, we know $\epsilon ^{*}\gtrsim \tilde{\epsilon }\asymp _{M}\sqrt{k ^{*}}$. The proof is completed.
\end{proof}

\subsubsection{Proof of Corollary \ref{corollary: minimax rate for constrained logistic regression}}\label{subsubsection: proof of minimax rate for constrained logistic regression}

\begin{proof}
    This corollary originates directly from Lemma \ref{lemma: critical radius for ellipsoid} and Theorem \ref{theorem: main theorem for theoretical algorithm}. We first note that $\text{diam}(K) ^{2}=4a _{n}$. When $k ^{*}>0$, the results are immediate from Lemma \ref{lemma: critical radius for ellipsoid} and the fact that $k ^{*}+1\asymp k ^{*}$. When $k ^{*}=0$, it implies that $a _{n}\le \frac{1}{\kappa (M)}\lesssim _{M}1$. Setting $\tilde{\epsilon }=\sqrt{a _{n}}$, we have $\tilde{\epsilon }^{2}\kappa (M)=a _{n}\kappa (M)\le 1 \lesssim _{M}\log N ^{\text{loc}}(\tilde{\epsilon },c)$ for sufficiently large $c$. Therefore we have $\epsilon ^{*}\gtrsim \sqrt{a _{n}}\asymp \text{diam}(K) $ and consequently we know from Theorem \ref{theorem: main theorem for theoretical algorithm} that the minimax rate of the $\ell _{2}$ estimation error is $\text{diam}(K) ^{2}\asymp _{M}a _{n}\lesssim _{M}(k ^{*}+1)\wedge a _{n}$ in this case. The proof is completed.
\end{proof}

\end{document}